\documentclass[11pt]{amsart}
\usepackage{amssymb}
\usepackage{amsmath}
\usepackage{amscd}
\usepackage{amsthm}
\usepackage{enumerate}
\usepackage{tikz-cd}
\usepackage{graphicx,color}
\usepackage{booktabs,array}
\usepackage[colorlinks=true,linkcolor=blue,citecolor=blue]{hyperref}
\usepackage[margin=3cm,marginpar=2.5cm]{geometry}

\usepackage{tikz}
\usetikzlibrary{shapes.geometric}
\usetikzlibrary{angles}
\usepackage{caption, subcaption}
\usetikzlibrary{calc,decorations.pathreplacing}
\usetikzlibrary{arrows}
\usetikzlibrary{intersections}

\tikzset{
  wfib/.style={circle,draw,fill=white,inner sep=0pt,minimum size=6pt},
  wcpt/.style={rectangle,draw,fill=white,inner sep=0pt,minimum size=6pt},
  blk/.style={circle,draw,fill=black,inner sep=0pt,minimum size=5pt},
  hor/.style={line width=1.4pt},
  nhor/.style={line width=0.5pt},
  hornotY/.style={line width=1.4pt,dash pattern=on 3pt off 2pt},
  cont/.style={line width=1.4pt,dotted},
  compbg/.style={line width=13pt,gray!22,line cap=round,line join=round},
  typ/.style={font=\scriptsize,inner sep=1.5pt},
}

\renewcommand{\bar}[1]{\overline{#1}}

\newcommand*{\colvec}[2]{\begin{pmatrix}#1\\#2\end{pmatrix}}
\newtheorem{thm}{Theorem}[section]
\numberwithin{equation}{section}
\usepackage{etoolbox}
\newcounter{eqanchor}

\AtBeginEnvironment{equation}{\stepcounter{eqanchor}}
\newtheorem{prop}[thm]{Proposition}
\newtheorem{lem}[thm]{Lemma}
\newtheorem{cor}[thm]{Corollary}

\theoremstyle{definition}
\newtheorem{defn}[thm]{Definition}
\newtheorem{rem}[thm]{Remark}
\newtheorem{conv}[thm]{Convention}
\newtheorem{exmp}[thm]{Example}
\newtheorem{ques}[thm]{Question}

\newtheorem{notation}[thm]{Notation}

\title[]{Subgroup distortion in free-by-cyclic groups with linearly growing monodromy} \date{}

\author{Hoang Thanh Nguyen}
\address{Hoang Thanh Nguyen, Department of Mathematics, FPT University, Da Nang, Vietnam}
\email{hoangnt63@fe.edu.vn}



\begin{document}

\begin{abstract}
Let $G=F\rtimes_\phi\mathbb Z$ be the mapping torus of a linearly growing UPG automorphism $\phi$ of a finitely generated free group $F$. We compute the distortion of every finitely generated subgroup $H$ of $G$, and we show that it is linear, quadratic or exponential. In contrast, for each $d\ge2$, Dison and Riley constructed an automorphism of polynomial degree $d$ whose mapping torus contains a free subgroup with distortion equivalent to the Ackermann function $A_{d+1}$, which grows faster than every iterated exponential function. Using the natural splitting of $G$ over $\mathbb Z^2$ subgroups and its topological representative, we associate to $H$ a family of \emph{horizontal components}, together with \emph{dilation homomorphisms} defined from slopes along these components. We prove that $H$ is undistorted in $G$ if and only if every horizontal component of $H$ consists of only one piece. Otherwise, the distortion of $H$ is quadratic if all dilation homomorphisms are trivial, and exponential if one of them is nontrivial. Along the way, we show that the distortion of $H$ is determined, up to superadditive closure, by the distortions of the stabilizers of its horizontal components.
\end{abstract}

 \maketitle

\section{Introduction}\label{sec:intro}

 A finitely generated group $G$ can be considered as a metric space when we equip the group with the word metric from a finite generating set. 
 Gromov has been successful in promoting this idea to study finitely generated
groups. Different finite generating sets of $G$ give different metrics on $G$. However, these metric spaces are unique up to quasi-isometry. Here, a map between two metric spaces is called a quasi-isometry if it is coarsely bi-Lipschitz and coarsely surjective.

A basic goal in geometric group theory, proposed by Gromov, is to quantify the failure of the coarsely bi-Lipschitz property of the inclusion $H \to G$ of
a finitely generated subgroup $H$ in a finitely generated group $G$.
 More precisely, let $\mathcal{A}$ and $\mathcal{S}$ be finite generating sets of $G$ and $H$, respectively. The \emph{distortion} of $H$ in $G$ is the function
\[
\Delta^G_H(n)=\max\bigl\{\,|h|_{\mathcal S}\ \big|\ h\in H,\ |h|_{\mathcal A}\le n\,\bigr\}.
\]
Up to the natural equivalence $\simeq$ of Section~\ref{subsec:distortion}, the function $\Delta_{H}^{G}$ does not depend on the choice of the finite generating sets $\mathcal{A}$ and $\mathcal{S}$.

The subgroup $H$ is called \textit{undistorted} in $G$ if $\Delta^G_H(n)\simeq n$, equivalently, if the inclusion $H\hookrightarrow G$ is a quasi-isometric embedding. Undistorted subgroups are those whose large-scale geometry is preserved inside $G$. Standard examples are quasiconvex subgroups of hyperbolic groups and subgroups acting cocompactly on convex subcomplexes of $\mathrm{CAT}(0)$ cube complexes on which the ambient group acts properly and cocompactly. 

\begin{ques}\label{ques:main}
Given a group $G$, which finitely generated subgroups of $G$ are undistorted, and how large is the distortion of those that are not?
\end{ques}

The set of distortion functions can be very rich, even in groups with good geometric properties. For instance, Dani and Riley \cite{DR25} recently showed that for all integers $p>q>0$ and $k>0$, every non-elementary torsion-free hyperbolic group $H$ is a subgroup of a hyperbolic group $G$ with $\Delta^G_H(n)\simeq\exp^k(n^{p/q})$, where $\exp^k$ denotes the $k$-fold iterated exponential function. Free-by-cyclic groups with polynomially growing monodromy also exhibit very large distortion. For example, the hydra groups of Dison and Riley \cite{DR13} contain free subgroups whose distortion grows like the Ackermann functions (see Section~\ref{subsec:beyond}). In contrast, we show that in the linear case only three distortion functions occur.

In this paper we study Question~\ref{ques:main} for free-by-cyclic groups $G=F\rtimes_\phi\mathbb Z$, where $F$ is a finitely generated free group and $\phi$ is a linearly growing UPG automorphism. We call such groups \textit{linear free-by-cyclic groups}. We obtain the following results.
\begin{itemize}
\item We compute the distortion of every finitely generated subgroup of $G$. It is linear, quadratic or exponential, and the three cases are distinguished by explicit geometric invariants, namely the horizontal components of the subgroup and their dilation homomorphisms (see Theorem~\ref{thm:intro_main}).

\item We show that, up to superadditive closure, the distortion of a finitely generated subgroup is determined by the distortions of finitely many explicit subgroups, namely the stabilizers of its ``horizontal components'' (see Theorem~\ref{thm:intro_reduction}).
\end{itemize}

By \cite{AM22,DT23}, in the form given by Hagen and Wise \cite[Proposition~4.11]{HW25}, the group $G$ splits as a graph of groups over a bipartite graph $\Delta$ (see Proposition~\ref{prop:split-over-tori}), with the following properties.
\begin{itemize}
\item The white vertex groups are of the form $F_w\times\mathbb Z$, with $F_w$ free of rank at least two.
\item The black vertex groups are copies of $\mathbb Z^2$.
\item Every edge group is the whole adjacent black vertex group.
\end{itemize}
We work with the topological representative $M_\varphi$ of this splitting constructed by Hagen and Wise \cite[Section~4.3.1]{HW25} (see Section~\ref{subsection:topological_representative}). It is a graph of spaces whose white vertex spaces are products $X_w\times S^1$ of a graph and a circle, and whose black vertex spaces are tori. The torus at a black vertex $b$ is glued to the adjacent white vertex spaces with pairwise distinct twisting numbers $k(e)$ (see Proposition~\ref{prop:split-over-tori}(\ref{prop:split-over-tori8})). We do not assume that the black vertices of $\Delta$ have valence two (when every black vertex has valence two, then such a mapping torus is called ``unbranched'' in \cite{BGGH25}).

\subsection{Horizontal components} 

The distortion of a finitely generated subgroup $H\le G$ is governed by two invariants which we introduce in this paper, the \textit{horizontal components} of $H$ and their \textit{dilation homomorphisms}. They are defined for every finitely generated subgroup of $G$.  They are modelled on the almost fiber part and the spirality of horizontal surfaces in graph manifolds \cite{RW98,Liu17,HN19,Ngu19}. In our setting, however, there is no surface, and the invariants have to be built from the subgroup itself.

The new difficulties are described in Section~\ref{subsec:comparison}. In Section~\ref{sec:examples} we show that the same two invariants recover known distortion phenomena, which were previously proved by different methods.

Let $H\le G$ be finitely generated. If $H$ is conjugate into a vertex group, then $H$ is undistorted in $G$ (see Remark~\ref{rem:whyundistorted}). So we assume that $H$ is not conjugate into a vertex group. Let $T$ be the Bass--Serre tree of the splitting of $G$, and let $T_H\subset T$ be the minimal $H$-invariant subtree. To define horizontal components, we need the following terminology.
\begin{itemize}
\item A white vertex $\tilde w$ of $T_H$ is of \textit{fiber type} if $H\cap Z(G_{\tilde w})$ is trivial, where $Z(G_{\tilde w})\cong\mathbb Z$ is the centre of $G_{\tilde w}$. Otherwise $\tilde w$ is of \textit{compact-fiber type} (see Definition~\ref{defn:fibertype}).
\item An edge $(\tilde b,\tilde w)$ of $T_H$ is \textit{horizontal} if $H\cap G_{\tilde b}$ is infinite cyclic and $\tilde w$ is of fiber type.
\item A black vertex is \textit{horizontal} if it has at least two horizontal edges (see Definition~\ref{defn:horizontal_black}).
\end{itemize}

\begin{defn}\label{defn:intro_horizontal}
The \textit{horizontal components} of $H$ are the connected components $T_i$ of the subgraph $\widehat Y\subset T_H$ spanned by the horizontal edges at horizontal black vertices (see Definitions~\ref{defn:hatYandOmegaTi} and~\ref{defn:horizontal_component}). A horizontal component $T_i$ \textit{consists of one piece} if at most one of its white vertices has two or more neighbours in $T_i$.
\end{defn}

Horizontal components are the analogues of the components of the almost fiber part of a surface.
In particular, every white vertex of a horizontal component is of fiber type.  We refer the reader to Figure~\ref{fig:hatY} in Section~\ref{sec:model} for an illustration of these notions.
To each horizontal component $T_i$ we attach two objects:
\begin{itemize}
\item its stabilizer $H_i=\mathrm{Stab}_H(T_i)$,
\item a \textit{dilation homomorphism} $\omega_i\colon H_i\to(\mathbb Q_{>0},\times)$ (see Definition~\ref{defn:dilation}).
\end{itemize}
The dilation homomorphism is built from slopes. Let $\tilde b$ be a horizontal black vertex and let $u$ generate $H\cap G_{\tilde b}$. Roughly speaking, the \textit{slope} of an ordered pair of horizontal edges at $\tilde b$ compares the horizontal coordinates of $u$ as seen from the two corresponding white vertices (see Definition~\ref{defn:slope}). For $h\in H_i$, the value $\omega_i(h)$ is the product of the slopes at the black vertices of the geodesic in $T_i$ from a white vertex $x$ to $hx$.

\subsection{Main results}
Our main theorem computes the distortion of every finitely generated subgroup of $G$.

\begin{thm} [Theorem~\ref{thm:main}] \label{thm:intro_main}
Let $\phi$ be a linearly growing UPG automorphism of a finitely generated free group $F$, let $G=F\rtimes_\phi\mathbb Z$, and let $H\le G$ be a finitely generated subgroup which is not conjugate into a vertex group. Then
\[
\Delta^G_H(n)\ \simeq\ \begin{cases}
n & \text{if every horizontal component consists of one piece},\\
n^2 & \text{if some horizontal component does not consist of one}\\
 & \text{piece, and every dilation homomorphism is trivial},\\
2^n & \text{if some dilation homomorphism is nontrivial}.
\end{cases}
\]
In particular, $H$ is undistorted in $G$ if and only if every horizontal component of $H$ consists of one piece.
\end{thm}

An immediate consequence of Theorem~\ref{thm:intro_main} is the following.

\begin{cor}[Corollary~\ref{cor:spectrum}]\label{cor:intro_spectrum}
Let $G$ be a linear free-by-cyclic group. The distortion function of every finitely generated subgroup of $G$ is equivalent to $n$, to $n^2$, or to $2^n$.
\end{cor}

\medskip
The proof of Theorem~\ref{thm:intro_main} has three parts. The first part is a characterization of the horizontal components which contribute to the distortion, together with lower bounds.

\begin{thm}[Lower bounds, Theorems~\ref{thm:lower_bound} and~\ref{thm:dichotomy} and Lemma~\ref{lem:one_piece}]\label{thm:intro_lower}
Let $G$ and $H$ be as in Theorem~\ref{thm:intro_main}, and let $T_i$ be a horizontal component of $H$. The following are equivalent:
\begin{enumerate}
\item $T_i$ consists of one piece,
\item $T_i$ is bounded,
\item $H_i$ fixes a vertex of $T_i$,
\item $H_i$ is undistorted in $G$.
\end{enumerate}
If these conditions fail, then
\[
\Delta^G_H(n)\ \succeq\ \Delta^G_{H_i}(n)\ \succeq\ n^2.
\]
\end{thm}

The second part is a reduction theorem, which is the analogue of \cite[Theorem~4.1]{Ngu19}.

\begin{thm}[Reduction to horizontal components, Theorem~\ref{thm:reduction}]\label{thm:intro_reduction}
Let $G$ and $H$ be as in Theorem~\ref{thm:intro_main}. Let $T_1,\dots,T_r$ be representatives of the $H$-orbits of horizontal components, and let
\[
f(n)=\max\{n,\Delta^G_{H_1}(n),\dots,\Delta^G_{H_r}(n)\}.
\]
Then
\[
f\ \preceq\ \Delta^G_H\ \preceq\ \overline f,
\]
where $\overline f$ is the superadditive closure of $f$.
\end{thm}

The third part establishes upper bounds for the stabilizers of the horizontal components. 

By Theorem~\ref{thm:intro_reduction} and an elementary estimate for superadditive closures (see Lemma~\ref{lem:closure_powers}), the upper bounds in Theorem~\ref{thm:intro_main} reduce to establishing the corresponding bounds for the stabilizers, namely $\Delta^G_{H_i}(n)\preceq n^2$ when $\omega_i$ is trivial and $\Delta^G_{H_i}(n)\preceq 2^n$ in general (see Theorem~\ref{thm:upper_Hi}).

The proof follows the \textit{height} of a point above the model space along a path in $\widetilde M_\varphi$. Within a white vertex space, the height varies in a Lipschitz manner. When the path crosses a black vertex space, the height is multiplied by the inverse of a slope, up to a bounded error (see Lemma~\ref{lem:height_black}). It follows that along a path of length $n$, the heights are bounded by a linear function of $n$ when the slope products are bounded, which is the case when $\omega_i$ is trivial. In general, they are bounded by an exponential function of $n$. This argument is in the spirit of the upper bound in \cite[Section~6]{HN19}, where the dilation of a horizontal surface is interpreted in terms of slopes of lines in the JSJ planes.

\subsection{\texorpdfstring{Relation with $3$-manifolds and comparison with \cite{HN19} and \cite{Ngu19}}{Relation with 3-manifolds and comparison with HN19 and Ngu19}}\label{subsec:comparison}

When every black vertex of $\Delta$ has valence $2$ (the \textit{unbranched} case of \cite{HW25}, \cite{BGGH25}), the splitting of $G$ resembles the JSJ decomposition of a graph manifold. The white vertex spaces play the role of Seifert fibered blocks, and the black tori play the role of JSJ tori. Distortion of subgroups of $3$-manifold groups is well understood.
\begin{itemize}
\item Hruska and Nguyen \cite{HN19} proved that a horizontal surface in a graph manifold is quadratically or exponentially distorted, according to whether its dilation (or spirality) invariant introduced in \cite{RW98,Liu17} is trivial or not.
\item Nguyen \cite{Ngu19} showed that the distortion of a $\pi_1$-injective, properly immersed surface in a nongeometric $3$-manifold is determined by the components of its almost fiber part. In particular, the surface is undistorted if and only if each such component consists of only one piece.
\item Nguyen and Sun \cite{NS20} studied the distortion of arbitrary finitely generated subgroups of $3$-manifold groups.
\end{itemize}

A linear free-by-cyclic group $G$ is in general not a $3$-manifold group, and our setting differs from that of surfaces in graph manifolds in several ways, which we now describe.

Our overall strategy follows that of \cite{HN19,Ngu19,NS20}. The approaches in these works build respectively on \cite{RW98,Liu17,Sun20}. In our setting, however, we have to construct the horizontal components and the dilation homomorphisms from scratch. We expect these constructions to be useful in other contexts as well. For example, we believe that horizontal components introduced in this paper can be used to characterize which subgroups of linear free-by-cyclic groups are separable as in \cite{Liu17}, \cite{Sun20}.

However, several parts of arguments of \cite{HN19}, \cite{Ngu19}, \cite{NS20} rely essentially on the surface structure and require different approaches in our setting. The main new ingredients are as follows.

\begin{enumerate}
\item \textit{A model space for an arbitrary subgroup.} In \cite{HN19,Ngu19}, the subgroup is the fundamental group of an immersed surface, whose universal cover provides a natural geometric model for the group action. An arbitrary finitely generated subgroup $H\le G$ comes with no such model. We therefore construct a model space $\Omega_H\subset\widetilde M_\varphi$ directly (see Section~\ref{sec:model}). In a white vertex space of fiber type, the corresponding piece of $\Omega_H$ is a band around the graph of an equivariant Lipschitz function defined on a subtree obtained from a Stallings core (see Section~\ref{subsec:horizontalpiece_white} and Definition~\ref{defn:band}). These pieces are joined by strips over the black pieces.

\item \textit{Subgroups which are not surface subgroups.} A $\pi_1$-injective surface meets the JSJ tori in essential curves. In contrast, a finitely generated subgroup $H\le G$ may meet a black vertex group trivially, or in a subgroup commensurable with the centre of an adjacent white vertex group. This gives rise to the four types of black vertices described in Lemma~\ref{lem:black_types}. When $H\cap G_{\tilde b}=1$, the gate estimate of \cite{Ngu19} does not apply. We establish the corresponding estimate using a new argument based on the pairwise distinctness of the twisting numbers at a black vertex (see Lemma~\ref{lem:local_gates}(\ref{local_gates_2})).

\item \textit{Branching.} A JSJ torus has exactly two sides, whereas a black vertex of the underlying graph $\Delta$ may have arbitrary valence. Thus, slopes are associated to ordered pairs of edges at a black vertex and satisfy a cocycle relation (see Lemma~\ref{lem:slope}(\ref{lem:slope_4})). Moreover, the dilation homomorphism is well defined only after establishing a median property in $T_i$ (see Lemma~\ref{lem:basisPi}).

\item \textit{Graphs instead of surfaces.} The white vertex spaces are products of graphs and circles, and the tori are attached along immersed circuits. These circuits are not boundary components and need not be disjoint. For the lower bound, we replace the surface geometry used in \cite{HN19} with arguments involving trees. The main ingredients are a bounded overlap property for the axes of two elements whose cyclic subgroups intersect trivially (see Lemma~\ref{lem:axisintersection}) and a gate property for pairs of lines in a tree (see Lemma~\ref{lem:boundedgates}). For the upper bound, the gates in the tree factors arise from the maximality of the cyclic edge groups and the fact that two distinct black neighbours of a white vertex have stabilizers with cyclic intersection (see Lemma~\ref{lem:tree_gates}).
\end{enumerate}

\subsection{Beyond linear UPG automorphisms}\label{subsec:beyond}
The UPG assumption in Corollary~\ref{cor:intro_spectrum} causes no loss of generality. Suppose that $[\phi]$ is linearly growing. By \cite{BFH00}, some positive power $[\phi]^m$, with $m\ge1$, is UPG and remains linearly growing. The subgroup $G'=F\rtimes\langle t^m\rangle$ has index $m$ in $G=F\rtimes_\phi\mathbb Z$ and is isomorphic to the mapping torus of a UPG representative of $[\phi]^m$ (see Section~\ref{subsec:free-by-cyclic}). For any finitely generated subgroup $H\le G$, the subgroup $H'=H\cap G'$ has finite index in $H$. Moreover, $\Delta^G_H\simeq\Delta^{G'}_{H'}$, since the inclusions $H'\hookrightarrow H$ and $G'\hookrightarrow G$ are quasi-isometries. It follows that whenever $[\phi]$ is linearly growing, the distortion of every finitely generated subgroup of $F\rtimes_\phi\mathbb Z$ is equivalent to one of $n$, $n^2$ or $2^n$.

The situation is quite different for automorphisms of higher polynomial degree. For $k\ge2$, Dison and Riley \cite{DR13} consider the group
\[
G_k=\bigl\langle a_1,\dots,a_k,t\ \bigm|\ t^{-1}a_1t=a_1,\ t^{-1}a_it=a_ia_{i-1}\ (1<i\le k)\bigr\rangle ,
\]
which is the mapping torus of a UPG automorphism of $F(a_1,\dots,a_k)$ with polynomial growth of degree $k-1$. They  consider its free subgroup $H_k=\langle a_1t,\dots,a_kt\rangle$ of rank $k$ and prove that $\Delta^{G_k}_{H_k}\simeq A_k$, where $A_k$ is the $k$-th Ackermann function \cite[Theorem~1.3]{DR13}. Thus, for an automorphism of degree $d\ge2$, a finitely generated subgroup can have distortion equivalent to $A_{d+1}$, which is not bounded by any fixed iterate of the exponential function.

At present we do not know the answer to Question~\ref{ques:main} when $\phi$ is polynomially growing of degree at least $2$. However, we expect the linear case treated in this paper to serve as the base step of an induction on the degree. This expectation is based on the cyclic hierarchy described in \cite[Proposition~4.3]{BGGH25}.

\subsection{Examples}\label{subsec:intro_examples}
Let $G = F \rtimes_{\phi} \mathbb Z$ with $\phi$ polynomially growing of degree $d \ge 0$. Then $\Delta_{F}^{G}(n) \simeq n^{d+1}$ by combining \cite{Ger94} and \cite{Mac02}.
The exponential distortion of the hydra subgroup is due to Dison and Riley \cite{DR13}. These results were obtained by unrelated methods.

In Section~\ref{sec:examples} we recover the second result, and the case $d=1$ of the first one, from our invariants. For subgroups of the fiber, each dilation homomorphism is trivial and hence the distortion is at most quadratic (see Corollary~\ref{cor:intro_fiber}). For the hydra subgroup $H_2$ there is a single horizontal component, and its dilation homomorphism takes the value $2$ (see Example~\ref{ex:hydra}). Hence the distortion is exponential by our main result.

\subsection{Outline} Section~\ref{sec:preliminary} contains background on free-by-cyclic groups, Bass--Serre trees and subgroup distortion, together with the splitting of linear free-by-cyclic groups and its topological representative from \cite{HW25}. In Section~\ref{sec:model}, we construct the model space and the horizontal components of a subgroup, and we show that the stabilizers of the horizontal components are undistorted in the subgroup. In Section~\ref{sec:slopes}, we introduce slopes and dilations and prove the basic properties that are used later, in particular the key identity. Sections~\ref{sec:labelling}, \ref{sec:rho} and~\ref{sec:overlap} are devoted to the lower bounds. We construct a family of elements of $H_i$ whose length in $G$ is at most linear, and we prove Theorems~\ref{thm:lower_bound} and~\ref{thm:dichotomy}. In Section~\ref{sec:undistortion}, we prove the reduction theorem (see Theorem~\ref{thm:reduction}). In Section~\ref{sec:upper}, we prove the upper bounds for the stabilizers (see Theorem~\ref{thm:upper_Hi}) and the main theorem (see Theorem~\ref{thm:main}). Finally, in Section~\ref{sec:examples}, we recover known examples from our invariants.

\subsection*{AI declaration} The author used Claude (Opus 5.5) to assist with checking and streamlining the calculations, proofreading and improving the English, organizing the paper, and generating the figures. Claude also suggested Examples 10.1 and 10.2. The author takes full responsibility for the mathematical arguments and the final content of the paper.

\section{Preliminaries}
\label{sec:preliminary}

\subsection{Free-by-cyclic groups}\label{subsec:free-by-cyclic}

Let $F$ be a free group of finite rank and let $\phi\in\mathrm{Aut}(F)$.
The \textit{mapping torus} of $\phi$ is the semidirect product
\[
G=F\rtimes_\phi\mathbb Z=\bigl\langle F,\,t \bigm| txt^{-1}=\phi(x)\ \text{for all } x\in F\bigr\rangle .
\]
Groups of this form are called \textit{free-by-cyclic}. Let $q\colon G\to\mathbb Z$
be the homomorphism with $q(F)=0$, $q(t)=1$ and its kernel is $F$. If $g\in F$
and $c_g$ denotes conjugation by $g$, then $t\mapsto gt$ induces an isomorphism
$F\rtimes_{c_g\circ\phi}\mathbb Z\to F\rtimes_\phi\mathbb Z$. Hence the isomorphism
type of $G$ depends only on the outer class $[\phi]\in\mathrm{Out}(F)$.
Topologically, if $f\colon X\to X$ is a homotopy equivalence of a finite graph
inducing $\phi$ on $\pi_1X\cong F$, then the mapping torus
$X\times[0,1]/(x,1)\sim(f(x),0)$ has fundamental group isomorphic to $G$, for a
suitable choice of the stable letter $t$.

Fix a basis of $F$ and for $g\in F$, let $\|g\|$ be the length of a cyclically
reduced word representing the conjugacy class of $g$. The outer automorphism
$[\phi]$ is \textit{polynomially growing} if for every $g\in F$ there is a
polynomial $P_g$ with $\|\phi^n(g)\|\le P_g(n)$ for all $n\ge1$. It is
\textit{linearly growing} if, in addition, $\|\phi^n(g)\|\le C_g\,n$ for every
$g\in F$ and some constant $C_g$, and $\|\phi^n(g)\|$ is unbounded for some
$g\in F$. These notions do not depend on the basis or on the representative of
$[\phi]$. Following \cite{BFH00}, a polynomially growing $[\phi]$ is
\textit{UPG} (unipotent polynomially growing) if the induced automorphism of
$H_1(F;\mathbb Z)$ is unipotent. Every polynomially growing outer automorphism
has a positive power that is UPG \cite{BFH00}. We say that $G=F\rtimes_\phi\mathbb Z$
is a \textit{linear free-by-cyclic group} if $[\phi]$ is linearly growing and UPG.

\subsection{Graphs of groups and Bass--Serre trees}\label{subsec:bass-serre}

A graph $\Delta$ has vertex set $\mathrm{Vert}(\Delta)$ and a set
$\mathrm{Edge}(\Delta)$ of oriented edges. An edge $e$ goes from $e^-$ to $e^+$.
A graph of groups $\mathcal G$ over $\Delta$ assigns groups $G_v$ and $G_e$ to
vertices and edges, together with monomorphisms $G_e\to G_{e^\pm}$. Its
fundamental group $\pi_1\mathcal G$ acts without inversions on the Bass--Serre
tree $T$ with quotient $\Delta$ where vertex and edge stabilizers are the conjugates
of the $G_v$ and $G_e$ \cite{Ser80}.

Let $H$ be a group acting on a tree $T$ without inversions. An element $h\in H$ is \textit{elliptic} if it fixes a vertex, and \textit{loxodromic} otherwise. A loxodromic element $h$ has a unique \textit{axis} $A_h\subset T$, which is a bi-infinite geodesic on which $h$ acts as a translation of length $\tau(h)>0$. If $H$ is finitely generated and every element of $H$ is elliptic, then $H$ fixes a vertex of $T$ \cite[Corollary~7.3]{Bas93}. If $H$ fixes no vertex of $T$, then the union $T_H$ of the axes of the loxodromic elements of $H$ is the unique minimal $H$-invariant subtree \cite[Proposition~7.5]{Bas93}. The action of $H$ on $T_H$ is minimal, so if $H$ is finitely generated, then $H\backslash T_H$ is a finite graph \cite[Proposition~7.9(b)]{Bas93}. By Bass--Serre theory, $H$ is then the fundamental group of the graph of groups over $H\backslash T_H$ with vertex groups $H_{\tilde v}=\mathrm{Stab}_H(\tilde v)$ and edge groups $H_{\tilde e}=\mathrm{Stab}_H(\tilde e)$.

A \textit{graph of spaces} over $\Delta$ assigns spaces $M^v$ and $M^e$ to
vertices and edges, and is obtained from $\bigsqcup_v M^v$ by attaching
$M^e\times[-1,1]$ along maps $M^e\times\{\pm1\}\to M^{e^\pm}$ that are
$\pi_1$-injective. Its fundamental group is the fundamental group of the
associated graph of groups \cite{SW79}, and its universal cover is a tree of
spaces over the Bass--Serre tree.

 \subsection{Linear mapping torus decomposition}
  \label{subsection:lineartopological}

  In the linear case, we employ a splitting over $\mathbb{Z}^2$ subgroups constructed in \cite[Proposition~5.2.2]{AM22} and \cite[Theorem~2.7]{DT23}. We adopt the formulation of these results from \cite[Proposition~4.11]{HW25}, where various properties from \cite{AM22}, \cite{DT23} are summarized. We only list the properties of \cite[Proposition~4.11]{HW25} that we use in this paper.

\begin{prop} [Linear mapping torus decomposition]
    \label{prop:split-over-tori}
Let $\phi$ be a linearly growing UPG  automorphism of the finite rank free group $F$ and let $G=F\rtimes_\phi\langle t\rangle$ be its mapping torus.  Then there is a finite connected graph $\Delta$ such that $F$ and $G$ split as graphs of groups with underlying graph $\Delta$, and the following all hold.
\begin{enumerate}
 \item \label{prop:split-over-tori1} $\Delta$ is bipartite, and $v\in \mathrm{Vert}(\Delta)$ is \textbf{black} if $F_v$ has rank $1$ and \textbf{white} otherwise, all edges are oriented from their black to their white vertex, and edge monomorphisms to black vertices are surjective. If $v$ is black, then $F_v=\langle r_v\rangle$ for some infinite-order $r_v$ and $r_e$ is sent by the edge map to $r_v$ for each $e$ with $e^-=v$.\label{item:linear-bipartite}

\item \label{prop:split-over-tori2} For each $v\in \mathrm{Vert}(\Delta)$, the vertex group $F_v$ has finite rank at least $1$, and $G_v=F_v\times \langle t_v\rangle$, where $t_v\in Ft$.

\item \label{prop:split-over-tori4} For each $e\in\mathrm{Edge}(\Delta)$, the group $F_e=\langle r_e\rangle$ is a maximal infinite cyclic subgroup of $F$ and $G_e=\langle r_e\rangle \times\langle t_e\rangle$.  Let $r_e$ denote the image of $r_e$ in $F_{e^-}$ and let $p_e$ denote its image in $F_{e^+}$.  Then the image of $t_e$ in $G_{e^-}$ is $t_{e^-}$ and the image of $t_e$ in $G_{e^+}$ is $p_e^{k(e)}t_{e^+}$ for some $k(e)\in\mathbb{Z}$.

\item If $e\in\mathrm{Edge}(\Delta)$ and $v\in \mathrm{Vert}(\Delta)$ is incident to $e$, then the centre of $G_v$ is contained in $\mathrm{image}(G_e\to G_v)$. \label{prop:split-over-tori5} 

\item Let $b\in \mathrm{Vert}(\Delta)$ be black and let $v,v'\in \mathrm{Vert}(\Delta)$ be white vertices joined to $b$ by distinct edges $e,e'$.  Then $\langle t_{b}r_e^{-k(e)}, t_{b}r_{e'}^{-k(e')}\rangle$ has finite index in $G_b$, so $k(e)\neq k(e')$.\label{prop:split-over-tori8}

\item If $\tilde v\in \mathrm{Vert}(T)$ is white and $\tilde b,\tilde b'\in \mathrm{Vert}(T)$ are distinct black vertices adjacent to $\tilde v$, then $\mathrm{Stab}_G(\tilde b)\cap\mathrm{Stab}_G(\tilde b')$ is cyclic.\label{prop:split-over-tori9}
\end{enumerate}
\end{prop}

 \subsection{Linear-growth topological representatives} \label{subsection:topological_representative}
Let $G$ be a linear mapping torus and let $\mathcal{G}$ be its graph of groups decomposition given by Proposition~\ref{prop:split-over-tori}.
We recall from \cite[Section~4.3.1]{HW25} a graph of spaces $M_\varphi$ with $\pi_1(M_\varphi)\cong G$, which we use throughout the paper. It has the following properties.

\begin{enumerate}
    \item The underlying graph of $M_\varphi$ is $\Delta$ that is bipartite.  Let $\mathrm{Vert}_\bullet(\Delta)$ be the set of black vertices and $\mathrm{Vert}_\circ(\Delta)$ the white vertices.  No element of $\mathrm{Vert}_\bullet(\Delta)$ is a leaf of $\Delta$.\label{item:multitwist-verts}

    \item For each $b\in \mathrm{Vert}_\bullet(\Delta)$, let $X_b$ be a graph homeomorphic to $S^1$, with one vertex $x_b$.  We identify $F_b$ with $\pi_1(X_b,x_b)$. \label{item:multitwist-vertspace-black}

    \item For each $w\in\mathrm{Vert}_\circ(\Delta)$, let $X_w$ be a finite graph with a base vertex $x_w$.  We identify $F_w$ with $\pi_1(X_w,x_w)$. \label{item:multitwist-vertspace-white}

    \item For each $e\in \mathrm{Edge}(\Delta)$, we recall that $b: = e^-\in\mathrm{Vert}_\bullet(\Delta)$ and $w: = e^+\in\mathrm{Vert}_\circ(\Delta)$. Let 
    $$
    \iota_e^{-} \colon X_b \times\{-1\}\to X_{b}
    $$ be the identity and let 
    $$
    \iota_{e}^{+} \colon X_{b}\times\{+1\}\to X_{w}
    $$ be obtained by subdividing $X_{b}\times\{+1\}$ so it is a combinatorial circle, and then mapping to the graph $X_w$ by a combinatorial immersion.  
    
    We let $p_e\to X_{w}$ be the immersed circuit which represents a conjugacy class of maximal cyclic subgroups in $F_w$ and has the property that $\iota_{e}^+$ corresponds to $p_e$. These maps induce the edge monomorphisms from Proposition \ref{prop:split-over-tori}. 
    We let $r_e\to X_b$ be the path corresponding to $\iota_e^-$, i.e. the identity map on an oriented $1$--cycle. \label{item:multitwist-edges}

\item The \textit{total space $Y$} is formed from $$X:=\bigsqcup_{v\in \mathrm{Vert}(\Delta)}X_v$$ by attaching, for each edge $e$, the cylinder $X_{e^-}\times[-1,1]$ using the given attaching maps $\iota_e^\pm$.  This makes $Y$ a $2$--complex whose $1$--skeleton consists of $X$, together with a copy of $\mathrm{Edge}(\Delta)$.  Moreover, $\pi_1Y\cong F$. 

\item Let $\Phi:=[\phi]\in\mathrm{Out}(F)$. There is a homotopy equivalence $\varphi\colon Y\to Y$, preserving the decomposition of $Y$ into the spaces $X_v$ and the cylinders, which induces on $F\cong\pi_1Y$ a Dehn twist (multitwist) automorphism representing $\Phi$ \cite[Section~4.3.1]{HW25}. We use the following graph of spaces model $M_\varphi$ for the mapping torus of $\varphi$, whose fundamental group is $G$ \cite[Section~4.3.1]{HW25}.

\item For each $v\in\mathrm{Vert}(\Delta)$, let $S_v$ be a $1$--cycle and let $M^v=X_v\times S_v$  (So, $M^v$ is a torus when $v\in\mathrm{Vert}_\bullet(\Delta)$.)  Then $M_\varphi$  splits as a graph of spaces with underlying graph $\Delta$ as follows.
For each $v\in\mathrm{Vert}(\Delta)$, the vertex space is $M^v$.
        For each $e = (b, w)\in \mathrm{Edge}(\Delta)$ oriented from a black vertex $b$ to a white vertex $w$, the edge space $M^e$ is $X_{b}\times S_{b}$. 
         We attach $M^e\times[-1,1]$ by
  identifying $M^e\times \{-1\}$ with $X_{b}\times S_{b}$ using the identity map, and gluing $M^e\times \{+1\}\cong X_{b}\times S_{b}$ to $X_{w}\times S_{w}$ using the map 
  \[
  I_e \colon M^{e} \times \{+1\} \to X_{w} \times S_{w}
  \]defined as follows.  The torus $M^e=X_{b}\times S_{b}$ is spanned by the loops $r_e$ (traverse $X_{b}$ once in the positive direction) and $t_e$ (traverse $S_{b}$ once in the positive direction). They represent the generators $r_b$ and $t_b$ of $G_b$. We identify $X_b\times S_b$ with $S^1\times S^1=(\mathbb R/\mathbb Z)^2$ by the product of the arc-length parametrizations of the two circles of length $1$, so that $r_e \mapsto \colvec{1}{0}$ and $t_e \mapsto \colvec{0}{1}$, where the first factor of $S^1\times S^1$ corresponds to $X_b$ and the second to $S_b$.  Now map $S^1\times S^1 \to X_{w}\times S_{w}$ by first applying the homeomorphism $S^1\times S^1\to S^1\times S^1$ given by $\begin{pmatrix}
            1& k(e)\\
            0&1
        \end{pmatrix}$
        and then the map $p_e\times \mathrm{id}$, where $p_e$ maps the first factor onto the cycle $p_e\subset X_w$. Thus $r_e\mapsto\colvec{1}{0}\mapsto p_e$ and $t_e\mapsto\colvec{k(e)}{1}\mapsto p_e^{k(e)}t_w$, in agreement with Proposition~\ref{prop:split-over-tori}(\ref{prop:split-over-tori4}).\label{item:multitwist-decomp}
\end{enumerate}

\subsection{Distortion}
\label{subsec:distortion}
For non-decreasing functions $f,g\colon[0,\infty)\to[0,\infty)$ we write $f\preceq g$ if there is a constant $C\ge1$ such that $f(x)\le C\,g(Cx+C)+Cx+C$ for all $x\ge0$, and we write $f\simeq g$ if $f\preceq g$ and $g\preceq f$. For example, $x^2+\Omega^{x}\simeq 2^x$ for every $\Omega>1$.
Let $G$ be a finitely generated group with a finite generating set $S$, and let $H\le G$ be a subgroup with a finite generating set $S'$. We write $|\cdot|_G=|\cdot|_S$ and $|\cdot|_H=|\cdot|_{S'}$ for the word lengths. The \textit{distortion} of $H$ in $G$ is the function
\[
\Delta_H^G(n):=\max\bigl\{\,|h|_H:\ h\in H,\ |h|_G\le n\,\bigr\}.
\]
Up to $\simeq$, $\Delta_H^G$ does not depend on the choice of the finite generating sets. The subgroup $H$ is called \textit{undistorted} in $G$ if $\Delta_H^G(n)\simeq n$.

We shall also apply the Milnor--\v{S}varc lemma to spaces which we do not know to be proper, and to actions which we only know to be co-bounded. The usual proof (see \cite[Proposition~I.8.19]{BH99} and \cite{BDM07}) gives the following version.
 
\begin{lem}[Milnor--\v{S}varc, co-bounded version]\label{lem:svarc-milnor}
Let $\Gamma$ be a group acting by isometries on a length space $(X,d)$, and let $o\in X$. Suppose that the following holds.
\begin{enumerate}
\item\label{lem:svarc-milnor_1} there is $R>0$ with $\Gamma\cdot B(o,R)=X$ (i.e., the action is \textit{co-bounded}),
\item\label{lem:svarc-milnor_2} for every $r\ge0$ the set $\{g\in\Gamma:\ d(o,go)\le r\}$ is finite.
\end{enumerate}
Then $S:=\{g\in\Gamma:\ d(o,go)\le 2R+1\}$ is a finite generating set of $\Gamma$, and
\[
d(o,go)\le(2R+1)\,|g|_S\qquad\text{and}\qquad |g|_S\le d(o,go)+2\qquad\text{for all }g\in\Gamma .
\]
In particular, $\Gamma$ is finitely generated and the orbit map $g\mapsto go$ is a quasi-isometry from $\Gamma$, with any word metric with respect to a finite generating set, to $(X,d)$.
\end{lem}

\begin{proof}
The proof of \cite[Proposition~I.8.19]{BH99} applies verbatim. Properness of $X$ is not needed, and condition (\ref{lem:svarc-milnor_2}) is used only to show that $S$ is finite. The bound $|g|_S\le d(o,go)+2$ comes from subdividing a path from $o$ to $go$ of length at most $d(o,go)+1$ into at most $d(o,go)+2$ subpaths of length at most $1$.
\end{proof}

\subsection{Choice of metrics}
\label{subsec:metrics}
Since we compute the distortion of a subgroup $H$ in $G=\pi_1(M_\varphi)$ by using the geometry of the universal cover $\widetilde M_\varphi$ (see Lemma~\ref{lem:distortionuniversal}), we need to discuss the metric on $M_\varphi$. We note that the choice of length metric does not affect the distortion, so we will choose a convenient metric.

\subsubsection{Metric on the topological representative}
Recall from Section~\ref{subsection:lineartopological} that for each vertex $v \in \mathrm{Vert}(\Delta) = \mathrm{Vert}_{\bullet}(\Delta) \cup \mathrm{Vert}_{\circ} (\Delta)$, the vertex space is $M^{v} = X_v \times S_v$ where $X_v$ is a graph and $S_v$ is a circle. We equip $X_v$ with the graph metric in which every edge has length $1$, and we equip $S_v$ with the metric of a circle of length $1$. We then equip $M^v$ with the $\ell^1$-product metric
\[
d_v\bigl((x,s),(x',s')\bigr):=d_{X_v}(x,x')+d_{S_v}(s,s').
\]
For each edge $e =(b, w)$ oriented from a black vertex $b$ to a white vertex $w$, recall that the edge space is $M^e = X_b \times S_b$.

We equip $M^e$ with the product metric and $M^{e} \times [-1,1]$ with the product metric as well.
With respect to these metrics, the gluing map $I_e \colon M^{e} \times \{+1\} \to X_w \times S_w$ (see Section~\ref{subsection:topological_representative}) is Lipschitz, since it is the composition of a linear shear and of $p_e\times\mathrm{id}$, where $p_e$ is a combinatorial immersion of a circle of length $1$ subdivided into $|p_e|$ edges. It need not be injective, because $p_e$ is only an immersed circuit. Its lifts to universal covers, however, are embeddings.

Since the underlying graph $\Delta$ of the splitting of $G$ is finite, there is a constant $\Lambda\ge1$ such that every $I_e$ is $\Lambda$-Lipschitz. The product metrics on $M^v$ and on $M^e\times[-1,1]$ induce a length metric on the space $M_\varphi$, denoted by $d$ (see \cite[Section~3.1]{BBI01}). The metric $d$ on $M_\varphi$ induces a metric on the universal cover $\widetilde{M}_\varphi$, which is also denoted by $d$ (by abuse of notation).

Let $\widetilde\pi\colon\widetilde M_\varphi\to T$ be the map of (\ref{D1}) below. For a vertex $\tilde v$ of $T$ over $v$, the lifted product metric on $\widetilde M^{\tilde v}=\widetilde\pi^{-1}(\tilde v)$ is denoted by $d_{\tilde v}$. Thus $d_{\tilde v}$ is the $\ell^1$-product metric on $\widetilde M^{\tilde w}=\widetilde X_w\times\mathbb R$ when $\tilde v=\tilde w$ is white, and on $\widetilde M^{\tilde b}=\widetilde X_b\times\mathbb R\cong\mathbb R^2$ when $\tilde v=\tilde b$ is black. By Lemma~\ref{lem:K_universal} below, there is a constant $K\ge1$ such that for every vertex $\tilde v$ of $T$ and all $x, y \in \widetilde{M}^{\tilde v}$ we have
\begin{equation}\label{eq:K}\tag{$\ast$}
\frac{1}{K} d_{\tilde v}(x, y) \le d(x, y) \le K d_{\tilde v}(x,y).
\end{equation}

For a path $\gamma$ in $\widetilde M_\varphi$ we write $\ell_d(\gamma)$ for its length with respect to $d$, and for a path $\gamma$ in $\widetilde M^{\tilde v}$ we write $\ell_{d_{\tilde v}}(\gamma)$ for its length with respect to $d_{\tilde v}$. Since lengths are computed from distances, (\ref{eq:K}) gives
\begin{equation}\label{eq:Klength}\tag{$\ast\ast$}
\frac1K\,\ell_{d_{\tilde v}}(\gamma)\le\ell_d(\gamma)\le K\,\ell_{d_{\tilde v}}(\gamma)\qquad\text{for every path $\gamma$ in $\widetilde M^{\tilde v}$.}
\end{equation}
Since $d$ is the length metric induced by the metrics on the pieces, every path contained in a single piece (a vertex space, or a cylinder $\widetilde M^{\tilde b}\times[-1,1]$) has $d$-length at most its length in the metric of that piece. In particular, for every edge $\tilde e=(\tilde b,\tilde w)$ of $T$ and every $z\in\widetilde M^{\tilde b}$, the segment $\{z\}\times[-1,1]$ of the cylinder of $\tilde e$ has length $2$, so
\begin{equation}\label{eq:vertical-segment}\tag{$\dagger$}
d\bigl(z,\iota_{\tilde e}(z)\bigr)\le 2\le 2K .
\end{equation}
Moreover, for the $\ell^1$-product metric on a product $X\times Y$ of length spaces, the length of a path $(\sigma,s)$ is $\ell(\sigma)+\ell(s)$.

The following result is a direct application of the Milnor--\v{S}varc lemma \cite[Proposition~I.8.19]{BH99}.

\begin{lem}\label{lem:distortionuniversal}
Let $x_0\in\widetilde M_\varphi$. There are constants $\lambda_G\ge1$ and $\epsilon_G \ge0$ such that
\[
\frac1{\lambda_G}|g|_G-\epsilon_G \le d(x_0,gx_0) \le \lambda_G \,|g|_G+\epsilon_G \quad\text{for all }g\in G .
\]
\end{lem}

The following lemma compares $d$ with the metrics $d_{\tilde v}$ on the vertex spaces.

\begin{lem}\label{lem:K_universal}
There is a constant $K\ge\Lambda$ such that for every vertex $\tilde v$ of $T$ and all $x,y\in\widetilde M^{\tilde v}$,
\[
\frac1{K}\,d_{\tilde v}(x,y) \le\ d(x,y) \le\ d_{\tilde v}(x,y).
\]
\end{lem}

\begin{proof}
The right-hand inequality holds because $d$ is the length metric obtained by gluing the vertex spaces and the cylinders. Since $\Delta$ is finite, it suffices to treat one vertex in each $G$-orbit. If $d(x,y)<1$, then, since $\widetilde\pi$ is $1$-Lipschitz and the edges of $T$ have length $2$ (see (\ref{D1})), paths from $x$ to $y$ of length less than $1$ stay in $\widetilde M^{\tilde v}$ and the adjacent half-cylinders. Retracting these half-cylinders onto $\widetilde M^{\tilde v}$, which is $\Lambda$-Lipschitz since $\iota_{\tilde e}$ is, gives $d_{\tilde v}(x,y)\le\Lambda\,d(x,y)$. Here the retraction is $(z,s)\mapsto z$ if $\tilde v$ is black and $(z,s)\mapsto\iota_{\tilde e}(z)$ if $\tilde v$ is white. If $d(x,y)\ge1$, then, since $G_{\tilde v}$ is undistorted in $G$ (see \cite[Lemma~4.1]{Mut24}), Lemmas~\ref{lem:svarc-milnor} and~\ref{lem:distortionuniversal} give constants $\lambda,B\ge1$ with $d_{\tilde v}(x,y)\le\lambda\,d(x,y)+B\le(\lambda+B)\,d(x,y)$. Thus $K:=\max\{\Lambda,\lambda+B\}$ works.
\end{proof}

\section{The model spaces for subgroups} \label{sec:model}
In this section, we fix the topological representative $M_\varphi$ of Section~\ref{subsection:lineartopological}. Let $G$ be the fundamental group of $M_\varphi$ and let $H$ be a finitely generated subgroup of $G$. In this section, we are going to construct a subspace $\Omega_H$ of the universal cover $\widetilde{M}_\varphi$ that will serve as a coarse geometric model for $H$. This will help us to reduce the problem of distortion of $H$ in $G$ to a geometric comparison at the level of the universal cover.

Let $w$ be a white vertex in the underlying graph $\Delta$ of $M_\varphi$. Recall from Proposition~\ref{prop:split-over-tori} that the vertex group $G_w$ of $G$ is of the form $G_w = F_w \times \langle t_w \rangle$ with $\mathrm{rk}(F_w) \ge 2$, and its centre is $Z(G_w) = \langle t_w \rangle$.

Let $p_w \colon G_w \to F_w$ be the projection given by $f t_{w}^k \mapsto f$ and $q_w \colon G_w \to \mathbb Z$ given by $f t_{w}^k \mapsto k$. In particular, we have $\mathrm{ker}(p_w) = \langle t_w \rangle$ and $\mathrm{ker}(q_w) = F_w$. Since $t_w\in Ft$ and $F_w\le F$, the map $q_w$ is the restriction to $G_w$ of the projection $G=F\rtimes\langle t\rangle\to\mathbb Z$.

Let $T$ be the Bass--Serre tree associated to the linear mapping torus decomposition given by Proposition~\ref{prop:split-over-tori}.
A vertex $\tilde v$ of $T$ is called \textit{white} (resp.\ \textit{black}) if its image in $\Delta$ is a white (resp.\ black) vertex of $\Delta$. The edges of $T$ inherit the orientation of the edges of $\Delta$, from the black vertex to the white vertex. Throughout the paper, we write an edge of $T$ as $\tilde e=(\tilde b,\tilde w)$, with $\tilde b$ black and $\tilde w$ white. 

For each vertex $v$ of $\Delta$ let $\bar v$ be the vertex of $T$ over $v$ with $\mathrm{Stab}_G(\bar v)=G_v$. Let $\tilde w$ be a white vertex of $T$ and let $w$ be its image in $\Delta$. Then $\tilde w=g\bar w$ for some $g\in G$, and $g$ is determined by $\tilde w$ up to right multiplication by $G_w$. We fix such an element $g=g_{\tilde w}$ for every white vertex $\tilde w$. Let $G_{\tilde w}:=\mathrm{Stab}_{G}( \tilde w) = g G_w g^{-1}$. We extend $p_w, q_w$ to $G_{\tilde w}$ via conjugation by $g$, and by abuse of notation, we still write $p_w, q_w$. Replacing $g$ by $gc$ with $c\in G_w$ changes $p_w$ by the inner automorphism of $F_w$ given by $p_w(c)$, and does not change $q_w$.

We define
\[
H_{\tilde w}:=H\cap G_{\tilde w},\quad P_{\tilde w}:=p_{w}(H_{\tilde w})\le F_w, \quad  J_{\tilde w}:=H_{\tilde w}\cap\ker(p_w)=H_{\tilde w}\cap Z(G_{\tilde w}).
\]

We then have the following short exact sequence
\[
1 \to J_{\tilde w} \to H_{\tilde w} \xrightarrow{p_w} P_{\tilde w} \to 1
\]
Note that $J_{\tilde w}=H\cap Z(G_{\tilde w})$ does not depend on the choice of $g_{\tilde w}$ and if $h \in H$ then we have $J_{ h \cdot \tilde{w}} = h \cdot J_{\tilde w} \cdot h^{-1}$.

\begin{defn}\label{defn:fibertype}
We say that the white vertex $\tilde w$ is of \textit{fiber type} if $J_{\tilde w}$ is trivial. Otherwise, we say that $\tilde w$ is of \textit{compact-fiber type}.
\end{defn}

Vertices of $T_H$ inherit their colours from $T$. Throughout the paper we draw them as follows. 

$$\circ \,\, \text{white, fiber type}, \quad \square \,\, \text{white, compact-fiber type} \quad \bullet \,\, \text{black} $$

We list below some facts that we often use in this paper.

\noindent \textbf{Ambient facts}:

\begin{enumerate}
\renewcommand{\theenumi}{D\arabic{enumi}}
\renewcommand{\labelenumi}{(\theenumi)}
    \item \label{D1} There is a $G$-equivariant, $1$-Lipschitz map $\widetilde\pi:\widetilde M_\varphi\to T$, where each edge of $T$ has length $2$. The map $\widetilde\pi$ collapses each vertex space to the corresponding vertex, and sends the point $(z,t)$ of the cylinder $\widetilde M^{\tilde b}\times[-1,1]$ of an edge $\tilde e=(\tilde b,\tilde w)$ to the point of $\tilde e$ at distance $1+t$ from $\tilde b$. The edges have length $2$ because the factor $[-1,1]$ of the cylinders has length $2$, and $\widetilde\pi$ is $1$-Lipschitz because the cylinders carry the $\ell^1$-product metric.
    
    Let $\widetilde M^{\tilde v}:=\widetilde\pi^{-1}(\tilde v)$ where $\tilde v$ is a vertex of $T$. Then $\widetilde M^{\tilde b}=\mathbb R^2$ when $\tilde b$ is  a black vertex and $\widetilde M^{\tilde w}=\widetilde X_w\times\mathbb R$ when $\tilde w$ is a white vertex. For an edge $\tilde e=(\tilde b,\tilde w)$, the space $\widetilde M^{\tilde b}\times[-1,1]$ is glued by $(z,-1)\sim z$ and $(z,1)\sim\iota_{\tilde e}(z)\in\widetilde M^{\tilde w}$, and $g\circ\iota_{\tilde e}=\iota_{g\tilde e}\circ g$ for all $g\in G$. 

\item \label{D2} $G$ acts on $(\widetilde M_\varphi,d)$ isometrically, properly and cocompactly.

\item \label{D3} For every vertex $\tilde v$ of $T$, the metrics $d$ and $d_{\tilde v}$ on $\widetilde M^{\tilde v}$ satisfy (\ref{eq:K}), and lengths of paths in $\widetilde M^{\tilde v}$ satisfy (\ref{eq:Klength}), where $K\ge1$ is the constant of Section~\ref{subsec:metrics}.

\item \label{D4} $G_{\tilde b}\cong\mathbb Z^2$ acts on $\widetilde M^{\tilde b}\cong\mathbb R^2$ by translations. For a white vertex $\tilde w$, we identify $\widetilde M^{\tilde w}$ with $\widetilde X_w\times\mathbb R$ via $g_{\tilde w}$. Then for $g\in G_{\tilde w}$ we have $g\cdot(\tilde x,s)=(p_w(g)\tilde x,\ s+q_w(g))$.

\item \label{D5} $G_{\tilde b}=G_{\tilde e}\le G_{\tilde w}$ for every edge $\tilde e=(\tilde b,\tilde w)$ (see Proposition~\ref{prop:split-over-tori}(\ref{prop:split-over-tori1})). Hence $H_{\tilde b}\le H_{\tilde w}$ and $\iota_{\tilde e}(uz)=u\cdot\iota_{\tilde e}(z)$ for all $u\in G_{\tilde b}$.

\end{enumerate}

\subsection{Horizontal pieces}
In order to construct the model space $\Omega_H$ for $H$ (see Definition~\ref{defn:model-space}), we introduce the following collection of ``horizontal pieces'' (see Sections~\ref{subsec:horizontalpiece_black} and~\ref{subsec:horizontalpiece_white}) in vertex spaces. Roughly speaking, the model space is obtained by gluing these horizontal pieces together.

\begin{rem} \label{rem:whyundistorted}
    If $H$ is conjugate into a vertex group $G_v$, then $H$ is undistorted in $G$. Indeed, every finitely generated subgroup of a vertex group $G_v$ is undistorted in $G_v$ \cite[Lemma~8.6]{Lin25}, and $G_v$ is undistorted in $G$ \cite[Lemma~4.1]{Mut24}. It follows that $H$ is undistorted in $G$.
\end{rem}

By Remark~\ref{rem:whyundistorted}, from now on we assume that $H$ is not conjugate into any vertex group of $G$.

Let $H \curvearrowright T$ be the action induced from the action $G \curvearrowright T$. Since $H$ is not conjugate into a vertex group, $H$ fixes no vertex of $T$. Since $H$ is finitely generated, $H$ contains an element acting loxodromically on $T$, and there is a unique minimal $H$-invariant subtree $T_H\subset T$. The tree $T_H$ is the union of the axes of the loxodromic elements of $H$. In particular, $T_H$ has no vertex of valence one. The quotient graph $\Lambda_H := H\backslash T_H$, which we call \textit{the underlying graph} of $H$, is finite. For each vertex $\tilde v$ in $T_H$ and for each edge $\tilde e$ in $T_H$, we define $H_{\tilde v} : = \mathrm{Stab}_{H}(\tilde v)$ and $H_{\tilde e} : = \mathrm{Stab}_{H} (\tilde{e})$. 

According to \cite[Theorem~1.3]{HW21}, if a finitely generated group splits as a graph of groups with finitely generated edge groups, then its vertex groups are also finitely generated. Applying this to the splitting of $H$, whose edge groups $H_{\tilde e}\le G_{\tilde e}\cong\mathbb Z^2$ are finitely generated, we see that each vertex group $H_{\tilde v}$ is finitely generated.

We first fix representatives $\tilde{b}_1, \ldots, \tilde{b}_{n_\bullet}$ of the $H$-orbits of black vertices of $T_H$ and representatives $\tilde{w}_1, \ldots, \tilde{w}_{n_\circ}$ of the $H$-orbits of white vertices of $T_H$.

\subsubsection{Construction of horizontal piece $\mathcal{H}(\tilde{v})$ when $\tilde{v}$ is black} \label{subsec:horizontalpiece_black}
For each $1\le j\le n_\bullet$, we fix a point $*_j\in\widetilde M^{\tilde b_j}\cong\mathbb R^2$ and let $\mathcal H(\tilde b_j)$ be the convex hull of the orbit $H_{\tilde b_j}\cdot *_j$. Since $H_{\tilde b_j}\le\mathbb Z^2$ acts on $\mathbb R^2$ by translations, $\mathcal H(\tilde b_j)$ is a point, a line or the whole plane, according to whether $H_{\tilde b_j}$ has rank $0$, $1$ or $2$. It is $H_{\tilde b_j}$-invariant. Since every black vertex of $T_H$ has the form $h\tilde b_j$ with $h\in H$, the set $\mathcal H(h\tilde b_j):=h\,\mathcal H(\tilde b_j)$ is well defined. Moreover, $\mathcal H(\tilde b_j)=H_{\tilde b_j}\cdot D_{\tilde b_j}$ for a compact set $D_{\tilde b_j}$, namely a point, a segment or a parallelogram. For $\tilde b=h\tilde b_j$ we fix such an $h$ and let $D_{\tilde b}:=hD_{\tilde b_j}$.

\subsubsection{Construction of horizontal piece $\mathcal{H}(\tilde{v})$ when $\tilde{v}$ is white} \label{subsec:horizontalpiece_white}

Let $1\le i\le n_\circ$, and let $w_i$ be the image of $\tilde w_i$ in $\Delta$. Let $X_{P_{\tilde w_i}}:=\widetilde X_{w_i}/P_{\tilde w_i}$ be the covering graph of $X_{w_i}$ corresponding to $P_{\tilde w_i}$, and let $\bar\rho\colon\widetilde X_{w_i}\to X_{P_{\tilde w_i}}$ be the quotient map. Let $C_{\tilde w_i}\subset X_{P_{\tilde w_i}}$ be the Stallings core with respect to a base vertex \cite{Sta83} (see also \cite[Lemma~2.1]{Lin25}), which is a single vertex if $P_{\tilde w_i}$ is trivial. Since $P_{\tilde w_i}$ is finitely generated, $C_{\tilde w_i}$ is a finite connected subgraph carrying $\pi_1(X_{P_{\tilde w_i}})\cong P_{\tilde w_i}$.

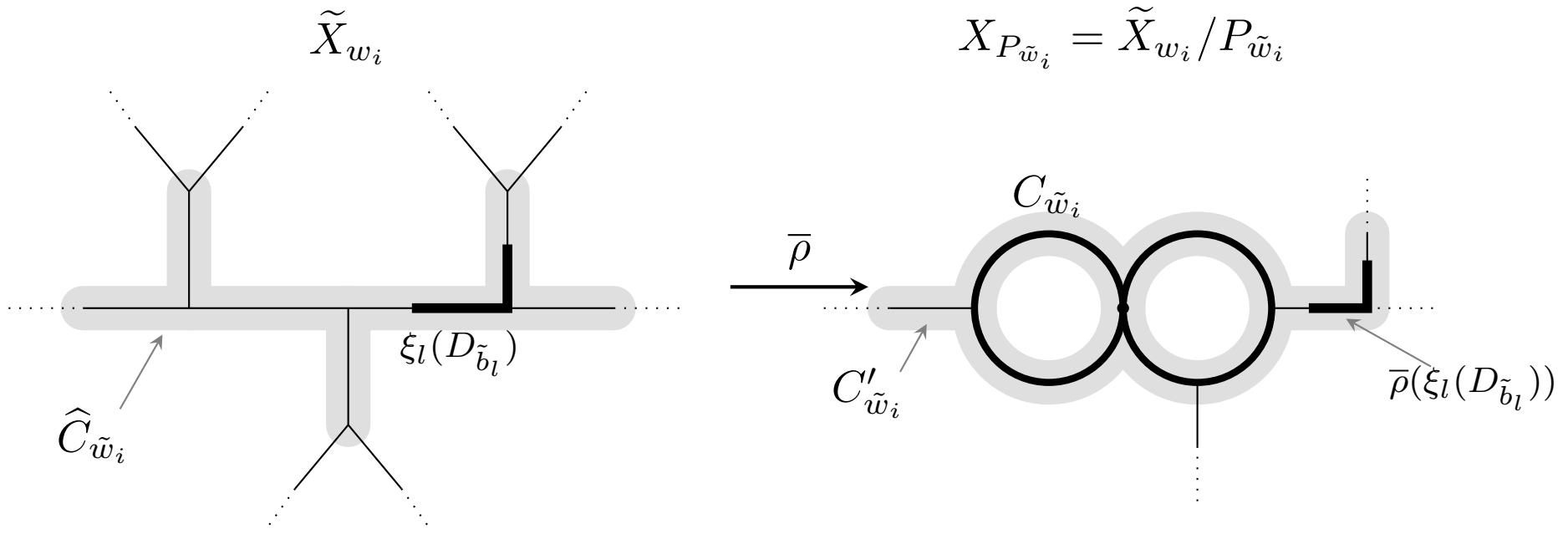
\begin{figure}[htb]
\centering
\begin{tikzpicture}[scale=0.84,>=stealth]
\begin{scope}
\draw[line width=10pt,gray!25,line cap=round,line join=round] (0.2,0)--(5.2,0);
\draw[line width=10pt,gray!25,line cap=round,line join=round] (1.2,0)--(1.2,1.1);
\draw[line width=10pt,gray!25,line cap=round,line join=round] (2.7,0)--(2.7,-1.1);
\draw[line width=10pt,gray!25,line cap=round,line join=round] (4.2,0)--(4.2,1.1);
\draw (0.2,0)--(5.2,0);
\draw (1.2,0)--(1.2,1.1); \draw (2.7,0)--(2.7,-1.1); \draw (4.2,0)--(4.2,1.1);
\draw[dotted] (-0.5,0)--(0.2,0); \draw[dotted] (5.2,0)--(5.9,0);
\draw (1.2,1.1)--(0.7,1.7); \draw (1.2,1.1)--(1.7,1.7);
\draw[dotted] (0.7,1.7)--(0.45,2.05); \draw[dotted] (1.7,1.7)--(1.95,2.05);
\draw (2.7,-1.1)--(2.2,-1.7); \draw (2.7,-1.1)--(3.2,-1.7);
\draw[dotted] (2.2,-1.7)--(1.95,-2.05); \draw[dotted] (3.2,-1.7)--(3.45,-2.05);
\draw (4.2,1.1)--(3.7,1.7); \draw (4.2,1.1)--(4.7,1.7);
\draw[dotted] (3.7,1.7)--(3.45,2.05); \draw[dotted] (4.7,1.7)--(4.95,2.05);
\draw[line width=2.2pt] (3.3,0)--(4.2,0)--(4.2,0.6);
\node[font=\scriptsize,below] at (3.75,-0.05) {$\xi_l(D_{\tilde b_l})$};
\node[font=\small] at (0.3,-1.2) {$\widehat C_{\tilde w_i}$};
\draw[->,gray] (0.55,-0.95)--(0.95,-0.25);
\node[font=\small] at (2.7,2.55) {$\widetilde X_{w_i}$};
\end{scope}
\draw[->,thick] (6.3,0.2)--(7.6,0.2) node[midway,above,font=\small] {$\bar\rho$};
\begin{scope}[shift={(10,0)}]
\draw[line width=10pt,gray!25,line cap=round,line join=round] (1.4,0)--(2.3,0)--(2.3,0.7);
\draw[line width=10pt,gray!25,line cap=round,line join=round] (-1.4,0)--(-2.2,0);
\draw[line width=10pt,gray!25] (-0.7,0) circle (0.7);
\draw[line width=10pt,gray!25] (0.7,0) circle (0.7);
\draw[line width=1.6pt] (-0.7,0) circle (0.7);
\draw[line width=1.6pt] (0.7,0) circle (0.7);
\fill (0,0) circle (1.6pt);
\draw (1.4,0)--(2.3,0)--(2.3,0.7); \draw (-1.4,0)--(-2.2,0);
\draw[dotted] (2.3,0.7)--(2.3,1.3); \draw[dotted] (2.3,0)--(3.0,0);
\draw[dotted] (-2.2,0)--(-2.9,0);
\draw (0.7,-0.7)--(0.7,-1.3); \draw[dotted] (0.7,-1.3)--(0.7,-1.9);
\draw[line width=2.2pt] (1.75,0)--(2.3,0)--(2.3,0.45);
\node[font=\scriptsize] at (3.3,-0.75) {$\bar\rho(\xi_l(D_{\tilde b_l}))$};
\draw[->,gray] (2.9,-0.55)--(2.1,-0.1);
\node[font=\small] at (-0.7,1.05) {$C_{\tilde w_i}$};
\node[font=\small] at (-2.4,-0.8) {$C'_{\tilde w_i}$};
\draw[->,gray] (-2.1,-0.6)--(-1.85,-0.15);
\node[font=\small] at (0,2.55) {$X_{P_{\tilde w_i}}=\widetilde X_{w_i}/P_{\tilde w_i}$};
\end{scope}
\end{tikzpicture}
\caption{The subtree $\widehat C_{\tilde w_i}$ (left, shaded) is the preimage of a finite connected subgraph $C'_{\tilde w_i}$ (right, shaded) containing the Stallings core $C_{\tilde w_i}$ (thick).}
\label{fig:Chat}
\end{figure}

Let $\tilde e_l=(\tilde b_l,\tilde w_i)$ with $1\le l\le r$ be representatives of the $H_{\tilde w_i}$-orbits of edges of $T_H$ at $\tilde w_i$, and write $\iota_{\tilde e_l}=(\xi_l,\tau_l)$, where $\xi_l$ takes values in the tree $\widetilde X_{w_i}$ and $\tau_l$ in $\mathbb R$. The sets $\bar\rho(\xi_l(D_{\tilde b_l}))$ are compact. So there is a finite connected subgraph $C'_{\tilde w_i}\subset X_{P_{\tilde w_i}}$ containing $C_{\tilde w_i}$ and all of them, and we let (see Figure~\ref{fig:Chat})
\[
\widehat C_{\tilde w_i}:=\bar\rho^{-1}(C'_{\tilde w_i})\subset\widetilde X_{w_i}.
\]
Since $C'_{\tilde w_i}$ is connected and contains the core, $\widehat C_{\tilde w_i}$ is a $P_{\tilde w_i}$-invariant subtree. The group $P_{\tilde w_i}$ acts freely on it with compact quotient $C'_{\tilde w_i}$, and $\widehat C_{\tilde w_i}$ contains every $\xi_l(D_{\tilde b_l})$.

\smallskip
\noindent Recall that a white vertex is of form $\circ$ (fiber type) or of form $\square$ (compact-fiber type). Let $\tilde w_i$ be of form $\circ$. Then $J_{\tilde w_i}=1$, so $p_{w_i}$ restricts to an isomorphism $H_{\tilde w_i}\to P_{\tilde w_i}$, and we define the homomorphism
\[
\vartheta_i:=q_{w_i}\circ\bigl(p_{w_i}|_{H_{\tilde w_i}}\bigr)^{-1}\colon P_{\tilde w_i}\to\mathbb Z ,
\]
that is, $\vartheta_i(p_{w_i}(g))=q_{w_i}(g)$ for $g\in H_{\tilde w_i}$. Since $P_{\tilde w_i}$ acts freely on $\widehat C_{\tilde w_i}$ with finitely many orbits of vertices, we may choose orbit representatives $\tilde y_1,\dots,\tilde y_n$ and define $f_i(p\,\tilde y_k):=\vartheta_i(p)$ for $p\in P_{\tilde w_i}$, extended affinely over the edges. Then
\begin{equation}\label{eq:f_i_equivariant}
f_i(p\tilde x)=f_i(\tilde x)+\vartheta_i(p)\quad\text{for all }\tilde x\in\widehat C_{\tilde w_i},\ p\in P_{\tilde w_i}.
\end{equation}
Since $\widehat C_{\tilde w_i}$ has finitely many $P_{\tilde w_i}$-orbits of edges, (\ref{eq:f_i_equivariant}) implies that $f_i$ is $L_i$-Lipschitz for some $L_i\ge0$.

Moreover, we have $\xi_l(\mathcal H(\tilde b_l))\subset\widehat C_{\tilde w_i}$ since $\xi_l(D_{\tilde b_l})\subset\widehat C_{\tilde w_i}$ and $\xi_l(uz)=p_{w_i}(u)\,\xi_l(z)$ for $u\in H_{\tilde b_l}$ (see (\ref{D4}) and (\ref{D5})). So we may define the continuous function
\[
\delta_l(z):=\tau_l(z)-f_i(\xi_l(z)),\qquad z\in\mathcal H(\tilde b_l).
\]
Let $R\ge0$ be the maximum of $|\delta_l|$ over the compact sets $D_{\tilde b_l}$, and let $L\ge0$ be the maximum of the constants $L_i$, both taken over all $i$ such that $\tilde w_i$ is of fiber type and all $l$ (we let $R=L=0$ if there is no such $i$).

\begin{defn}\label{defn:band}
For a white representative $\tilde w_i$ of fiber type, the \textit{band} of $\tilde w_i$ is
\[
\mathrm{Band}_R(\tilde w_i):=\bigl\{(\tilde x,s)\in\widehat C_{\tilde w_i}\times\mathbb R\ :\ |s-f_i(\tilde x)|\le R\bigr\}.
\]
The \textit{horizontal piece} of $\tilde w_i$ is $\mathcal H(\tilde w_i):=\mathrm{Band}_R(\tilde w_i)$ if $\tilde w_i$ is of form $\circ$, and $\mathcal H(\tilde w_i):=\widehat C_{\tilde w_i}\times\mathbb R$ if $\tilde w_i$ is of form $\square$. For $h\in H$ we let $\mathcal H(h\tilde w_i):=h\,\mathcal H(\tilde w_i)$.
\end{defn}

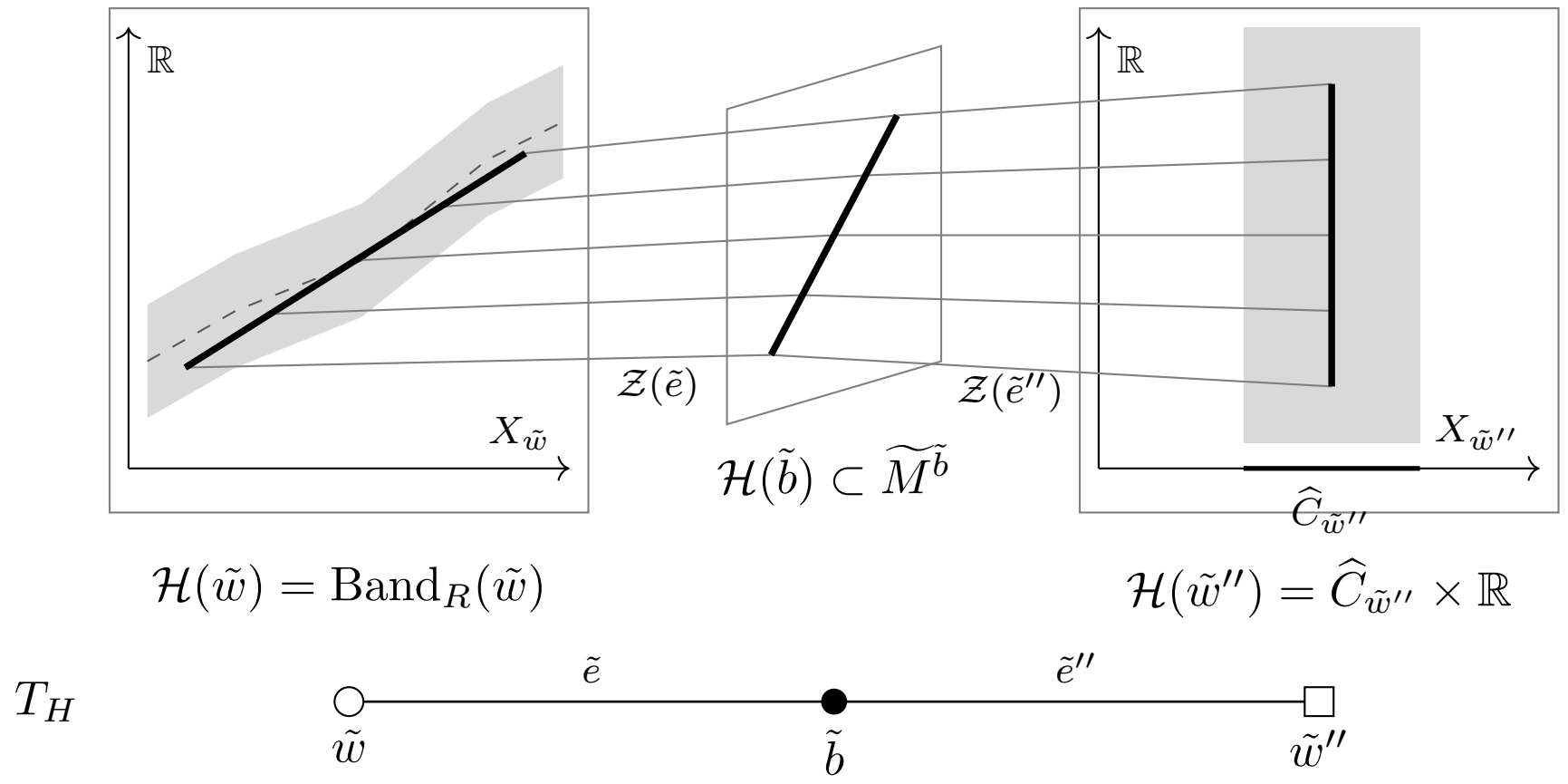
\begin{figure}[htb]
\centering
\begin{tikzpicture}[scale=0.92]
\draw[gray] (0,0) rectangle (3.8,4);
\draw[->] (0.15,0.35)--(3.65,0.35) node[above left,font=\scriptsize] {$X_{\tilde w}$};
\draw[->] (0.15,0.35)--(0.15,3.85) node[below right,font=\scriptsize] {$\mathbb R$};
\fill[gray!30] plot coordinates {(0.3,0.75) (1.0,1.15) (2.0,1.55) (3.0,2.35) (3.6,2.65)}
   -- plot coordinates {(3.6,3.55) (3.0,3.25) (2.0,2.45) (1.0,2.05) (0.3,1.65)} -- cycle;
\draw[gray!70!black,dashed] plot coordinates {(0.3,1.2) (1.0,1.6) (2.0,2.0) (3.0,2.8) (3.6,3.1)};
\node[font=\small] at (1.9,-0.6) {$\mathcal H(\tilde w)=\mathrm{Band}_R(\tilde w)$};
\draw[gray] (4.9,0.7)--(6.6,1.2)--(6.6,3.7)--(4.9,3.2)--cycle;
\node[font=\small] at (5.75,0.3) {$\mathcal H(\tilde b)\subset\widetilde M^{\tilde b}$};
\draw[gray] (7.7,0) rectangle (11.5,4);
\fill[gray!30] (9.0,0.55) rectangle (10.4,3.85);
\draw[->] (7.85,0.35)--(11.35,0.35) node[above left,font=\scriptsize] {$X_{\tilde w''}$};
\draw[->] (7.85,0.35)--(7.85,3.85) node[below right,font=\scriptsize] {$\mathbb R$};
\draw[line width=1pt] (9.0,0.35)--(10.4,0.35);
\node[font=\scriptsize,below] at (9.7,0.35) {$\widehat C_{\tilde w''}$};
\node[font=\small] at (9.6,-0.6) {$\mathcal H(\tilde w'')=\widehat C_{\tilde w''}\times\mathbb R$};
\coordinate (B1) at (5.25,1.25); \coordinate (B2) at (6.25,3.15);
\coordinate (L1) at (0.6,1.15);  \coordinate (L2) at (3.3,2.85);
\coordinate (R1) at (9.7,1.0);   \coordinate (R2) at (9.7,3.4);
\foreach \s in {0,0.25,0.5,0.75,1} {
  \draw[gray,thin] ($(B1)!\s!(B2)$)--($(L1)!\s!(L2)$);
  \draw[gray,thin] ($(B1)!\s!(B2)$)--($(R1)!\s!(R2)$);
}
\draw[line width=1.4pt] (B1)--(B2);
\draw[line width=1.4pt] (L1)--(L2);
\draw[line width=1.4pt] (R1)--(R2);
\node[font=\scriptsize] at (4.35,1.0) {$\mathcal Z(\tilde e)$};
\node[font=\scriptsize] at (7.15,0.95) {$\mathcal Z(\tilde e'')$};
\begin{scope}[shift={(0,-1.5)}]
\draw[nhor] (1.9,0)--(5.75,0) node[midway,above,font=\scriptsize] {$\tilde e$};
\draw[nhor] (5.75,0)--(9.6,0) node[midway,above,font=\scriptsize] {$\tilde e''$};
\node[wfib] at (1.9,0) {}; \node[blk] at (5.75,0) {}; \node[wcpt] at (9.6,0) {};
\node[font=\small,below=2pt] at (1.9,0) {$\tilde w$};
\node[font=\small,below=2pt] at (5.75,0) {$\tilde b$};
\node[font=\small,below=2pt] at (9.6,0) {$\tilde w''$};
\node[font=\small] at (-0.5,0) {$T_H$};
\end{scope}
\end{tikzpicture}
\caption{A portion of the model space $\Omega_H$ over a path $\tilde w-\tilde b-\tilde w''$ of $T_H$, with $\tilde w$ of fiber type, $\tilde w''$ of compact-fiber type and $\tilde b$ of type $(V)$.}
\label{fig:model}
\end{figure}

We refer the reader to Figure~\ref{fig:model} for an illustration of horizontal pieces associated to white vertices. 

\begin{prop} \label{prop:horizontal pieces}
 Let $\{ \mathcal{H}(\tilde v) \subset \widetilde {M}^{\tilde v} \}_{\tilde v \in V(T_H)}$ be the collection of horizontal pieces constructed in Sections~\ref{subsec:horizontalpiece_black} and~\ref{subsec:horizontalpiece_white}. Then there exists a constant $A \ge 1$ such that the following hold.
 \begin{enumerate}
     \item \label{prop:horizontal pieces_1} $\mathcal{H}(h \cdot \tilde v) = h \cdot \mathcal{H}(\tilde v)$ for $h \in H$.

     \item \label{prop:horizontal pieces_2} For every vertex $\tilde v$ of $T_H$ there exists a compact subset $D_{\tilde v}\subset\mathcal H(\tilde v)$ with $H_{\tilde v} \cdot D_{\tilde v} = \mathcal{H}(\tilde v)$.

     \item  \label{prop:horizontal pieces_3}   $\iota_{\tilde e} (\mathcal{H}(\tilde b)) \subset \mathcal{H}(\tilde w)$ for each edge $\tilde e  = (\tilde b, \tilde w)$ of $T_H$.

     \item  \label{prop:horizontal pieces_4}  For all points $x, y$ in $\mathcal{H}(\tilde v)$, there exists a path $\beta$ in $\mathcal{H}(\tilde v)$ connecting $x$ to $y$ such that $\ell_d(\beta) \le A\, d(x, y)$.
 \end{enumerate}
\end{prop}

\begin{proof}
The black horizontal pieces were treated in Section~\ref{subsec:horizontalpiece_black}, so we consider the white horizontal pieces. Recall from (\ref{D4}) that $g\in H_{\tilde w_i}$ acts by $g(\tilde x,s)=(p\tilde x,s+q_{w_i}(g))$, where $p=p_{w_i}(g)\in P_{\tilde w_i}$.

(\ref{prop:horizontal pieces_1}): Since $p\widehat C_{\tilde w_i}=\widehat C_{\tilde w_i}$, the set $\widehat C_{\tilde w_i}\times\mathbb R$ is $H_{\tilde w_i}$-invariant. So is $\mathrm{Band}_R(\tilde w_i)$, by (\ref{eq:f_i_equivariant}) and $\vartheta_i(p)=q_{w_i}(g)$. Hence $\mathcal{H}(h \cdot \tilde v) = h \cdot \mathcal{H}(\tilde v)$ for $h \in H$.

(\ref{prop:horizontal pieces_2}): Let $K_i\subset\widehat C_{\tilde w_i}$ be a compact set with $P_{\tilde w_i}K_i=\widehat C_{\tilde w_i}$. In the fiber-type case, $D_{\tilde w_i}:=\mathrm{Band}_R(\tilde w_i)\cap(K_i\times\mathbb R)$ has the required property by (\ref{prop:horizontal pieces_1}). In the compact-fiber case, $J_{\tilde w_i}$ is generated by $t_{\tilde w_i}^a$ for some $a\ge1$, where $t_{\tilde w_i}$ generates the centre of $G_{\tilde w_i}$, and $t_{\tilde w_i}^a$ acts by $(\tilde x,s)\mapsto(\tilde x,s+a)$, and $D_{\tilde w_i}:=K_i\times[0,a]$ has the required property. For $\tilde w=h\tilde w_i$ we let $D_{\tilde w}:=hD_{\tilde w_i}$.

(\ref{prop:horizontal pieces_3}): Since $\iota_{h\tilde e}\circ h=h\circ\iota_{\tilde e}$, it suffices by (\ref{prop:horizontal pieces_1}) to consider only the edges $\tilde e_l$. In the compact-fiber case, this is the inclusion $\xi_l(\mathcal H(\tilde b_l))\subset\widehat C_{\tilde w_i}$ shown above. In the fiber-type case, (\ref{D4}), (\ref{D5}) and (\ref{eq:f_i_equivariant}) show that $\delta_l$ is $H_{\tilde b_l}$-invariant. So $|\delta_l|\le R$ on $H_{\tilde b_l}D_{\tilde b_l}=\mathcal H(\tilde b_l)$, that is, $\iota_{\tilde e_l}(\mathcal H(\tilde b_l))\subset\mathrm{Band}_R(\tilde w_i)$.

(\ref{prop:horizontal pieces_4}): By (\ref{prop:horizontal pieces_1}) it suffices to consider the representatives. Let $K\ge1$ be the constant given by Section~\ref{subsec:metrics} and let $L$ be the constant given above Definition~\ref{defn:band}. Let $A:=2(1+L)K^2$ . 

If $\tilde v$ is black, let $\beta$ be the straight segment from $x$ to $y$ in the convex set $\mathcal H(\tilde v)$, which is a $d_{\tilde v}$-geodesic. If $\tilde v=\tilde w_i$ is white, write $x=(\tilde x_1,s_1)$ and $y=(\tilde x_2,s_2)$, and let $\sigma\colon[0,1]\to\widehat C_{\tilde w_i}$ be the geodesic from $\tilde x_1$ to $\tilde x_2$, parametrized proportionally to arc length.

In the compact-fiber case, let $\beta(\tau):=(\sigma(\tau),(1-\tau)s_1+\tau s_2)$. In the fiber-type case, we interpolate the offsets from $f_i$. That is, we let $u_k:=s_k-f_i(\tilde x_k)$ and $\beta(\tau):=\bigl(\sigma(\tau),f_i(\sigma(\tau))+(1-\tau)u_1+\tau u_2\bigr)$. Then $\beta$ stays in $\mathrm{Band}_R(\tilde w_i)$. For simplicity, we write $d_X$ for the metric of $\widetilde X_{w_i}$. Since $f_i$ is $L$-Lipschitz, it follows that
\[
\ell_{d_{\tilde w_i}}(\beta)\le(1+L)\,d_X(\tilde x_1,\tilde x_2)+|u_1-u_2|\le(1+2L)\,d_X(\tilde x_1,\tilde x_2)+|s_1-s_2| .
\]
In all cases $\ell_{d_{\tilde v}}(\beta)\le2(1+L)\,d_{\tilde v}(x,y)$, and (\ref{D3}) gives $\ell_d(\beta)\le A\,d(x,y)$.
\end{proof}

\subsection{The model space}
For an edge $\tilde e=(\tilde b,\tilde w)$ of $T$, let $\Phi_{\tilde e}\colon\widetilde M^{\tilde b}\times[-1,1]\to\widetilde M_\varphi$ be the map of the cylinder of $\tilde e$ given by the gluing in (\ref{D1}), so that $\Phi_{\tilde e}(z,-1)=z$ and $\Phi_{\tilde e}(z,1)=\iota_{\tilde e}(z)$. For an edge $\tilde e=(\tilde b,\tilde w)$ of $T_H$, the \emph{strip} of $\tilde e$ is
\[
\mathcal Z(\tilde e):=\Phi_{\tilde e}\bigl(\mathcal H(\tilde b)\times[-1,1]\bigr).
\]
Its two ends are $\mathcal H(\tilde b)$ and $\iota_{\tilde e}(\mathcal H(\tilde b))$, and the latter lies in $\mathcal H(\tilde w)$ by Proposition~\ref{prop:horizontal pieces}(\ref{prop:horizontal pieces_3}).

\begin{defn}\label{defn:model-space}
The \emph{model space for $H$} is
\[
\Omega_H:=\bigcup_{\tilde v\in V(T_H)}\mathcal H(\tilde v)\ \cup\bigcup_{\tilde e\in E(T_H)}\mathcal Z(\tilde e)\ \subset\ \widetilde M_\varphi ,
\]
equipped with the induced length metric $\widehat d_{\Omega_H}$. That is, $\widehat d_{\Omega_H}(x,y)$ is the infimum of the $d$-lengths of paths in $\Omega_H$ from $x$ to $y$.
\end{defn}

Thus $\Omega_H$ is obtained from the horizontal pieces by attaching a strip along each edge of $T_H$. Figure~\ref{fig:model} shows a portion of $\Omega_H$ over a path of length two in $T_H$.

Recall from (\ref{D1}) that $\widetilde\pi\colon\widetilde M_\varphi\to T$ is a $G$-equivariant $1$-Lipschitz map. We list some basic properties of the model space $\Omega_H$ of $H$ in the following lemma.

\begin{lem} \label{lem:basisalmostfiber}
Let $\Pi:=\widetilde\pi|_{\Omega_H}:\Omega_H\to T_H$, and fix a base point $x_0\in\Omega_H$. Then:
\begin{enumerate}
    \item \label{lem:basisalmostfiber1} $h\,\Omega_H=\Omega_H$ for all $h\in H$, and $H$ acts on $(\Omega_H,\widehat d_{\Omega_H})$ by isometries.
    \item \label{lem:basisalmostfiber2} $\Omega_H\cap\Pi^{-1}(\tilde v)=\mathcal H(\tilde v)$ for every vertex $\tilde v$.
    \item \label{lem:basisalmostfiber3} $d\le\widehat d_{\Omega_H}<\infty$, and $(\Omega_H,\widehat d_{\Omega_H})$ is a length space.
    \item \label{lem:basisalmostfiber4} For every $r\ge0$ the set $\{h\in H:\ \widehat d_{\Omega_H}(x_0,hx_0)\le r\}$ is finite.
    \item \label{lem:basisalmostfiber5} There is $r_0>0$ with $H\cdot B_{\widehat d}(x_0,r_0)=\Omega_H$.
\end{enumerate}
\end{lem}

\begin{proof}
(\ref{lem:basisalmostfiber1}): By Proposition~\ref{prop:horizontal pieces}(\ref{prop:horizontal pieces_1}) and the relation $h\circ\iota_{\tilde e}=\iota_{h\tilde e}\circ h$ of (\ref{D1}), we have $h\,\mathcal H(\tilde v)=\mathcal H(h\tilde v)$ and $h\,\mathcal Z(\tilde e)=\mathcal Z(h\tilde e)$. Hence $h\,\Omega_H=\Omega_H$, and $h$ preserves $\widehat d_{\Omega_H}$ because it is a $d$-isometry mapping paths in $\Omega_H$ to paths in $\Omega_H$.

(\ref{lem:basisalmostfiber2}): Each piece $\mathcal H(\tilde v')$ lies in $\widetilde\pi^{-1}(\tilde v')$, and the only points of a strip over vertices of $T$ are its two ends, which lie in pieces (see Proposition~\ref{prop:horizontal pieces}(\ref{prop:horizontal pieces_3})).

(\ref{lem:basisalmostfiber3}): Clearly $d\le\widehat d_{\Omega_H}$, and $\widehat d_{\Omega_H}$ is a length metric by definition. Points of a piece are joined inside it by rectifiable paths (see Proposition~\ref{prop:horizontal pieces}(\ref{prop:horizontal pieces_4})). Every point of a strip $\mathcal Z(\tilde e)$ lies on a vertical segment $t\mapsto\Phi_{\tilde e}(z,t)$ of length at most $2$ (see (\ref{eq:vertical-segment})) joining the two adjacent pieces. Since $T_H$ is connected, $\widehat d_{\Omega_H}<\infty$.

(\ref{lem:basisalmostfiber4}): By (\ref{lem:basisalmostfiber3}), this set is contained in $\{g\in G:\ d(x_0,gx_0)\le r\}$, which is finite by (\ref{D2}).

(\ref{lem:basisalmostfiber5}): Let $\tilde v_1,\dots,\tilde v_n$ represent the finitely many $H$-orbits of vertices of $T_H$, and let $D:=\bigcup_jD_{\tilde v_j}$ (see Proposition~\ref{prop:horizontal pieces}(\ref{prop:horizontal pieces_2})). Then every piece is contained in $H\cdot D$, and $D$ has finite $\widehat d_{\Omega_H}$-diameter by Proposition~\ref{prop:horizontal pieces}(\ref{prop:horizontal pieces_4}). By (\ref{lem:basisalmostfiber1}), every point of a piece lies within $\widehat d_{\Omega_H}$-distance $\rho:=\sup_{x\in D}\widehat d_{\Omega_H}(x_0,x)<\infty$ of $Hx_0$. By the vertical segments in (\ref{lem:basisalmostfiber3}), every point of a strip lies within distance $2$ of a piece. Hence $r_0:=\rho+3$ works.
\end{proof}

\begin{defn}\label{defn:vertical}
Let $\tilde w\in T$ be a white vertex over $w\in\Delta$, and write $G_{\tilde w}=gG_wg^{-1}$ with $G_w=F_w\times\langle t_w\rangle$ and $g=g_{\tilde w}$. The \textit{vertical element} of $\tilde w$ is $t_{\tilde w}:=g\,t_w\,g^{-1}\in G_{\tilde w}$.
\end{defn}

\begin{lem} \label{lem:basis_vertical}
The element $t_{\tilde w}$ is well defined, and:
\begin{enumerate}
    \item \label{lem:basis_vertical_1} $t_{\tilde w}$ commutes with every element of $G_{\tilde w}$, and $\langle t_{\tilde w}\rangle=Z(G_{\tilde w})=\ker p_w$.
    \item \label{lem:basis_vertical_2} $t_{h\tilde w}=h\,t_{\tilde w}\,h^{-1}$ for every $h\in G$.
    \item \label{lem:basis_vertical_3} $t_{\tilde w}$ acts on $\widetilde M^{\tilde w}=\widetilde X_w\times\mathbb R$ by $t_{\tilde w}\cdot(\tilde x,s)=(\tilde x,s+1)$.
    \item \label{lem:basis_vertical_4} If $\tilde w$ is of fiber type, then $H\cap\langle t_{\tilde w}\rangle=1$.
\end{enumerate}
\end{lem}

\begin{proof}
Replacing $g$ by $gc$ with $c\in G_w$ gives $gct_wc^{-1}g^{-1}=gt_wg^{-1}$, since $t_w$ is central in $G_w$; so $t_{\tilde w}$ is well defined. The same computation with $hg$ in place of $g$ gives (\ref{lem:basis_vertical_2}). Item (\ref{lem:basis_vertical_1}) follows from $Z(G_w)=\langle t_w\rangle=\ker p_w$, and item (\ref{lem:basis_vertical_3}) from (\ref{D4}), since $p_w(t_{\tilde w})=1$ and $q_w(t_{\tilde w})=1$. For (\ref{lem:basis_vertical_4}), $H\cap\langle t_{\tilde w}\rangle\subset H\cap G_{\tilde w}\cap\ker p_w=J_{\tilde w}$, which is trivial when $\tilde w$ is of fiber type.
\end{proof}

\begin{defn} [Coordinates]\label{defn:coordinates}
Let $\tilde b\in T$ be black, and let $b\in\Delta$ be the image of $\tilde b$ under the map $T \to \Delta$. Write $\tilde b=g\bar b$ with $g\in G$. The stabilizer of $\tilde b$ is $G_{\tilde b}=gG_bg^{-1}$ with $G_b=\langle r_b\rangle\times\langle t_b\rangle\cong\mathbb Z^2$. The isomorphism
$G_{\tilde b}\to G_b,\quad x\mapsto g^{-1}xg$
does not depend on the choice of $g$, because $G_b$ is abelian. So every $x\in G_{\tilde b}$ has well-defined \textit{coordinates} $(X,Y)\in\mathbb Z^2$ with $$g^{-1}xg=r_b^{X}t_b^{Y}.$$   
\end{defn}

\begin{lem} \label{lem:slope_vertical}
    Let $\tilde e=(\tilde b,\tilde w)$ be an edge of $T$ over $e=(b,w)$. Then $t_{\tilde w}\in G_{\tilde b}$, and its coordinates are $(-k(e),\,1)$.
\end{lem}

\begin{proof}
 By (\ref{D5}), $G_{\tilde b}=G_{\tilde e}\le G_{\tilde w}$. By Proposition~\ref{prop:split-over-tori}(\ref{prop:split-over-tori4}), under the identifications $G_{\tilde b}\cong G_b$ and $G_{\tilde w}\cong G_w$ this inclusion is, up to an inner automorphism of $G_w$,
$$r_b\mapsto p_e,\quad t_b\mapsto p_e^{k(e)}t_w .$$
Hence the element of $G_{\tilde b}$ with coordinates $(-k(e),1)$ is sent to a conjugate of
$$p_e^{-k(e)}\cdot p_e^{k(e)}t_w=t_w,$$
that is, to $t_{\tilde w}$, because $t_w$ is central in $G_w$.
\end{proof}

\begin{lem}[Types of black vertices]\label{lem:black_types}
Let $\tilde b$ be a black vertex of $T_H$. Then exactly one of the following
holds:
\begin{enumerate}
\item[(2)] $H_{\tilde b}\cong\mathbb Z^2$, and every neighbour of
$\tilde b$ in $T$ is of compact-fiber type;
\item[$(V)$] $H_{\tilde b}\cong\mathbb Z$, and exactly one neighbour
$\tilde w''$ of $\tilde b$ in $T$ is of compact-fiber type. Moreover,
$H_{\tilde b}\cap\langle t_{\tilde w''}\rangle$ has finite index in
$H_{\tilde b}$;
\item[$(H)$] $H_{\tilde b}\cong\mathbb Z$, and every neighbour of
$\tilde b$ in $T$ is of fiber type;
\item[(0)] $H_{\tilde b}=1$, and every neighbour of $\tilde b$ in $T$
is of fiber type.
\end{enumerate}
We say that $\tilde b$ is of type (2), $(V)$, $(H)$, or (0), respectively.
\end{lem}

\begin{figure}[htb]
\centering
\begin{tikzpicture}[scale=0.86]
\begin{scope}
\draw[nhor] (0,0)--(-1.2,0.8); \draw[nhor] (0,0)--(1.2,0.8); \draw[nhor] (0,0)--(0,-1.3);
\node[wfib] at (-1.2,0.8) {}; \node[wfib] at (1.2,0.8) {}; \node[wfib] at (0,-1.3) {};
\node[blk] at (0,0) {};
\node[font=\small,right=2pt] at (0,0.05) {$\tilde b$};
\node[font=\small] at (0,1.5) {type $(H)$};
\node[font=\footnotesize,align=center] at (0,-2.2) {$H_{\tilde b}\cong\mathbb Z$\\ all neighbours of form $\circ$};
\end{scope}
\begin{scope}[shift={(5.2,0)}]
\draw[nhor] (0,0)--(-1.2,0.8); \draw[nhor] (0,0)--(1.2,0.8); \draw[nhor] (0,0)--(0,-1.3);
\node[wcpt] at (-1.2,0.8) {}; \node[wcpt] at (1.2,0.8) {}; \node[wcpt] at (0,-1.3) {};
\node[blk] at (0,0) {};
\node[font=\small,right=2pt] at (0,0.05) {$\tilde b$};
\node[font=\small] at (0,1.5) {type (2)};
\node[font=\footnotesize,align=center] at (0,-2.2) {$H_{\tilde b}\cong\mathbb Z^2$\\ all neighbours of form $\square$};
\end{scope}
\begin{scope}[shift={(10.4,0)}]
\draw[nhor] (0,0)--(-1.2,0.8); \draw[nhor] (0,0)--(1.2,0.8); \draw[nhor] (0,0)--(0,-1.3);
\node[wfib] at (-1.2,0.8) {}; \node[wfib] at (1.2,0.8) {}; \node[wcpt] at (0,-1.3) {};
\node[blk] at (0,0) {};
\node[font=\small,right=2pt] at (0,0.05) {$\tilde b$};
\node[font=\small,right=3pt] at (0,-1.3) {$\tilde w''$};
\node[font=\small] at (0,1.5) {type $(V)$};
\node[font=\footnotesize,align=center] at (0,-2.2) {$H_{\tilde b}\cong\mathbb Z$\\ exactly one neighbour\\ of form $\square$};
\end{scope}
\end{tikzpicture}
\caption{The types $(H)$, (2) and $(V)$ of a black vertex $\tilde b$ of $T_H$ (see Lemma~\ref{lem:black_types}). Type (0) looks like type $(H)$, with $H_{\tilde b}=1$.}
\label{fig:black_types}
\end{figure}
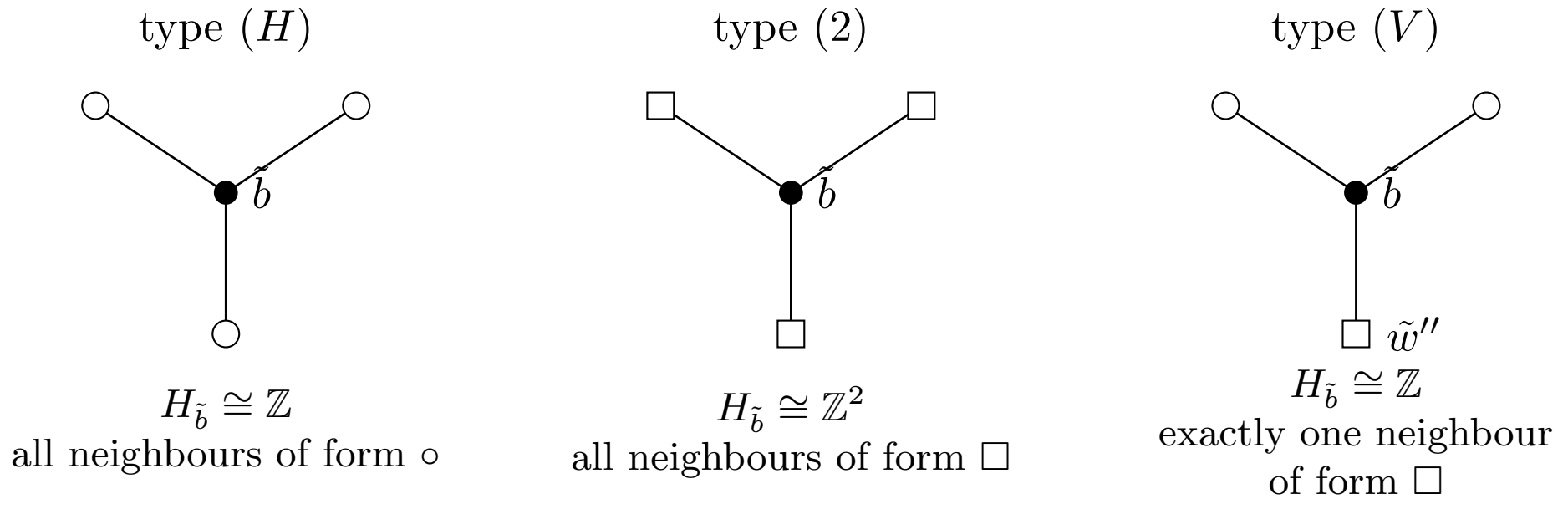

Figure~\ref{fig:black_types} illustrates the three types of black vertices with $H_{\tilde b}\neq1$.

\begin{proof}
Since $H_{\tilde b}\le G_{\tilde b}\cong\mathbb Z^2$, the rank of
$H_{\tilde b}$ is $0$, $1$, or $2$. Let $\tilde w$ be a neighbour of
$\tilde b$ in $T$, lying over an edge $e$. If $\tilde w$ is of
compact-fiber type, then
$1\neq J_{\tilde w}=\langle t_{\tilde w}^a\rangle$ for some $a\ge1$.
By Lemma~\ref{lem:slope_vertical}, we have $t_{\tilde w}\in G_{\tilde b}$,
and hence
\[
t_{\tilde w}^a\in H\cap G_{\tilde b}=H_{\tilde b}.
\]

Since $G_{\tilde b}=G_{\tilde e}$ for every edge $\tilde e$ of $T$ at $\tilde b$ (see Proposition~\ref{prop:split-over-tori}(\ref{prop:split-over-tori1})), each edge of $\Delta$ at $b$ lifts to exactly one edge of $T$ at $\tilde b$. So distinct edges of $T$ at $\tilde b$ lie over distinct edges of $\Delta$. Suppose that two distinct neighbours $\tilde w$ and $\tilde w'$ of
$\tilde b$, lying over distinct edges $e$ and $e'$, are of compact-fiber
type. Then $H_{\tilde b}$ contains both $t_{\tilde w}^a$ and
$t_{\tilde w'}^{a'}$ for some $a,a'\ge1$. Their coordinate vectors,
$a(-k(e),1)$ and $a'(-k(e'),1)$, are linearly independent because
$k(e)\neq k(e')$. It follows that $H_{\tilde b}\cong\mathbb Z^2$.

Conversely, suppose that $H_{\tilde b}\cong\mathbb Z^2$. Then
$H_{\tilde b}$ has finite index in $G_{\tilde b}$, so it contains a
positive power of $t_{\tilde w}$ for every neighbour $\tilde w$ of
$\tilde b$. Thus every neighbour of $\tilde b$ is of compact-fiber type.

If exactly one neighbour $\tilde w''$ of $\tilde b$ is of compact-fiber
type, then $H_{\tilde b}$ is nontrivial and contains
$t_{\tilde w''}^a$ for some $a\ge1$. Since $H_{\tilde b}$ has rank $1$,
the subgroup $\langle t_{\tilde w''}^a\rangle$ has finite index in
$H_{\tilde b}\cong\mathbb Z$.

Finally, if every neighbour of $\tilde b$ is of fiber type, then
$H_{\tilde b}$ has rank at most $1$.
\end{proof}

\begin{defn}[Horizontal edges and horizontal black vertices] \label{defn:horizontal_black}
An edge $\tilde e=(\tilde b,\tilde w)$ of $T_H$ is called \textit{horizontal} if $H_{\tilde b}$ is infinite cyclic and $\tilde w$ is of fiber type (i.e., a white vertex of form $\circ$). A black vertex $\tilde b$ of $T_H$ is called \textit{horizontal} if at least two edges of $T_H$ at $\tilde b$ are horizontal. We denote $$\mathcal{B}_H : = \{ \text{horizontal black vertices of type $(V)$ or $(H)$ of $T_H$} \} $$
\end{defn}

Since $J_{h\tilde w}=hJ_{\tilde w}h^{-1}$ and $H_{h\tilde b}=hH_{\tilde b}h^{-1}$, the set of horizontal edges and the set $\mathcal B_H$ are $H$-invariant. When every black vertex of $\Delta$ has valence two, a black vertex $\tilde b$ of $T_H$ is horizontal if and only if $H_{\tilde b}$ is infinite cyclic and both white vertices adjacent to $\tilde b$ are of fiber type.

By Lemma~\ref{lem:black_types}, an edge $\tilde e=(\tilde b,\tilde w)$ of $T_H$ is horizontal if and only if $\tilde b$ is of type $(V)$ or $(H)$ and $\tilde w$ is of fiber type. In particular, every horizontal black vertex is of type $(V)$ or $(H)$.

The black vertices of types $(V)$ and $(H)$ are the ones that contribute to the distortion of $H$, as we will show later. Let $\Pi\colon\Omega_H\to T_H$ be the map of Lemma~\ref{lem:basisalmostfiber}.

\begin{defn}\label{defn:hatYandOmegaTi}
Let $\widehat Y\subset T_H$ be the union of all closed horizontal edges of $T_H$ whose black endpoint lies in $\mathcal B_H$. It is possible that $\widehat Y$ is disconnected. Let $T_i$
be a connected component of $\widehat Y$, and set
\[
\Omega_{T_i}:=\Omega_H\cap\Pi^{-1}(T_i).
\]
We equip $\Omega_{T_i}$ with the metric
\[
\widehat d_{\Omega_{T_i}}(x,y)
:=\inf\{\ell_d(\gamma):\ \gamma\subset\Omega_{T_i}
\text{ is a path from }x\text{ to }y\}.
\]
\end{defn}

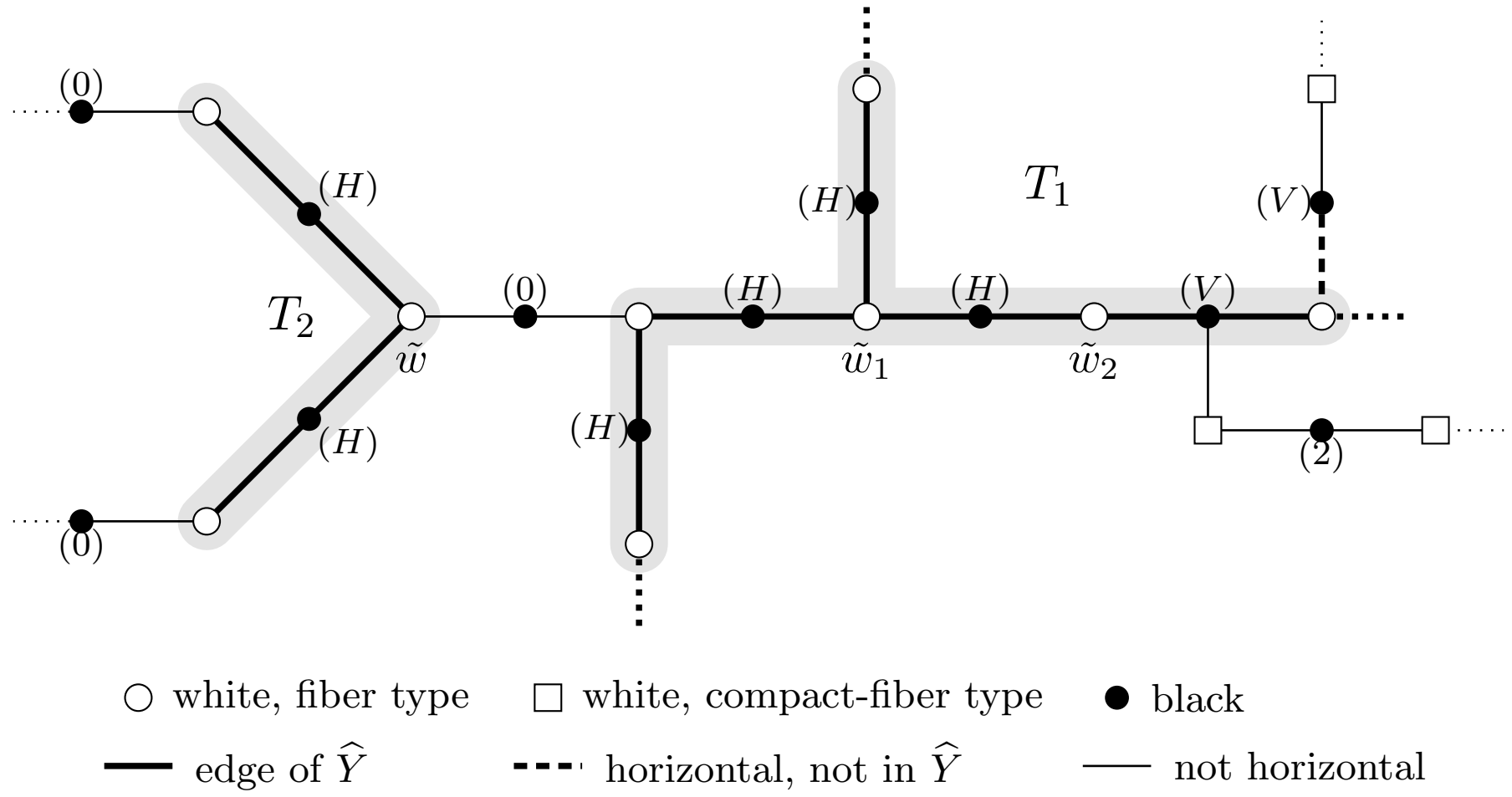
\begin{figure}[htb]
\centering
\begin{tikzpicture}[scale=0.9,
  wfib/.style={circle,draw,fill=white,inner sep=0pt,minimum size=6pt},
  wcpt/.style={rectangle,draw,fill=white,inner sep=0pt,minimum size=6pt},
  blk/.style={circle,draw,fill=black,inner sep=0pt,minimum size=5pt},
  hor/.style={line width=1.4pt},
  nhor/.style={line width=0.5pt},
  hornotY/.style={line width=1.4pt,dash pattern=on 3pt off 2pt},
  cont/.style={line width=1.4pt,dotted},
  compbg/.style={line width=13pt,gray!22,line cap=round,line join=round},
  typ/.style={font=\scriptsize,inner sep=1.5pt}]
\coordinate (c)  at (0,0);
\coordinate (b6) at (-0.9,0.9);  \coordinate (w8) at (-1.8,1.8);
\coordinate (b7) at (-0.9,-0.9); \coordinate (w9) at (-1.8,-1.8);
\coordinate (z8) at (-2.9,1.8);  \coordinate (z9) at (-2.9,-1.8);
\coordinate (b5) at (1,0);
\coordinate (w1) at (2,0);  \coordinate (b1) at (3,0);
\coordinate (w2) at (4,0);  \coordinate (b2) at (5,0);
\coordinate (w3) at (6,0);  \coordinate (b4) at (7,0);
\coordinate (w5) at (8,0);
\coordinate (b3) at (4,1);  \coordinate (w4) at (4,2);
\coordinate (b10) at (2,-1); \coordinate (w12) at (2,-2);
\coordinate (w6) at (7,-1); \coordinate (b8) at (8,-1); \coordinate (w10) at (9,-1);
\coordinate (b9) at (8,1);  \coordinate (w11) at (8,2);
\draw[compbg] (w8)--(b6)--(c)--(b7)--(w9);
\draw[compbg] (w12)--(b10)--(w1)--(b1)--(w2)--(b2)--(w3)--(b4)--(w5);
\draw[compbg] (w2)--(b3)--(w4);
\node at (-1.05,0) {$T_2$};
\node at (5.6,1.15) {$T_1$};
\draw[hor] (w8)--(b6)--(c)--(b7)--(w9);
\draw[hor] (w12)--(b10)--(w1)--(b1)--(w2)--(b2)--(w3)--(b4)--(w5);
\draw[hor] (w2)--(b3)--(w4);
\draw[cont] (w4)--++(0,0.75);
\draw[cont] (w12)--++(0,-0.75);
\draw[cont] (w5)--++(0.75,0);
\draw[hornotY] (w5)--(b9);
\draw[nhor] (b9)--(w11);
\draw[nhor] (c)--(b5)--(w1);
\draw[nhor] (b4)--(w6)--(b8)--(w10);
\draw[nhor] (w8)--(z8);  \draw[nhor] (w9)--(z9);
\draw[nhor,dotted] (z8)--++(-0.6,0); \draw[nhor,dotted] (z9)--++(-0.6,0);
\draw[nhor,dotted] (w11)--++(0,0.6); \draw[nhor,dotted] (w10)--++(0.6,0);
\foreach \v in {c,w8,w9,w1,w2,w3,w5,w4,w12} \node[wfib] at (\v) {};
\foreach \v in {w6,w10,w11} \node[wcpt] at (\v) {};
\foreach \v in {b6,b7,b1,b2,b3,b10,b4,b5,b8,b9,z8,z9} \node[blk] at (\v) {};

\node[typ,above right] at (b6) {$(H)$};
\node[typ,below right] at (b7) {$(H)$};
\node[typ,above]       at (z8) {(0)};
\node[typ,below]       at (z9) {(0)};
\node[typ,above]       at (b5) {(0)};
\node[typ,above]       at (b1) {$(H)$};
\node[typ,above]       at (b2) {$(H)$};
\node[typ,left]        at (b3) {$(H)$};
\node[typ,left]        at (b10) {$(H)$};
\node[typ,above]       at (b4) {$(V)$};
\node[typ,below]       at (b8) {(2)};
\node[typ,left]        at (b9) {$(V)$};

\node[font=\small,below=2pt] at (c) {$\tilde w$};
\node[font=\small,below=2pt] at (w2) {$\tilde w_1$};
\node[font=\small,below=2pt] at (w3) {$\tilde w_2$};

\begin{scope}[shift={(-2.4,-3.35)},font=\footnotesize]
\node[wfib] at (0,0) {}; \node[right] at (0.15,0) {white, fiber type};
\node[wcpt] at (3.6,0) {}; \node[right] at (3.75,0) {white, compact-fiber type};
\node[blk] at (8.6,0) {}; \node[right] at (8.75,0) {black};
\draw[hor] (-0.3,-0.6)--(0.3,-0.6); \node[right] at (0.35,-0.6) {edge of $\widehat Y$};
\draw[hornotY] (3.3,-0.6)--(3.9,-0.6); \node[right] at (3.95,-0.6) {horizontal, not in $\widehat Y$};
\draw[nhor] (8.3,-0.6)--(8.9,-0.6); \node[right] at (8.95,-0.6) {not horizontal};
\end{scope}
\end{tikzpicture}
\caption{A portion of $T_H$ and of $\widehat Y$, with the types of the black vertices (see Lemma~\ref{lem:black_types}). The shaded subtrees lie in two horizontal components $T_1$ and $T_2$. The dashed edge is horizontal, but its black endpoint has only one horizontal edge, so it is not in $\widehat Y$.}
\label{fig:hatY}
\end{figure}

We refer the reader to Figure~\ref{fig:hatY} for an illustration of Definition~\ref{defn:hatYandOmegaTi}.

\begin{rem}\label{rem:no_square}
The subgraph $\widehat Y$ does not contain any white vertex of form $\square$, since every white vertex of $\widehat Y$ is the endpoint of a horizontal edge.
\end{rem}

\begin{lem}\label{lem:vertices_of_Ti}
The subgraph $\widehat Y$ is $H$-invariant. Let $T_i$ be a connected component of $\widehat Y$. Then:
\begin{enumerate}
    \item \label{lem:vertices_of_Ti_1} every black vertex $\tilde b$ of $T_i$ belongs to $\mathcal B_H$,
    \item \label{lem:vertices_of_Ti_2} every white vertex of $T_i$ is of form $\circ$ (i.e., of fiber type), and every edge of $T_i$ is horizontal.
\end{enumerate}
\end{lem}

\begin{proof}
The subgraph $\widehat Y$ is $H$-invariant because $\mathcal B_H$, the set of horizontal edges, and $T_H$ are all $H$-invariant.
For (\ref{lem:vertices_of_Ti_1}), every black vertex $\tilde b$ of $\widehat Y$ belongs to $\mathcal B_H$ by construction, and the edges of $T_i$ at $\tilde b$ are precisely the horizontal edges of $T_H$ at $\tilde b$.
For (\ref{lem:vertices_of_Ti_2}), every edge of $\widehat Y$ is horizontal, and every white vertex of $\widehat Y$ is the white endpoint of a horizontal edge, so it is of fiber type.
\end{proof}

\begin{defn} \label{defn:horizontal_component} The connected components $T_i$ of $\widehat Y$ (see Definition~\ref{defn:hatYandOmegaTi}) are called the \textit{horizontal components} of $H$.
    
\end{defn}

\begin{lem} \label{lem:basifactsofTi}
    Let $T_i$ be a connected component of $\widehat Y$, and $H_i:=\mathrm{Stab}_H(T_i)$.
    \begin{enumerate}
        \item \label{lem:basifactsofTi_1} If $h\in H$ and $hT_i\cap T_i\neq\emptyset$, then $h\in H_i$. In particular $H_{\tilde v}\le H_i$ for every vertex $\tilde v\in T_i$.

        \item \label{lem:basifactsofTi_2}  $T_i$ has finitely many $H_i$-orbits of vertices.
    \end{enumerate}
\end{lem}

\begin{proof}
(\ref{lem:basifactsofTi_1}): $H$ permutes the components of $\widehat Y$, so $hT_i$ is a component, and if it meets $T_i$, then $hT_i=T_i$.

(\ref{lem:basifactsofTi_2}): By (\ref{lem:basifactsofTi_1}), two vertices of $T_i$ in the same $H$-orbit are in the same $H_i$-orbit, so $V(T_i)/H_i$ injects into the finite set $V(\Lambda_H)$.
\end{proof}

Using Lemma~\ref{lem:basifactsofTi} and the arguments in the proof of Lemma~\ref{lem:basisalmostfiber}, we have the following remark.

\begin{rem} \label{rem:remarkonHi}
The metric $\widehat d_{\Omega_{T_i}}$ is a finite length metric on $\Omega_{T_i}$, $H_i$ acts on $(\Omega_{T_i},\widehat d_{\Omega_{T_i}})$ by isometries, and for $x_0\in\Omega_{T_i}$:
\begin{itemize}
    \item $\{h\in H_i:\ \widehat d_{\Omega_{T_i}}(x_0,hx_0)\le r\}$ is finite for every $r\ge0$,

    \item $H_i\cdot B(x_0,r_0)=\Omega_{T_i}$ for some $r_0>0$.
\end{itemize}
\end{rem}

The following lemma shows that the metrics $\widehat d_{\Omega_H}$ and $\widehat d_{\Omega_{T_i}}$ on $\Omega_{T_i}$ are comparable.

\begin{lem} \label{lem:Omega_i_undistorted}
For all $x,y\in\Omega_{T_i}$, we have
$$\widehat d_{\Omega_H}(x,y) \le \widehat d_{\Omega_{T_i}}(x,y) \le A\,\widehat d_{\Omega_H}(x,y),$$
where $A$ is the constant of Proposition~\ref{prop:horizontal pieces}.   
\end{lem}

\begin{proof}
The inequality $\widehat d_{\Omega_H}(x,y) \le \widehat d_{\Omega_{T_i}}(x,y)$ follows from the inclusion $\Omega_{T_i}\subset\Omega_H$ and the definitions of the two metrics. It remains to prove the right-hand inequality. 

\begin{figure}[htb]
\centering
\begin{tikzpicture}[scale=0.9,>=stealth,
  wfib/.style={circle,draw,fill=white,inner sep=0pt,minimum size=6pt},
  blk/.style={circle,draw,fill=black,inner sep=0pt,minimum size=5pt},
  hor/.style={line width=1.4pt},
  nhor/.style={line width=0.5pt},
  cont/.style={line width=1.4pt,dotted}]

\fill[gray!30] (0,0.8)--(0.9,3.2)--(5.6,3.4)--(7.0,1.4)--(5.8,-0.6)--(1.0,-0.8)--cycle;
\node at (2.6,2.55) {$\mathcal H(\tilde v_0)$};

\coordinate (P1) at (0.5,-0.82); \coordinate (P2) at (6.3,-0.58);
\coordinate (Q1) at (0.5,-2.02); \coordinate (Q2) at (6.3,-1.78);
\fill[gray!10] (P1)--(P2)--(Q2)--(Q1)--cycle;
\foreach \s in {0.1,0.3,0.5,0.7,0.9} \draw[gray!60,thin] ($(P1)!\s!(P2)$)--($(Q1)!\s!(Q2)$);
\draw[thick] (P1)--(P2) node[right,font=\scriptsize] {$\iota_{\tilde e_1}(\mathcal H(\tilde b_1))$};
\draw[thick] (Q1)--(Q2) node[right,font=\scriptsize] {$\mathcal H(\tilde b_1)$};
\node[font=\scriptsize] at (-0.25,-1.42) {$\mathcal Z(\tilde e_1)$};

\coordinate (R1) at (5.3,3.83); \coordinate (R2) at (7.3,0.97);
\coordinate (S1) at (6.28,4.52); \coordinate (S2) at (8.28,1.66);
\fill[gray!10] (R1)--(R2)--(S2)--(S1)--cycle;
\foreach \s in {0.1,0.3,0.5,0.7,0.9} \draw[gray!60,thin] ($(R1)!\s!(R2)$)--($(S1)!\s!(S2)$);
\draw[thick] (R1)--(R2) node[below,font=\scriptsize] {$\iota_{\tilde e_2}(\mathcal H(\tilde b_2))$};
\draw[thick] (S1)--(S2) node[right,font=\scriptsize] {$\mathcal H(\tilde b_2)$};
\node[font=\scriptsize] at (6.25,4.85) {$\mathcal Z(\tilde e_2)$};

\coordinate (x) at (1.0,1.6); \coordinate (y) at (4.9,0.4);
\coordinate (A1) at (2.0,-0.775); \coordinate (B1) at (3.6,-0.705);
\coordinate (A2) at (5.9,2.97); \coordinate (B2) at (6.7,1.83);
\draw[line width=1pt] (x) .. controls (1.2,0.6) and (1.8,0.0) .. (A1);
\draw[line width=1pt] (A1) .. controls (1.7,-2.9) and (3.9,-3.3) .. (B1);
\draw[line width=1pt] (B1) .. controls (4.4,0.4) and (5.0,2.4) .. (A2);
\draw[line width=1pt] (A2) .. controls (7.8,4.9) and (9.6,2.3) .. (B2);
\draw[line width=1pt] (B2) .. controls (6.2,1.75) and (5.6,0.6) .. (y);
\draw[line width=1pt,dashed] (A1) .. controls (2.4,-0.2) and (3.2,-0.2) .. (B1);
\draw[line width=1pt,dashed] (A2) .. controls (5.75,2.35) and (6.1,1.95) .. (B2);
\node[font=\small] at (2.8,0.0) {$\beta_1$};
\node[font=\small] at (5.55,2.0) {$\beta_2$};
\node[font=\small] at (2.8,-3.2) {$\gamma|_{[a_1,b_1]}$};
\node[font=\small] at (9.1,3.65) {$\gamma|_{[a_2,b_2]}$};
\foreach \p in {x,y,A1,B1,A2,B2} \fill (\p) circle (2pt);
\node[font=\small,left] at (x) {$x$};
\node[font=\small,below] at (y) {$y$};
\node[font=\scriptsize,above left] at (A1) {$\gamma(a_1)$};
\node[font=\scriptsize,above right] at (B1) {$\gamma(b_1)$};
\node[font=\scriptsize,above left] at (A2) {$\gamma(a_2)$};
\node[font=\scriptsize,below right] at (B2) {$\gamma(b_2)$};
\node[font=\small] at (1.3,3.9) {$\Omega_{T_i}$};
\draw[->,gray] (1.3,3.65)--(1.6,3.15);
\begin{scope}[shift={(0.6,-4.6)}]
\fill[gray!22,rounded corners=6pt] (0.3,-0.3) rectangle (5.1,0.3);
\coordinate (v0) at (2.7,0);
\coordinate (c1) at (1.5,0); \coordinate (c2) at (3.9,0);
\coordinate (d1) at (1.8,-1.2); \coordinate (d2) at (3.6,-1.2);
\draw[hor] (0.6,0)--(4.8,0);
\draw[cont] (0.6,0)--(0.1,0); \draw[cont] (4.8,0)--(5.3,0);
\draw[nhor] (v0)--(d1) node[midway,left,font=\scriptsize] {$\tilde e_1$};
\draw[nhor] (v0)--(d2) node[midway,right,font=\scriptsize] {$\tilde e_2$};
\draw[nhor,dotted] (d1)--++(-0.4,-0.6); \draw[nhor,dotted] (d1)--++(0.2,-0.7);
\draw[nhor,dotted] (d2)--++(-0.2,-0.7); \draw[nhor,dotted] (d2)--++(0.4,-0.6);
\node[wfib] at (v0) {}; \node[blk] at (c1) {}; \node[blk] at (c2) {};
\node[blk] at (d1) {}; \node[blk] at (d2) {};
\node[font=\small,above=2pt] at (v0) {$\tilde v_0$};
\node[font=\scriptsize,left] at (d1) {$\tilde b_1$};
\node[font=\scriptsize,right] at (d2) {$\tilde b_2$};
\node[font=\small] at (5.8,0.45) {$T_i$};
\draw[gray,->,thick] (2.62,-0.12) .. controls (2.3,-0.7) and (1.9,-1.6) .. (1.6,-1.95)
  .. controls (1.35,-2.2) and (1.1,-1.9) .. (1.4,-1.6) .. controls (1.8,-1.1) and (2.3,-0.6) .. (2.55,-0.1);
\node[font=\scriptsize,gray!50!black] at (0.4,-1.9) {$\bar\gamma|_{[a_1,b_1]}$};
\end{scope}
\end{tikzpicture}
\caption{The proof of Lemma~\ref{lem:Omega_i_undistorted} near a white vertex $\tilde v_0$ of $T_i$. An excursion $\gamma|_{[a_j,b_j]}$ of $\gamma$ outside $\Omega_{T_i}$ begins and ends in $\mathcal H(\tilde v_0)$ (Claim~1), and it is replaced by a path $\beta_j$ in $\mathcal H(\tilde v_0)$ (dashed).}
\label{fig:excursion}
\end{figure}
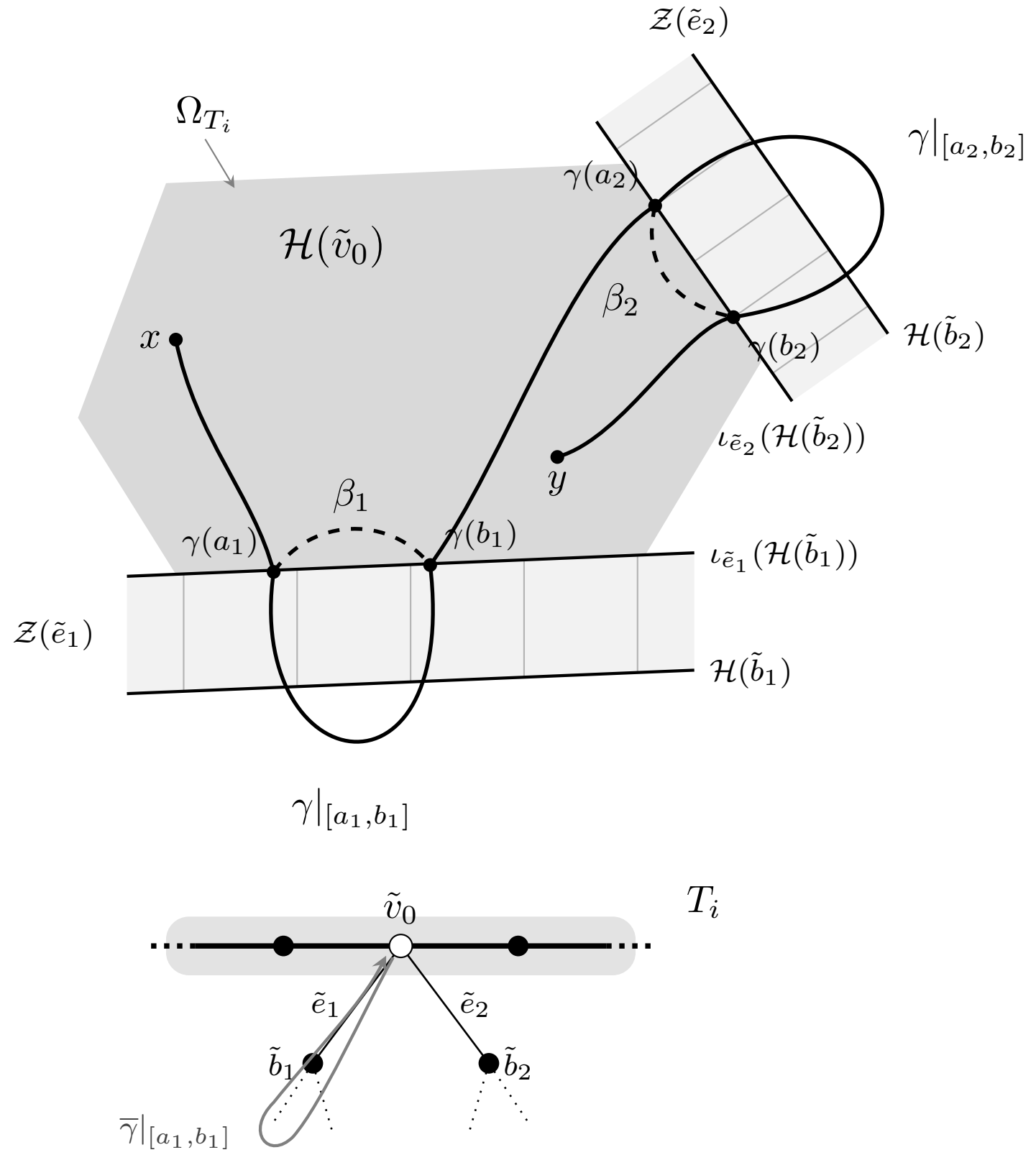

Fix $\epsilon>0$. Since $\widehat d_{\Omega_H}(x,y)$ is the infimum of the $d$-lengths of paths in $\Omega_H$ from $x$ to $y$, there exists a path $\gamma:[0,1]\to\Omega_H$ from $x$ to $y$ such that \[
\ell_d(\gamma)\le\widehat d_{\Omega_H}(x,y)+\epsilon \]

Let $\bar\gamma:=\Pi\circ\gamma \colon [0,1] \to T_H$, where $\Pi \colon \Omega_H \to T_H$ is the map given in Lemma~\ref{lem:basisalmostfiber}. The set $O:=\bar\gamma^{-1}(T_H\setminus T_i)$ is open in $[0,1]$ and contains neither $0$ nor $1$. Moreover, $\gamma([0,1]\setminus O) \subset \Omega_{T_i}$. Write $O =\bigsqcup_{j\in J}(a_j,b_j)$ as a disjoint union of open intervals, where $J$ is countable, and let $\ell_j:=\ell_d\bigl(\gamma|_{[a_j,b_j]}\bigr)$. Since these intervals are pairwise disjoint, \[\sum_{j\in J}\ell_j\le\ell_d(\gamma)<\infty\]

\textit{Claim~1} (see Figure~\ref{fig:excursion}): For each $j \in J$, we have $\gamma(a_j), \gamma(b_j) \in \Omega_{T_i}$.

Let $p:=\bar\gamma(a_j)$ and $q:=\bar\gamma(b_j)$. Both points lie in the tree $T_i$, while $\bar\gamma((a_j,b_j))\cap T_i=\emptyset$.
\begin{enumerate}[(a)]
    \item $p = q$. Indeed, if $p \neq q$, then the geodesic $[p,q]$ in $T_H$ lies in the subtree $T_i$, and every path in a tree from $p$ to $q$ contains $[p,q]$. Thus, $\bar\gamma$ would have to pass through a point of $T_i\setminus\{p,q\}$ at some time in $(a_j,b_j)$, a contradiction.

    \item $p$ is a vertex. Otherwise, $p$ lies in the interior of an edge $e\subset T_i$, and then $\bar\gamma(s)\in e\subset T_i$ for all $s>a_j$ sufficiently close to $a_j$, a contradiction.
\end{enumerate}
Write $\tilde v_0(j):=p$. By Lemma~\ref{lem:basisalmostfiber}(\ref{lem:basisalmostfiber2}), we have $\gamma(a_j),\gamma(b_j)\in\mathcal H(\tilde v_0(j))\subset\Omega_{T_i}$.

\textit{Claim~2}: There exists a path $\gamma'$ in $\Omega_{T_i}$ connecting $x$ to $y$ such that $$\widehat d_{\Omega_{T_i}}(x,y)\le\ell_d(\gamma') \le A\,\ell_d(\gamma)\le A\bigl(\widehat d_{\Omega_H}(x,y)+\epsilon\bigr).$$

By Proposition~\ref{prop:horizontal pieces}(\ref{prop:horizontal pieces_4}), for each $j\in J$ there exists a path $\beta_j\colon[a_j,b_j]\to\mathcal H(\tilde v_0(j))$ from $\gamma(a_j)$ to $\gamma(b_j)$ such that
$$\ell_d(\beta_j)\le A\,d\bigl(\gamma(a_j),\gamma(b_j)\bigr)\le A\,\ell_j.$$ Define $\gamma'\colon[0,1]\to\Omega_{T_i}$ by setting $\gamma':=\gamma$ on $[0,1]\setminus O$ and $\gamma':=\beta_j$ on $[a_j,b_j]$. The map $\gamma'$ is continuous. Indeed, it agrees with $\gamma$ on the closed set $[0,1]\setminus O$, each $\beta_j$ has the same endpoints as $\gamma|_{[a_j,b_j]}$, and the image of $\beta_j$ has $d$-diameter at most $A\ell_j$. Since $\sum_j\ell_j<\infty$, these diameters tend to $0$ along any sequence of distinct intervals $(a_j,b_j)$, which gives continuity at the accumulation points of these intervals. Then $\gamma'$ is a path in $\Omega_{T_i}$ from $x$ to $y$, and
\begin{align*}
\widehat d_{\Omega_{T_i}}(x,y)&\le\ell_d(\gamma')=\ell_d(\gamma)-\sum_{j\in J}\ell_j+\sum_{j\in J}\ell_d(\beta_j)\\
&\le\ell_d(\gamma)+(A-1)\sum_{j\in J}\ell_j\le A\,\ell_d(\gamma)\le A\bigl(\widehat d_{\Omega_H}(x,y)+\epsilon\bigr).
\end{align*} Letting $\epsilon\to0$ completes the proof.
\end{proof}

By Lemma~\ref{lem:basisalmostfiber}, Remark~\ref{rem:remarkonHi} and Lemma~\ref{lem:svarc-milnor}, the orbit maps of $H$ and $H_i$ to $\Omega_H$ and $\Omega_{T_i}$ are quasi-isometries. Together with Lemma~\ref{lem:Omega_i_undistorted}, this gives the following.

\begin{cor} \label{cor:Hiundistorted}
$H_i$ is finitely generated and undistorted in $H$.
\end{cor}

\begin{prop}\label{prop:distortion-Hi-H}
For every component $T_i$ of $\widehat Y$, we have $\Delta_{H_i}^{G}\preceq\Delta_{H}^{G}$.
\end{prop}

\begin{proof}
By Corollary~\ref{cor:Hiundistorted}, there are constants $\lambda\ge1$ and $\epsilon\ge0$ with $|h|_{H_i}\le\lambda\,|h|_H+\epsilon$ for all $h\in H_i$. If $h\in H_i$ and $|h|_G\le n$, then $|h|_H\le\Delta_H^G(n)$, so $|h|_{H_i}\le\lambda\,\Delta_H^G(n)+\epsilon$.
\end{proof}

\section{Slopes and Dilations of subgroups}\label{sec:slopes}
By Proposition~\ref{prop:distortion-Hi-H}, we know that $\Delta_{H_i}^{G} \preceq \Delta_{H}^{G}$ for every component $T_i$ of $\widehat Y$. In this section, we study how $H_i$ sits inside $G$. Once the distortion of each $H_i$ in $G$ is understood, this gives a lower bound for the distortion of $H$ in $G$.

Recall from the proof of Lemma~\ref{lem:black_types} that distinct edges of $T$ at a black vertex lie over distinct edges $e\neq e'$ of $\Delta$, so that $k(e)\neq k(e')$. 

We do not assume that black vertices of $\Delta$ have valence two. Recall that $T_i$ is a component of $\widehat Y$ and $H_i=\mathrm{Stab}_H(T_i)$. By Lemma~\ref{lem:vertices_of_Ti}, every black vertex of $T_i$ has at least two edges in $T_i$, all of them horizontal, and every white vertex of $T_i$ is of fiber type.

Let $w \in \Delta$ be a white vertex adjacent to the black vertex $b \in \Delta$, and let $e = (b, w)$. Let $\tilde e=(\tilde b,\tilde w)$ be an edge in $T$ over $e=(b,w)$. Recall from Definition~\ref{defn:coordinates} that every element $x\in G_{\tilde b}$ has well-defined coordinates $(X,Y)\in\mathbb Z^2$, namely $g^{-1}xg=r_b^Xt_b^Y$, where $\tilde b=g\bar b$ and $G_b=\langle r_b\rangle\times\langle t_b\rangle$. Recall also from Proposition~\ref{prop:split-over-tori} that the inclusion $G_{\tilde b}=G_{\tilde e}\hookrightarrow G_{\tilde w}$ is given by $r_b\mapsto p_e$ and $t_b\mapsto p_e^{k(e)}t_w$, hence we have
$$r_b^{X}t_b^{Y}\ \longmapsto\ p_e^{\,X+k(e)Y}\,t_w^{\,Y}.$$

\begin{defn} \label{defn:horizontal coordinate}
 The \textit{horizontal coordinate along $\tilde e$} of the group element $x \in G_{\tilde b}$ is $X + k(e) Y$. Hence we obtain the following map
$$\chi_{\tilde e}:G_{\tilde b}\to\mathbb Z \quad \text{given by} \,\, \chi_{\tilde e}(x):=X+k(e)\,Y.$$   
\end{defn}

Since $x\mapsto(X,Y)$ is an isomorphism $G_{\tilde b}\to\mathbb Z^2$, the map $\chi_{\tilde e}$ is a homomorphism. If $h\in G$, then $h\tilde b=(hg)\bar b$, and $(hg)^{-1}(hxh^{-1})(hg)=g^{-1}xg$. Hence $hxh^{-1}\in G_{h\tilde b}$ has the same coordinates as $x\in G_{\tilde b}$. Since $h\tilde e$ lies over the same edge $e$ as $\tilde e$, we get
\begin{equation}\label{eq:chi-equivariant}
\chi_{h\tilde e}(hxh^{-1})=\chi_{\tilde e}(x)\qquad\text{for all }h\in G,\ x\in G_{\tilde b}.
\end{equation}

\begin{lem} \label{lem:slope_vertical_extended}
    Let $\tilde e=(\tilde b,\tilde w)$ be an edge of $T$ over $e=(b,w)$. Let $t_{\tilde w}$ be the vertical element given by Definition~\ref{defn:vertical}. Then the following properties hold.
\begin{enumerate}
    \item \label{lem:slope_vertical_2}  $\chi_{\tilde e}(t_{\tilde w})=0$.
    \item  \label{lem:slope_vertical_3} If $\tilde e'=(\tilde b,\tilde w')$ is another edge at $\tilde b$, over $e'$, then $\chi_{\tilde e'}(t_{\tilde w})=k(e')-k(e)$.
\end{enumerate}
\end{lem}

\begin{proof} By Lemma~\ref{lem:slope_vertical}, $t_{\tilde w}$ has coordinates $(-k(e),1)$. Hence $\chi_{\tilde e}(t_{\tilde w})=-k(e)+k(e)=0$ and $\chi_{\tilde e'}(t_{\tilde w})=-k(e)+k(e')$.
\end{proof}

\begin{lem}\label{lem:chi_nonzero}
Let $\tilde b$ be a black vertex of $T_H$ with $H_{\tilde b}$ infinite cyclic, generated by $u$, and let $\tilde e=(\tilde b,\tilde w)$ be an edge at $\tilde b$ with $\tilde w$ of fiber type. Then $\chi_{\tilde e}(u)\neq0$.
\end{lem}

\begin{proof}
Let $(a,\gamma)$ be the coordinates of $u$. If $\chi_{\tilde e}(u)=a+k(e)\gamma=0$, then $(a,\gamma)=\gamma\,(-k(e),1)$, so $u=t_{\tilde w}^{\gamma}$ by Lemma~\ref{lem:slope_vertical}. Thus $u\in H\cap\langle t_{\tilde w}\rangle$. Since $\tilde w$ is of fiber type, $u=1$ by Lemma~\ref{lem:basis_vertical}(\ref{lem:basis_vertical_4}), contradicting that $u$ generates an infinite cyclic group.
\end{proof}

\subsection{Slopes and the dilation homomorphism}
Recall that a black vertex may be incident to more than two edges, so the notion of slope at a black vertex is not well-defined. Instead, we define the slope from one edge to another as follows.

\begin{defn} \label{defn:slope}
Let $\tilde b\in T_i$ be a black vertex, let $\tilde e,\tilde e'$ be two edges of $T_i$ at $\tilde b$, and let $u$ be a generator of $H_{\tilde b}$. The \textit{slope from $\tilde e$ to $\tilde e'$} is
$$\mathfrak{sl}(\tilde e\to\tilde e'):=\left|\frac{\chi_{\tilde e'}(u)}{\chi_{\tilde e}(u)}\right|\in\mathbb Q_{>0},$$
which is defined by Lemma~\ref{lem:chi_nonzero}.
We also set $\mathfrak{sl}(\tilde e\to\tilde e):=1$.
Explicitly, if $u$ has coordinates $(a,\gamma)$, then
$$\mathfrak{sl}(\tilde e\to\tilde e')=\left|\frac{a+k(e')\gamma}{a+k(e)\gamma}\right|.$$
\end{defn}

\begin{lem} \label{lem:slope}
In the setting of Definition~\ref{defn:slope}:
\begin{enumerate}
\item \label{lem:slope_1} $\mathfrak{sl}(\tilde e\to\tilde e')$ does not depend on the choice of the generator $u$.
    \item \label{lem:slope_2} $\mathfrak{sl}(\tilde e'\to\tilde e)=\mathfrak{sl}(\tilde e\to\tilde e')^{-1}$.
    \item \label{lem:slope_3} $\mathfrak{sl}(h\tilde e\to h\tilde e')=\mathfrak{sl}(\tilde e\to\tilde e')$ for all $h\in H_i$.
    \item \label{lem:slope_4} $\mathfrak{sl}(\tilde e\to\tilde e')\,\mathfrak{sl}(\tilde e'\to\tilde e'')=\mathfrak{sl}(\tilde e\to\tilde e'')$ for any three edges $\tilde e,\tilde e',\tilde e''$ of $T_i$ at $\tilde b$.
\end{enumerate}
\end{lem}

\begin{proof}
(\ref{lem:slope_1}): The other generator is $u^{-1}$, and $\chi_{\tilde e}(u^{-1})=-\chi_{\tilde e}(u)$, $\chi_{\tilde e'}(u^{-1})=-\chi_{\tilde e'}(u)$.

(\ref{lem:slope_2}) and (\ref{lem:slope_4}): These parts are immediate from the definition, since all slopes at $\tilde b$ are quotients of the numbers $|\chi_{\tilde e}(u)|$ for the same $u$.

(\ref{lem:slope_3}): If $u$ generates $H_{\tilde b}$, then $huh^{-1}$ generates $H_{h\tilde b}=hH_{\tilde b}h^{-1}$, and by (\ref{eq:chi-equivariant}), $\chi_{h\tilde e}(huh^{-1})=\chi_{\tilde e}(u)$ and $\chi_{h\tilde e'}(huh^{-1})=\chi_{\tilde e'}(u)$.
\end{proof}

\begin{defn} \label{defn:slope_product}
 Let $\tilde w,\tilde w'$ be white vertices in the tree $T_i$. For each black vertex $\tilde b$ on the geodesic $[\tilde w,\tilde w']$, let $\tilde e^-_{\tilde b}$ and $\tilde e^+_{\tilde b}$ be the edges by which the geodesic enters and leaves $\tilde b$. As edges of $T$, both are oriented from $\tilde b$ to a white vertex, and the signs only record the direction of the geodesic. We use this convention for entering and leaving edges throughout the paper. We define the \textit{slope product} from $\tilde w$ to $\tilde w'$ by
$$\mathfrak S(\tilde w,\tilde w'):=\prod_{\tilde b\in[\tilde w,\tilde w']\ \text{black}}\mathfrak{sl}\bigl(\tilde e^-_{\tilde b}\to\tilde e^+_{\tilde b}\bigr) \,\, \text{and}\,\, \mathfrak S(\tilde w,\tilde w):=1.$$
\end{defn}

\begin{lem} \label{lem:basisPi}
Let $\tilde w,\tilde w',\tilde w''$ be white vertices in $T_i$ and let $h\in H_i$. Then the following properties hold.
\begin{enumerate}
    \item \label{lem:basisPi_1} $\mathfrak S(\tilde w',\tilde w)=\mathfrak S(\tilde w,\tilde w')^{-1}$.
    \item \label{lem:basisPi_2} $\mathfrak S(\tilde w,\tilde w'')=\mathfrak S(\tilde w,\tilde w')\,\mathfrak S(\tilde w',\tilde w'')$.
    \item \label{lem:basisPi_3} $\mathfrak S(h\tilde w,h\tilde w')=\mathfrak S(\tilde w,\tilde w')$.
\end{enumerate}
\end{lem}

\begin{proof}
(\ref{lem:basisPi_1}): Reversing the geodesic interchanges $\tilde e^-_{\tilde b}$ and $\tilde e^+_{\tilde b}$, so the claim follows from Lemma~\ref{lem:slope}(\ref{lem:slope_2}).

(\ref{lem:basisPi_3}): The element $h$ maps $[\tilde w,\tilde w']$ onto $[h\tilde w,h\tilde w']$, preserving the order, so the claim follows from Lemma~\ref{lem:slope}(\ref{lem:slope_3}).

(\ref{lem:basisPi_2}): Let $c$ be the median of $\tilde w,\tilde w',\tilde w''$ in $T_i$, so that $[\tilde w,\tilde w']=[\tilde w,c]\cup[c,\tilde w']$, $[\tilde w',\tilde w'']=[\tilde w',c]\cup[c,\tilde w'']$ and $[\tilde w,\tilde w'']=[\tilde w,c]\cup[c,\tilde w'']$.

Every black vertex $\tilde b\neq c$ on these geodesics lies in the interior of one of $[\tilde w,c]$, $[c,\tilde w']$, or $[c,\tilde w'']$. In each geodesic containing it, the vertex is entered and left through the same edges, possibly in the opposite direction. The black vertices of $[c,\tilde w']$ other than $c$ contribute factors to both $\mathfrak S(\tilde w,\tilde w')$ and $\mathfrak S(\tilde w',\tilde w'')$, with the direction reversed in the latter. Hence, their factors cancel by Lemma~\ref{lem:slope}(\ref{lem:slope_2}). 

The factors contributed by the remaining black vertices $\neq c$ agree on both sides of (\ref{lem:basisPi_2}). If $c$ is white, it contributes no factor, and the proof is complete. If $c$ is black, then $c\notin\{\tilde w,\tilde w',\tilde w''\}$, and the edges $\tilde{e},\tilde{e}',  \tilde{e}''$ of $T_i$ at $c$ pointing towards $\tilde w,\tilde w',\tilde w''$ are pairwise distinct. The factors contributed by $c$ are $\mathfrak{sl}(\tilde{e} \to \tilde{e'}) \, \mathfrak{sl}(\tilde{e}' \to \tilde{e}'')$ on the right-hand side and $\mathfrak{sl}(\tilde{e} \to \tilde{e}'')$ on the left-hand side. These agree by Lemma~\ref{lem:slope}(\ref{lem:slope_4}).
\end{proof}

\begin{defn} \label{defn:dilation}
Fix a white vertex $\tilde w_0\in T_i$. The \textit{dilation} of $H_i$ is the following map
$$\omega_i:H_i\to(\mathbb Q_{>0},\times) \quad \text{given by} \quad h \mapsto \mathfrak S(\tilde w_0,\,h\tilde w_0).$$
\end{defn}

We establish some basic properties of the dilation of $H_i$.

\begin{prop} \label{prop:dilation_basis} For each $i$, let $\omega_i$ be the dilation function given by Definition~\ref{defn:dilation}. Then the following properties hold.
 \begin{enumerate}
     \item  \label{prop:dilation_basis_1} $\omega_i$ is a homomorphism.
     \item  \label{prop:dilation_basis_2} $\omega_i$ does not depend on the choice of $\tilde w_0$.
     \item \label{prop:dilation_basis_3} $\omega_i(h)=1$ whenever $h$ fixes a vertex of $T_i$. Hence $\omega_i$ factors through $\pi_1(\Lambda_i)$, where $\Lambda_i=H_i\backslash T_i$ is the quotient graph. More precisely, let $\mathfrak p_i\colon H_i\to\pi_1(\Lambda_i)$ be the natural surjection. 
     Then there is a homomorphism $\bar\omega_i$ such that the following diagram commutes.
\[
\begin{tikzcd}[column sep=large,row sep=large]
H_i \arrow[r,"\omega_i"] \arrow[d,two heads,"\mathfrak p_i"'] & (\mathbb Q_{>0},\times)\\
\pi_1(\Lambda_i) \arrow[ur,"\bar\omega_i"'] &
\end{tikzcd}
\]
In other words, $\omega_i(h)$ only depends on the loop in the finite graph $\Lambda_i$ obtained by projecting a path in $T_i$ from a vertex $\tilde v$ to $h\tilde v$.
     \item  \label{prop:dilation_basis_4} Let $h\in H_i$ act loxodromically on $T_i$, with axis $A_h$ oriented in the direction of translation. For a black vertex $\tilde b\in A_h$, let $\tilde e^-_{\tilde b},\tilde e^+_{\tilde b}$ be the edges of $A_h$ entering and leaving $\tilde b$. Then, for any set $\mathcal F$ of representatives of the $\langle h\rangle$-orbits of black vertices of $A_h$,

$$\omega_i(h)=\prod_{\tilde b\in\mathcal F}\mathfrak{sl}\bigl(\tilde e^-_{\tilde b}\to\tilde e^+_{\tilde b}\bigr).$$ 
 \end{enumerate}   
\end{prop}

\begin{proof}
(\ref{prop:dilation_basis_1}): For $h_1,h_2\in H_i$, the vertex $h_1\tilde w_0$ is white, so Lemma~\ref{lem:basisPi}(\ref{lem:basisPi_2}),(\ref{lem:basisPi_3}) gives
$$\omega_i(h_1h_2)=\mathfrak S(\tilde w_0,h_1\tilde w_0)\,\mathfrak S(h_1\tilde w_0,h_1h_2\tilde w_0)=\omega_i(h_1)\,\omega_i(h_2).$$

(\ref{prop:dilation_basis_2}): For another white vertex $\tilde w_1\in T_i$ and $h\in H_i$, Lemma~\ref{lem:basisPi} gives
$$\mathfrak S(\tilde w_1,h\tilde w_1)=\mathfrak S(\tilde w_1,\tilde w_0)\,\mathfrak S(\tilde w_0,h\tilde w_0)\,\mathfrak S(h\tilde w_0,h\tilde w_1)=\mathfrak S(\tilde w_1,\tilde w_0)\,\omega_i(h)\,\mathfrak S(\tilde w_1,\tilde w_0)^{-1}=\omega_i(h).$$

(\ref{prop:dilation_basis_3}): Suppose that $h\in H_i$ fixes a vertex $\tilde v\in T_i$. If $\tilde v$ is black, then $h\in H_{\tilde v}\le G_{\tilde e}$ for every edge $\tilde e$ at $\tilde v$ (see (\ref{D5})), so $h$ also fixes a white vertex of $T_i$. So we may assume that $\tilde v$ is white. By (\ref{prop:dilation_basis_2}), we may compute $\omega_i$ with base vertex $\tilde v$, and $\omega_i(h)=\mathfrak S(\tilde v,\tilde v)=1$.

Let $N\trianglelefteq H_i$ be the normal subgroup generated by the vertex stabilizers $H_{\tilde v}$, $\tilde v\in T_i$ (these lie in $H_i$ by Lemma~\ref{lem:basifactsofTi}(\ref{lem:basifactsofTi_1})). By what we just proved, $\omega_i$ is trivial on each $H_{\tilde v}$. Since $\mathrm{ker} (\omega_i)$ is a normal subgroup of $H_i$, it contains $N$ so $\omega_i$ is trivial on $N$. The action of $H_i$ on $T_i$ has no inversions, because $T_i$ is bipartite and the action preserves colours. By \cite[\S I.5.4, Corollary~1]{Ser80}, $H_i/N\cong\pi_1(\Lambda_i)$, and the quotient map $H_i\to H_i/N$ is the map $\mathfrak p_i$. Therefore $\omega_i$ induces a homomorphism
$$\bar\omega_i:\pi_1(\Lambda_i)\to(\mathbb Q_{>0},\times)\quad\text{with}\quad \omega_i=\bar\omega_i\circ\mathfrak p_i .$$

(\ref{prop:dilation_basis_4}): Let $h\in H_i$ be loxodromic. $T_i$ is an $\langle h\rangle$-invariant subtree, and every invariant subtree of a loxodromic element contains its axis. So $A_h\subset T_i$. As $A_h$ is a bi-infinite line in a bipartite tree, so it contains a white vertex $\tilde w_1$. By (\ref{prop:dilation_basis_2}),
$$\omega_i(h)=\mathfrak S(\tilde w_1,h\tilde w_1).$$ Since $\tilde w_1\in A_h$, the geodesic $[\tilde w_1,h\tilde w_1]$ is the segment of $A_h$ from $\tilde w_1$ to $h\tilde w_1$, traversed in the direction of translation. Let $\tilde b_1,\dots,\tilde b_m$ be its black vertices in order. Both endpoints are white, so each $\tilde b_j$ is interior, and the edges by which the geodesic enters and leaves $\tilde b_j$ are $\tilde e^-_{\tilde b_j}$ and $\tilde e^+_{\tilde b_j}$. By Definition~\ref{defn:slope_product},
$$\omega_i(h)=\prod_{j=1}^m\mathfrak{sl}\bigl(\tilde e^-_{\tilde b_j}\to\tilde e^+_{\tilde b_j}\bigr).$$

Note that $\{\tilde b_1,\dots,\tilde b_m\}$ is a set of representatives of the $\langle h\rangle$-orbits of black vertices of $A_h$.

Since $h$ preserves $A_h$ and its orientation, $h\tilde e^-_{\tilde b}=\tilde e^-_{h\tilde b}$ and $h\tilde e^+_{\tilde b}=\tilde e^+_{h\tilde b}$. By Lemma~\ref{lem:slope}(\ref{lem:slope_3}), the factor $\mathfrak{sl}(\tilde e^-_{\tilde b}\to\tilde e^+_{\tilde b})$ is the same for $\tilde b$ and $h\tilde b$. So the product does not change when $\{\tilde b_1,\dots,\tilde b_m\}$ is replaced by any other set $\mathcal F$ of representatives.
\end{proof}

The following lemma allows us to view slopes in a different way. 

\begin{lem} \label{lem:key_identity}
Let $\tilde b\in T_i$ be a black vertex, and let $\tilde e=(\tilde b,\tilde w)$ and $\tilde e'=(\tilde b,\tilde w')$ be two distinct edges of $T_i$ at $\tilde b$, lying over $e$ and $e'$. Let $u$ be a generator of the infinite cyclic group $H_{\tilde b}$. Let
$$\kappa:=k(e')-k(e),\quad \alpha:=\chi_{\tilde e'}(u),\quad \beta:=-\chi_{\tilde e}(u).$$ Then
\begin{enumerate}
    \item  \label{lem:key_identity_1}  $\kappa\neq0$, $\alpha\neq0$ and $\beta\neq0$.

    \item \label{lem:key_identity_2}  In $G_{\tilde b}$, we have
$u^{\kappa}=t_{\tilde w}^{\,\alpha}\;t_{\tilde w'}^{\,\beta}.$

\item \label{lem:key_identity_3}  The exponents are unique, that is if $u^{\kappa}=t_{\tilde w}^{\,\alpha_1}t_{\tilde w'}^{\,\beta_1}$ then $\alpha_1=\alpha$ and $\beta_1=\beta$.

\item    \label{lem:key_identity_4}   $\left|\dfrac{\alpha}{\beta}\right|=\mathfrak{sl}(\tilde e\to\tilde e')$.
\end{enumerate}
\end{lem}

We refer the reader to 
Figure~\ref{fig:key} for an illustration of  Lemma~\ref{lem:key_identity} in a concrete case.

\begin{proof}
(\ref{lem:key_identity_1}): Distinct edges of $T$ at $\tilde b$ lie over distinct edges of $\Delta$ at $b$, so $e\neq e'$, and $k(e)\neq k(e')$ by Proposition~\ref{prop:split-over-tori}(\ref{prop:split-over-tori8}). Thus $\kappa\ne0$. The inequalities $\alpha\neq0$, $\beta\neq0$ follow from Lemma~\ref{lem:chi_nonzero}.

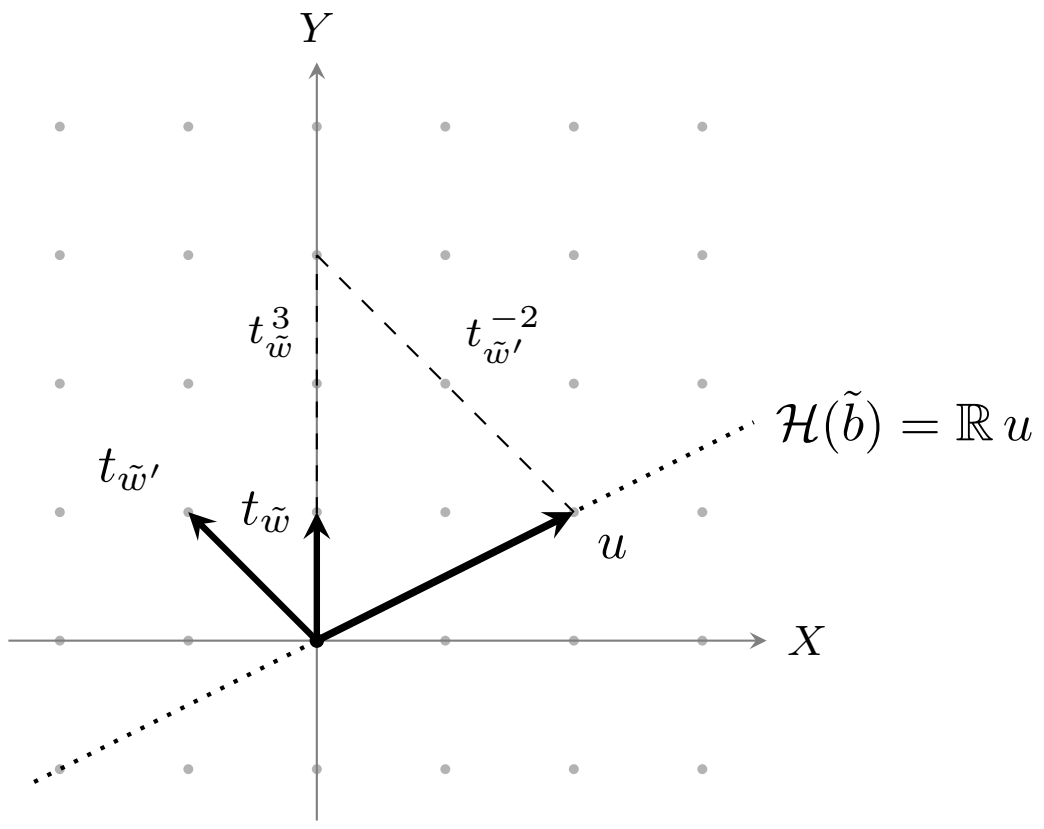
\begin{figure}[htb]
\centering
\begin{tikzpicture}[scale=0.85,>=stealth]
\foreach \x in {-2,...,3} \foreach \y in {-1,...,4} \fill[gray!60] (\x,\y) circle (1.1pt);
\draw[gray,->] (-2.4,0)--(3.5,0) node[right,font=\scriptsize,black] {$X$};
\draw[gray,->] (0,-1.4)--(0,4.5) node[above,font=\scriptsize,black] {$Y$};
\draw[dotted,thick] (-2.2,-1.1)--(3.4,1.7) node[right,font=\small] {$\mathcal H(\tilde b)=\mathbb R\,u$};
\draw[dashed] (0,0)--(0,3) node[pos=0.8,left,font=\scriptsize] {$t_{\tilde w}^{\,3}$};
\draw[dashed] (0,3)--(2,1) node[midway,above right,font=\scriptsize] {$t_{\tilde w'}^{\,-2}$};
\draw[->,line width=1.2pt] (0,0)--(0,1) node[left,font=\small] {$t_{\tilde w}$};
\draw[->,line width=1.2pt] (0,0)--(-1,1) node[above left,font=\small] {$t_{\tilde w'}$};
\draw[->,line width=1.4pt] (0,0)--(2,1) node[below right,font=\small] {$u$};
\fill (0,0) circle (1.6pt);
\end{tikzpicture}
\caption{The key identity in $G_{\tilde b}\cong\mathbb Z^2$ for $k(e)=0$, $k(e')=1$ and $u=(2,1)$. Here $t_{\tilde w}=(0,1)$, $t_{\tilde w'}=(-1,1)$ and $u=t_{\tilde w}^{\,3}\,t_{\tilde w'}^{\,-2}$, so $\mathfrak{sl}(\tilde e\to\tilde e')=3/2$.}
\label{fig:key}
\end{figure}

(\ref{lem:key_identity_2}): Let $(a,\gamma)$ be the coordinates of $u$, and write $k=k(e)$, $k'=k(e')$. Since $G_{\tilde b}\cong\mathbb Z^2$ is abelian, it suffices to compare coordinates.
On the left side,  $u^\kappa$ has coordinates $(\kappa a, \kappa \gamma) = \bigl((k'-k)a,\ (k'-k)\gamma\bigr)$.
On the right side, by Lemma~\ref{lem:slope_vertical}, $t_{\tilde w}$ and $t_{\tilde w'}$ have coordinates $(-k,1)$ and $(-k',1)$. With $\alpha=a+k'\gamma$ and $\beta=-(a+k\gamma)$, we have
    \begin{align*}
         \alpha(-k,1)+\beta(-k',1) &=\bigl(-k(a+k'\gamma)+k'(a+k\gamma),\ (a+k'\gamma)-(a+k\gamma)\bigr) \\
         & =\bigl((k'-k)a,\ (k'-k)\gamma\bigr). 
    \end{align*}
The two sides agree, hence $u^{\kappa}=t_{\tilde w}^{\,\alpha}\;t_{\tilde w'}^{\,\beta}.$

(\ref{lem:key_identity_3}): The vectors $(-k,1)$ and $(-k',1)$ have determinant $k'-k=\kappa\neq0$, so they are linearly independent, and the expression of a vector as an integer combination of them is unique.

(\ref{lem:key_identity_4}): By Definition~\ref{defn:slope}, we have $$\mathfrak{sl}(\tilde e\to\tilde e')=|\chi_{\tilde e'}(u)/\chi_{\tilde e}(u)|=|\alpha/\beta|.$$
\end{proof}

\section{Labelling the axis}\label{sec:labelling}

\begin{conv}\label{conv:fix_h}
Throughout Sections~\ref{sec:labelling}, \ref{sec:rho} and~\ref{sec:overlap}, we fix a component $T_i$ of $\widehat Y$ and an element $h\in H_i$ acting loxodromically on $T_i$. Let $\omega_i\colon H_i\to(\mathbb Q_{>0},\times)$ be the dilation of $H_i$ given by Definition~\ref{defn:dilation}. Since $\omega_i$ is a homomorphism (see Proposition~\ref{prop:dilation_basis}(\ref{prop:dilation_basis_1})), $\omega_i(h^{-1})=\omega_i(h)^{-1}$. Hence, by replacing $h$ by $h^{-1}$ if necessary, we assume that $\omega_i(h)\le1$.
\end{conv}

\subsection{Notation}  \label{subsec:axis_notation}
\begin{enumerate}
    \item Let $A_h\subset T_i$ be the axis of $h$ given by Proposition~\ref{prop:dilation_basis}(\ref{prop:dilation_basis_4}), oriented in the direction of translation. Since $T_i$ is bipartite and $h$ preserves colours, the combinatorial translation length of $h$ (the number of edges by which $h$ translates $A_h$) is even. Now write it as $2m$ with $m\ge1$. Label the vertices of $A_h$ in order as $$\dots,\ \tilde w_{-1},\ \tilde b_0,\ \tilde w_0,\ \tilde b_1,\ \tilde w_1,\ \tilde b_2,\ \tilde w_2,\ \dots$$
where  $\tilde b_j$ is black, $\tilde w_j$ is white, and each $\tilde b_j$ lies between $\tilde w_{j-1}$ and $\tilde w_j$. Then

$$h\tilde b_j=\tilde b_{j+m},\quad h\tilde w_j=\tilde w_{j+m}\quad(j\in\mathbb Z).$$ 

\item For each $j$, let $\tilde e_j^-:=(\tilde b_j,\tilde w_{j-1})$ and $\tilde e_j^+:=(\tilde b_j,\tilde w_j)$ be the edges by which $A_h$ enters and leaves $\tilde b_j$. Let $e_j^-,e_j^+$ be their images in $\Delta$. Since $h\tilde e_j^\pm=\tilde e_{j+m}^\pm$, we have $e_{j+m}^\pm=e_j^\pm$.

\item By Lemma~\ref{lem:vertices_of_Ti}, every $H_{\tilde b_j}$ is infinite cyclic and every $\tilde w_j$ is of fiber type. We fix generators $u_1,\dots,u_m$ of $H_{\tilde b_1},\dots,H_{\tilde b_m}$, and then extend to all $j\in\mathbb Z$ by
$u_{j+m}:=h\,u_j\,h^{-1}.$
Then $u_j$ generates $H_{\tilde b_j}$ for every $j$, because $H_{\tilde b_{j+m}}=hH_{\tilde b_j}h^{-1}$. In particular $u_{nm+1}=h^nu_1h^{-n}$ for $n\ge0$. We denote $t_j:=t_{\tilde w_j}$. By Lemma~\ref{lem:basis_vertical}(\ref{lem:basis_vertical_2}), $t_{j+m}=h\,t_j\,h^{-1}$.
\end{enumerate}

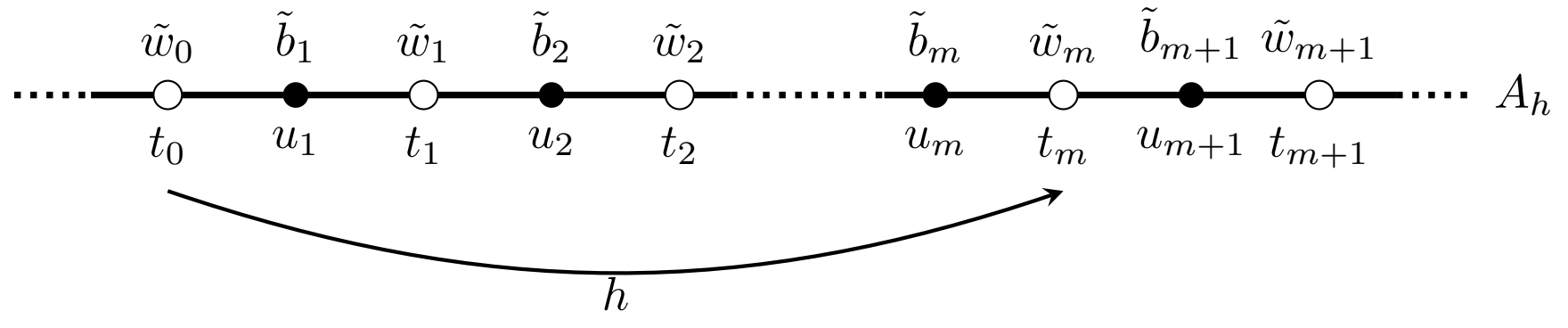
\begin{figure}[htb]
\centering
\begin{tikzpicture}[scale=0.95,>=stealth]
\draw[hor] (-0.6,0)--(4.4,0);
\draw[cont] (4.4,0)--(5.6,0);
\draw[hor] (5.6,0)--(9.6,0);
\draw[cont] (-1.2,0)--(-0.6,0);
\draw[cont] (9.6,0)--(10.2,0);
\foreach \x/\n in {0/0,2/1,4/2,7/m,9/m+1}{
  \node[wfib] at (\x,0) {};
  \node[font=\small,above=3pt] at (\x,0) {$\tilde w_{\n}$};
  \node[font=\small,below=3pt] at (\x,0) {$t_{\n}$};}
\foreach \x/\n in {1/1,3/2,6/m,8/m+1}{
  \node[blk] at (\x,0) {};
  \node[font=\small,above=3pt] at (\x,0) {$\tilde b_{\n}$};
  \node[font=\small,below=3pt] at (\x,0) {$u_{\n}$};}
\draw[->,thick] (0,-0.75) .. controls (2.5,-1.6) and (4.5,-1.6) .. (7,-0.75);
\node[font=\small] at (3.5,-1.55) {$h$};
\node[font=\small] at (10.6,0) {$A_h$};
\end{tikzpicture}
\caption{The axis $A_h\subset T_i$ with the labels of Section~\ref{subsec:axis_notation}. The element $h$ shifts all indices by $m$.}
\label{fig:axis}
\end{figure}

These labels are shown in Figure~\ref{fig:axis}.

Three facts are used constantly:
\begin{itemize}
    \item $u_j\in G_{\tilde b_j}\le G_{\tilde w_{j-1}}\cap G_{\tilde w_j}$ by (\ref{D5}), so $u_j$ fixes $\tilde w_{j-1}$, $\tilde b_j$ and $\tilde w_j$.

    \item $t_{j-1},t_j\in G_{\tilde b_j}$, by Lemma~\ref{lem:slope_vertical}.

    \item $t_j$ commutes with $u_j$ and with $u_{j+1}$, since both lie in $G_{\tilde w_j}$ (see Lemma~\ref{lem:basis_vertical}(\ref{lem:basis_vertical_1})).
\end{itemize}

\begin{lem} \label{lem:periodic_key_identity}
For each $j\in\mathbb Z$ let
$$\kappa_j:=k(e_j^+)-k(e_j^-),\quad \alpha_j:=\chi_{\tilde e_j^+}(u_j),\quad \beta_j:=-\chi_{\tilde e_j^-}(u_j).$$    Then
\begin{enumerate}
    \item  \label{lem:periodic_key_identity_1}   $\kappa_j,\alpha_j,\beta_j$ are nonzero integers, and
$u_j^{\kappa_j}=t_{j-1}^{\,\alpha_j}\;t_j^{\,\beta_j}.$
\item  \label{lem:periodic_key_identity_2}  $\left|\dfrac{\alpha_j}{\beta_j}\right|=\mathfrak{sl}\bigl(\tilde e_j^-\to\tilde e_j^+\bigr)$.
\item  \label{lem:periodic_key_identity_3}  The three sequences are $m$-periodic, that is $\kappa_{j+m}=\kappa_j$, $\alpha_{j+m}=\alpha_j$, $\beta_{j+m}=\beta_j$.

\item  \label{lem:periodic_key_identity_4} $\displaystyle\omega_i(h)=\prod_{j=1}^m\left|\frac{\alpha_j}{\beta_j}\right|$.
\end{enumerate}
\end{lem}

\begin{proof}
    (\ref{lem:periodic_key_identity_1}) and (\ref{lem:periodic_key_identity_2}): Apply Lemma~\ref{lem:key_identity} at the black vertex $\tilde b_j$ with $\tilde e:=\tilde e_j^-$, $\tilde w:=\tilde w_{j-1}$ and $\tilde e':=\tilde e_j^+$, $\tilde w':=\tilde w_j$.

(\ref{lem:periodic_key_identity_3}): $\kappa_{j+m}=\kappa_j$ because $e_{j+m}^\pm=e_j^\pm$. For $\alpha$ and $\beta$: since $u_{j+m}=hu_jh^{-1}$ and $\tilde e^\pm_{j+m}=h\tilde e^\pm_j$, (\ref{eq:chi-equivariant}) gives $\chi_{\tilde e_{j+m}^\pm}(u_{j+m})=\chi_{\tilde e_j^\pm}(u_j)$.

(\ref{lem:periodic_key_identity_4}): 
By Proposition~\ref{prop:dilation_basis}(\ref{prop:dilation_basis_4}), $\omega_i(h)$ is the product of $\mathfrak{sl}(\tilde e^-_{\tilde b}\to\tilde e^+_{\tilde b})$ over representatives of the $\langle h\rangle$-orbits of black vertices of $A_h$. The vertices $\tilde b_1,\dots,\tilde b_m$ are such representatives, and the edges entering and leaving $\tilde b_j$ are $\tilde e_j^-$ and $\tilde e_j^+$. Now use (\ref{lem:periodic_key_identity_2}).
\end{proof}

\begin{lem} \label{lem:consectutive_blacks}
    Let $\tilde w\in T_i$ be a white vertex and let $\tilde b\ne\tilde b'$ be two black vertices of $T_i$ adjacent to $\tilde w$. Then
$H_{\tilde b}\cap H_{\tilde b'}=1 .$
In particular $H_{\tilde b_j}\cap H_{\tilde b_{j+1}}=1$ for all $j$.
\end{lem}
\begin{proof}
    By Proposition~\ref{prop:split-over-tori}(\ref{prop:split-over-tori9}), the group $G_{\tilde b}\cap G_{\tilde b'}$ is cyclic. Let $c$ be a generator of $G_{\tilde b}\cap G_{\tilde b'}$. By Lemma~\ref{lem:slope_vertical}, $t_{\tilde w}$ lies in both $G_{\tilde b}$ and $G_{\tilde b'}$, so $t_{\tilde w}=c^{\,n}$ for some $n\neq0$, because $t_{\tilde w}\neq1$.
Let $x\in H_{\tilde b}\cap H_{\tilde b'}$, and write $x=c^{\,k}$. Then
$x^{\,n}=c^{\,kn}=t_{\tilde w}^{\,k}\in H\cap\langle t_{\tilde w}\rangle .$
Since $\tilde w$ is of fiber type (see Lemma~\ref{lem:vertices_of_Ti}(\ref{lem:vertices_of_Ti_2})), $H\cap\langle t_{\tilde w}\rangle=1$ by Lemma~\ref{lem:basis_vertical}(\ref{lem:basis_vertical_4}). So $t_{\tilde w}^{\,k}=1$, which forces $k=0$ because $t_{\tilde w}$ has infinite order. Hence $x=1$.
\end{proof}

For each $j$, let $w$ be the image of $\tilde w_j$ under the map $T\to\Delta$. Recall from Section~\ref{sec:model} that $p_w\colon G_w\to F_w$ is the projection given by $ft_w^k\mapsto f$, extended to $G_{\tilde w_j}$ via conjugation by $g_{\tilde w_j}$.

\begin{cor} \label{cor:trivialintersectionofpw}
    For each $j$, the elements $p_w(u_j)$ and $p_w(u_{j+1})$ of $F_{w}$ generate cyclic subgroups with trivial intersection.
\end{cor}

\begin{proof}
    Both $u_j$ and $u_{j+1}$ lie in $H_{\tilde w_j}$, and $p_w$ is injective on $H_{\tilde w_j}$ because $\tilde w_j$ is of fiber type. So $\langle p_w(u_j)\rangle\cap\langle p_w(u_{j+1})\rangle=p_w\bigl(\langle u_j\rangle\cap\langle u_{j+1}\rangle\bigr)=1$ by Lemma~\ref{lem:consectutive_blacks}.
\end{proof}

\subsection{Telescoping formula} \label{sec:telescoping}

By Convention~\ref{conv:fix_h} and Lemma~\ref{lem:periodic_key_identity}(\ref{lem:periodic_key_identity_4}), we define
\begin{equation}\label{eq:Omega}
\Omega:=\omega_i(h)^{-1}=\prod_{j=1}^m\left|\frac{\beta_j}{\alpha_j}\right|\ \ge1 .
\end{equation}
Thus $\Omega=1$ when the dilation of $h$ is trivial, and $\Omega>1$ otherwise.

The following lemma is \cite[Lemma~5.2]{HN19}, written in our notation. For the convenience of the reader, and since our indexing is differ from that of \cite{HN19}, we include a proof.

\begin{lem} \label{lem:HN19}
    Let $(\alpha_j)_{j\in\mathbb Z}$ and $(\beta_j)_{j\in\mathbb Z}$ be $m$-periodic sequences of nonzero integers, and assume that $\Omega:=\prod_{j=1}^m|\beta_j/\alpha_j|\ge1$. There is a constant $A_0>0$, depending only on these sequences, such that for every $\mu>0$ there exists a sequence of integers $(s_j)_{j\ge1}$ with:
    \begin{enumerate}
        \item \label{lem:HN19_1} $|s_j|\ge\mu$ for all $j\ge1$.
        \item \label{lem:HN19_2} $|\,\beta_js_j+\alpha_{j+1}s_{j+1}\,|\le A_0$ for all $j\ge1$.
        \item \label{lem:HN19_3} $\displaystyle\sum_{j=1}^{nm}|s_j|\ \succeq\ n^2+\Omega^{\,n}$, where the implicit constants depend on $\mu$.
    \end{enumerate}
\end{lem}

\begin{proof}
Let $A_0:=\max_j|\alpha_j|$ and
\[
\varpi:=\min\Bigl\{\prod_{i=j}^{j+l-1}\Bigl|\frac{\beta_i}{\alpha_{i+1}}\Bigr|\ \Bigm|\ 1\le j\le m,\ 0\le l\le m-1\Bigr\},
\]
where the empty product is $1$. Then $0<\varpi\le1$. Choose an integer $s_1$ with $|s_1|\ge\mu/\varpi$, and define inductively
\[
s_{j+1}:=-\operatorname{sgn}(\beta_js_j)\operatorname{sgn}(\alpha_{j+1})\,\Bigl\lceil\frac{|\beta_js_j|+1}{|\alpha_{j+1}|}\Bigr\rceil .
\]
By induction every $s_j$ is a nonzero integer. Let $x_j:=|\beta_js_j|$ and $Q_j:=|\alpha_{j+1}|\,\bigl\lceil(x_j+1)/|\alpha_{j+1}|\bigr\rceil$. Then $Q_j$ is an integer multiple of $|\alpha_{j+1}|$ with $x_j+1\le Q_j\le x_j+|\alpha_{j+1}|$. Since $\alpha_{j+1}s_{j+1}=-\operatorname{sgn}(\beta_js_j)\,Q_j$, we get
\[
1\ \le\ |\beta_js_j+\alpha_{j+1}s_{j+1}|=Q_j-x_j\ \le\ |\alpha_{j+1}|\ \le\ A_0,
\]
which is (\ref{lem:HN19_2}). Moreover
\begin{equation}\label{eq:sj-recursion}
|s_{j+1}|=\frac{Q_j}{|\alpha_{j+1}|}\ \ge\ \Bigl|\frac{\beta_j}{\alpha_{j+1}}\Bigr|\,|s_j|+\frac1{A_0}.
\end{equation}
Iterating (\ref{eq:sj-recursion}) gives $|s_{j+l}|\ge\bigl(\prod_{i=j}^{j+l-1}|\beta_i/\alpha_{i+1}|\bigr)|s_j|$ for $l\ge0$. For $l=m$, periodicity gives $\prod_{i=j}^{j+m-1}|\beta_i/\alpha_{i+1}|=\prod_{i=1}^m|\beta_i|/\prod_{i=1}^m|\alpha_i|=\Omega$, and keeping the last additive term of (\ref{eq:sj-recursion}) we obtain
\begin{equation}\label{eq:sj-period}
|s_{j+m}|\ \ge\ \Omega\,|s_j|+\frac1{A_0}\ \ge\ |s_j|+\frac1{A_0}\qquad(j\ge1).
\end{equation}

(\ref{lem:HN19_1}): Write $j=1+km+l$ with $k\ge0$ and $0\le l\le m-1$. By (\ref{eq:sj-period}), $|s_{1+km}|\ge|s_1|$, and by periodicity the product $\prod_{i=1+km}^{j-1}|\beta_i/\alpha_{i+1}|$ is at least $\varpi$. Hence $|s_j|\ge\varpi|s_{1+km}|\ge\varpi|s_1|\ge\mu$.

(\ref{lem:HN19_3}): By (\ref{eq:sj-period}) and (\ref{lem:HN19_1}), for $1\le r\le m$ and $k\ge0$ we have $|s_{r+km}|\ge k/A_0$ and $|s_{r+km}|\ge\Omega^k|s_r|\ge\Omega^k\mu$. Therefore
\[
\sum_{j=1}^{nm}|s_j|=\sum_{r=1}^m\sum_{k=0}^{n-1}|s_{r+km}|\ \ge\ \frac{m\,n(n-1)}{2A_0}
\]
and
\[
\sum_{j=1}^{nm}|s_j|\ \ge\ |s_{1+(n-1)m}|\ \ge\ \mu\,\Omega^{\,n-1}.
\]
Hence $$2\sum_{j=1}^{nm}|s_j|\ge\frac{m\,n(n-1)}{2A_0}+\mu\,\Omega^{\,n-1}$$ which gives (\ref{lem:HN19_3}).
\end{proof}

\section{A desirable collection of distorted elements}\label{sec:rho}

Our goal in this section is to construct a collection $\{ \rho_n \}$ of group elements in the subgroup $H$ so that the length of $\rho_n$ in the ambient group $G$ is bounded above by a linear function of $n$ (see Proposition~\ref{prop:Ambient_upper_bound}). In the following sections, we then show that  the length of $\rho_n$ in the subgroup $H_i$ is bounded from below by a quadratic or exponential function of $n$.

Let $u_1, \ldots, u_m$ be the group elements given in part~(3) of Section~\ref{subsec:axis_notation}.

\begin{defn} \label{defn:q_j} Let $\kappa_j,\alpha_j,\beta_j$ be the integers given by Lemma~\ref{lem:periodic_key_identity}, let $\mu$ be the constant fixed in (\ref{eq:choice-mu}) below, and let $(s_j)$ be the sequence given by Lemma~\ref{lem:HN19} for this $\mu$ (Lemma~\ref{lem:HN19} applies by (\ref{eq:Omega})).
    For each $j\ge1$, we define
$$\nu_j:=\kappa_j\,s_j\in\mathbb Z\setminus\{0\}, \quad q_j:=u_1^{\nu_1}\,u_2^{\nu_2}\cdots u_j^{\nu_j}\in H_i,\quad q_0:=1 .$$
For $n\ge1$, let $N:=nm$ and
$$g_n:=q_N=u_1^{\nu_1}\cdots u_N^{\nu_N},\quad \rho_n:=g_n\,u_{N+1}\,g_n^{-1}\in H_i .$$
These elements lie in $H_i$ because each $u_j\in H_{\tilde b_j}\le H_i$ (see Lemma~\ref{lem:basifactsofTi}(\ref{lem:basifactsofTi_1})).
\end{defn}

The constant $\mu$ in (\ref{eq:choice-mu}) is defined in terms of the constants $\theta$, $\Gamma_0$, $c_\ast$, $C_\ast$ of Section~\ref{sec:overlap}. These constants depend only on $K$, $h$ and $u_1,\dots,u_m$ (and on the choice of gates), and not on $\mu$, on the sequence $(s_j)$ or on $n$.

To prove Proposition~\ref{prop:Ambient_upper_bound}, we need several lemmas. Recall from Section~\ref{subsec:axis_notation} that  $t_j:=t_{\tilde w_j}$ and $t_{j+m} =h\,t_j\,h^{-1}$. 
The following lemma gives a formula for the element $q_j$ of Definition~\ref{defn:q_j}.
\begin{lem} \label{lem:Telescoping}
    For each $l\ge1$, let
$c_0:=\alpha_1s_1$ and $c_l:=\beta_ls_l+\alpha_{l+1}s_{l+1}.$
Then for every $j\ge1$ we have
$q_j=t_0^{\,c_0}\,t_1^{\,c_1}\,t_2^{\,c_2}\cdots t_{j-1}^{\,c_{j-1}}\cdot t_j^{\,\beta_js_j}.$
\end{lem}

\begin{proof}
 We prove the statement by induction on $j$. Recall from Definition~\ref{defn:q_j} that $\nu_j:=\kappa_j\,s_j$.

\textit{Base case $j=1$}: By Lemma~\ref{lem:periodic_key_identity}(\ref{lem:periodic_key_identity_1}), $u_1^{\kappa_1}=t_0^{\alpha_1}t_1^{\beta_1}$. The elements $t_0,t_1$ lie in the abelian group $G_{\tilde b_1}$, so raising to the power $s_1$ gives

$$q_1 = u_1^{\nu_1} =u_1^{\kappa_1s_1}=t_0^{\,\alpha_1s_1}\,t_1^{\,\beta_1s_1}=t_0^{\,c_0}\,t_1^{\,\beta_1s_1}.$$

\textit{Inductive step}:  Suppose that the statement holds for $q_j$. In the same way as in the base case, we have $u_{j+1}^{\nu_{j+1}}=t_j^{\,\alpha_{j+1}s_{j+1}}\,t_{j+1}^{\,\beta_{j+1}s_{j+1}}$. Therefore

$$q_{j+1}=q_j\,u_{j+1}^{\nu_{j+1}}=t_0^{\,c_0}\cdots t_{j-1}^{\,c_{j-1}}\cdot\underbrace{t_j^{\,\beta_js_j}\cdot t_j^{\,\alpha_{j+1}s_{j+1}}}_{=\ t_j^{\,c_j}}\cdot t_{j+1}^{\,\beta_{j+1}s_{j+1}}$$
which is the formula for $q_{j+1}$.
\end{proof}

\begin{cor} \label{cor:expression_rho_n}
    For $n\ge1$ and $N=nm$, we define
$P_n:=t_0^{\,c_0}\,t_1^{\,c_1}\cdots t_{N-1}^{\,c_{N-1}}\in G$ where $c_l$ is defined in Lemma~\ref{lem:Telescoping}. Let $g_n$ and $\rho_n$ be given by Definition~\ref{defn:q_j}.
Then the following properties hold.
\begin{enumerate}
    \item \label{cor:expression_rho_n_1}  $g_n=P_n\cdot t_N^{\,\beta_Ns_N}$,
    \item \label{cor:expression_rho_n_2}   $\rho_n=P_n\,u_{N+1}\,P_n^{-1}$,
    \item  \label{cor:expression_rho_n_3}  $|c_l|\le A_\mu:=\max\{A_0,\ |\alpha_1s_1|\}$ for all $0\le l\le N-1$, and $A_\mu$ does not depend on $n$.
\end{enumerate}
\end{cor}

\begin{proof}
    (\ref{cor:expression_rho_n_1}) This is Lemma~\ref{lem:Telescoping} with $j=N$, since $g_n = q_N$.

(\ref{cor:expression_rho_n_2}) By Section~\ref{subsec:axis_notation}, $u_{N+1}\in G_{\tilde b_{N+1}}\le G_{\tilde w_N}$, and $t_N$ is central in $G_{\tilde w_N}$. So $t_N^{\,x}\,u_{N+1}\,t_N^{-x}=u_{N+1}$ for every $x\in\mathbb Z$. Using (\ref{cor:expression_rho_n_1}), we have
\[
\rho_n=P_n\,t_N^{\,\beta_Ns_N}\,u_{N+1}\,t_N^{-\beta_Ns_N}\,P_n^{-1}=P_n\,u_{N+1}\,P_n^{-1}.
\]

(\ref{cor:expression_rho_n_3}) For $1\le l\le N-1$, $|c_l|\le A_0$ by Lemma~\ref{lem:HN19}(\ref{lem:HN19_2}). For $l=0$, $|c_0|=|\alpha_1s_1|$, and the sequence $(s_j)$, hence $s_1$, does not depend on $n$.
\end{proof}

\begin{notation}\label{not:points_xl}
Choose points $x_l\in\mathcal H(\tilde w_l)\subset\widetilde M^{\tilde w_l}$ for $0\le l\le m-1$, and extend to all $l\ge0$ by
$x_{l+m}:=h\,x_l .$
Since $h\,\mathcal H(\tilde w_l)=\mathcal H(h\tilde w_l)=\mathcal H(\tilde w_{l+m})$ (see Proposition~\ref{prop:horizontal pieces}(\ref{prop:horizontal pieces_1})), we have $x_l\in\mathcal H(\tilde w_l)$ for every $l\ge0$. In particular,
$x_{nm}=h^{\,n}x_0$ for each $n \ge 0$, and $x_0\in\mathcal H(\tilde w_0)\subset\Omega_{T_i}$.
We define
$$D_x:=\max_{0\le l\le m-1}d(x_l,x_{l+1}),\quad D_u:=d(x_0,u_1x_0).$$
Because $h$ is an isometry of $(\widetilde M_\varphi,d)$ and $x_{l+m}=hx_l$, the first maximum controls all indices: $$d(x_l,x_{l+1})\le D_x\quad\text{for all }l\ge0 .$$
\end{notation}

\begin{lem} \label{lem:upperboun_translate_x}
    For every $l\ge0$, every $x\in\widetilde M^{\tilde w_l}$ and every $c\in\mathbb Z$, we have
$d\bigl(x,\ t_l^{\,c}x\bigr)\le K\,|c|$
where $K$ is the constant of (\ref{D3}).
\end{lem}

\begin{proof}
    Write $x=(\tilde x,s)\in\widetilde X_w\times\mathbb R$. By Lemma~\ref{lem:basis_vertical}(\ref{lem:basis_vertical_3}), we have $t_l^{\,c}x=(\tilde x,s+c)$. The two points $x$ and $t_l^{\,c}x$ lie on the same vertical fibre, so their distance in the product metric $d_{\tilde w_l}$ is $|c|$. By (\ref{D3}),
$$d\bigl(x,t_l^{\,c}x\bigr)\le K\,d_{\tilde w_l}\bigl(x,t_l^{\,c}x\bigr)\le K|c|.$$
\end{proof}

We need the following lemma in the proof of Proposition~\ref{prop:Ambient_upper_bound}.

\begin{lem} \label{lem:upperbound_x0_xN}
    For every $n\ge1$, let $P_n$ be given by Corollary~\ref{cor:expression_rho_n} and let $N=nm$. Then
$d\bigl(x_0,\ P_nx_N\bigr)\ \le\ N\,(KA_\mu+D_x).$
\end{lem}

\begin{proof}
    For $0\le l\le N$, let

$$P^{(l)}:=t_0^{\,c_0}\,t_1^{\,c_1}\cdots t_{l-1}^{\,c_{l-1}},\quad P^{(0)}:=1,\quad z_l:=P^{(l)}x_l .$$

Then $z_0=x_0$ and $z_N=P_nx_N$. We bound each step $d(z_l,z_{l+1})$ for $0\le l\le N-1$.

Since $P^{(l+1)}=P^{(l)}t_l^{\,c_l}$ and $P^{(l)}$ is an isometry,

$$d(z_l,z_{l+1})=d\bigl(P^{(l)}x_l,\ P^{(l)}t_l^{\,c_l}x_{l+1}\bigr)=d\bigl(x_l,\ t_l^{\,c_l}x_{l+1}\bigr).$$

By the triangle inequality through the point $t_l^{\,c_l}x_l$, and since $t_l^{\,c_l}$ is an isometry, we have
$$d\bigl(x_l,\ t_l^{\,c_l}x_{l+1}\bigr)\le d\bigl(x_l,\ t_l^{\,c_l}x_l\bigr)+d\bigl(t_l^{\,c_l}x_l,\ t_l^{\,c_l}x_{l+1}\bigr)=d\bigl(x_l,\ t_l^{\,c_l}x_l\bigr)+d(x_l,x_{l+1}).$$

The first term is at most $K|c_l|\le KA_\mu$ by Lemma~\ref{lem:upperboun_translate_x} and Corollary~\ref{cor:expression_rho_n}(\ref{cor:expression_rho_n_3}), since $x_l\in\widetilde M^{\tilde w_l}$. The second is at most $D_x$. Hence
$$d(z_l,z_{l+1})\le KA_\mu+D_x$$ for $0 \le l \le N-1$.
Summing over $l$ gives

$$d(x_0,P_nx_N)=d(z_0,z_N)\le\sum_{l=0}^{N-1}d(z_l,z_{l+1})\le N\,(KA_\mu+D_x).$$
\end{proof}

Let $\rho_n$ be the element of $H_i$ given by Definition~\ref{defn:q_j}. The following result shows that the length of $\rho_n$ in the ambient group $G$ is bounded above by a linear function of $n$. 

\begin{prop} \label{prop:Ambient_upper_bound}
 For every $n\ge1$, we have
$$d\bigl(x_0,\ \rho_nx_0\bigr)\ \le\ 2m\,(KA_\mu+D_x)\,n+D_u .$$
Consequently there are constants $C_1,C_2>0$, independent of $n$, such that
$$|\rho_n|_G\le C_1\,n+C_2 .$$   
\end{prop}

\begin{proof}
    Write $u:=u_{N+1}$, so $\rho_n=P_n\,u\,P_n^{-1}$. Insert the two intermediate points $P_nx_N$ and $P_n\,u\,x_N$:

$$d(x_0,\rho_nx_0)\le\underbrace{d\bigl(x_0,\ P_nx_N\bigr)}_{(\mathrm I)}+\underbrace{d\bigl(P_nx_N,\ P_n\,u\,x_N\bigr)}_{(\mathrm{II})}+\underbrace{d\bigl(P_n\,u\,x_N,\ P_n\,u\,P_n^{-1}x_0\bigr)}_{(\mathrm{III})}.$$
We now estimate each of the terms (I), (II) and (III).

For (I): (I) is at most $N(KA_\mu+D_x)$ by Lemma~\ref{lem:upperbound_x0_xN}, using Corollary~\ref{cor:expression_rho_n}(\ref{cor:expression_rho_n_2}).

For (II): since $P_n$ is an isometry, $(\mathrm{II})=d(x_N,\ u\,x_N)$. By Section~\ref{subsec:axis_notation} and Notation~\ref{not:points_xl}, $u=u_{nm+1}=h^{\,n}u_1h^{-n}$ and $x_N=h^{\,n}x_0$. Hence

$$d(x_N,u\,x_N)=d\bigl(h^{\,n}x_0,\ h^{\,n}u_1x_0\bigr)=d(x_0,u_1x_0)=D_u .$$

For (III): applying the isometry $(P_nu)^{-1}$ and then $P_n$,

$$(\mathrm{III})=d\bigl(x_N,\ P_n^{-1}x_0\bigr)=d\bigl(P_nx_N,\ x_0\bigr)\le N(KA_\mu+D_x),$$
again by Lemma~\ref{lem:upperbound_x0_xN}.

Adding the three terms and using $N=nm$ gives $$d\bigl(x_0,\ \rho_nx_0\bigr)\ \le\ 2m\,(KA_\mu+D_x)\,n+D_u .$$

By Lemma~\ref{lem:distortionuniversal}, $|g|_G\le\lambda_G\,d(x_0,gx_0)+\lambda_G\epsilon_G$ for all $g\in G$. Apply this to $g=\rho_n$, and take
$$C_1:=2m\lambda_G(KA_\mu+D_x),\quad C_2:=\lambda_G(D_u+\epsilon_G).$$
\end{proof}

\subsection{The tree geodesic}
We keep $N=nm$, $q_j=u_1^{\nu_1}\cdots u_j^{\nu_j}$, $g_n=q_N$ and $\rho_n=g_n\,u_{N+1}\,g_n^{-1}$ as in Definition~\ref{defn:q_j}.
Define
$$W_j:=q_j\,\tilde w_j\quad(0\le j\le N),\qquad B_j:=q_{j-1}\,\tilde b_j\quad(1\le j\le N) .$$
These are vertices of $T_i$, because $q_j\in H_i$ and $H_i$ preserves $T_i$. Note $W_0=\tilde w_0$.

\begin{lem}\label{lem:adjacency}
    For $1\le j\le N$:
    \begin{enumerate}
        \item \label{lem:adjacency_1} $B_j=q_{j-1}\tilde b_j=q_j\tilde b_j$,
        \item \label{lem:adjacency_2} $W_{j-1}=q_{j-1}\tilde w_{j-1}=q_j\tilde w_{j-1}$,
        \item \label{lem:adjacency_3} $W_j=q_j\tilde w_j=q_{j-1}\tilde w_j$,
        \item \label{lem:adjacency_4} $B_j$ is adjacent to $W_{j-1}$ by the edge $q_{j-1}\tilde e_j^-=q_j\tilde e_j^-$, and to $W_j$ by the edge $q_{j-1}\tilde e_j^+=q_j\tilde e_j^+$.
    \end{enumerate}
\end{lem}
\begin{proof}
    By Section~\ref{subsec:axis_notation}, $u_j$ fixes $\tilde w_{j-1}$, $\tilde b_j$ and $\tilde w_j$, hence also the edges $\tilde e_j^-$ and $\tilde e_j^+$ between them. Since $q_j=q_{j-1}u_j^{\nu_j}$, applying $q_j$ or $q_{j-1}$ to any of these five objects gives the same result. This proves (1), (2), (3) and the equalities of edges in (4). The adjacency in (4) follows by applying $q_{j-1}$ to the edges $\tilde e_j^-=(\tilde b_j,\tilde w_{j-1})$ and $\tilde e_j^+=(\tilde b_j,\tilde w_j)$.
\end{proof}

\begin{lem} \label{lem:sequence_geodesic} The sequence
$$W_0,\ B_1,\ W_1,\ B_2,\ W_2,\ \dots,\ B_N,\ W_N,\ \ \rho_nB_N,\ \rho_nW_{N-1},\ \dots,\ \rho_nB_1,\ \rho_nW_0$$
is a geodesic in $T_i$ from $\tilde w_0$ to $\rho_n\tilde w_0$, with $4N$ edges.
\end{lem}

\begin{proof}
By Lemma~\ref{lem:adjacency}(\ref{lem:adjacency_4}), consecutive vertices of the sequence $W_0,\allowbreak B_1,\allowbreak W_1,\dots,\allowbreak B_N,\allowbreak W_N$ are adjacent. This path has no backtracking. Indeed, by Lemma~\ref{lem:adjacency}, the two neighbours $W_{j-1}=q_j\tilde w_{j-1}$ and $W_j=q_j\tilde w_j$ of $B_j$ are distinct, and the two neighbours $B_j=q_j\tilde b_j$ and $B_{j+1}=q_j\tilde b_{j+1}$ of $W_j$ are distinct. Applying the isometry $\rho_n$, the same holds for the reversed sequence $\rho_nW_N,\allowbreak\rho_nB_N,\dots,\allowbreak\rho_nB_1,\allowbreak\rho_nW_0$.

For the rest of the proof, we show that $\rho_nB_N\neq B_N$, so that there is no backtracking at the junction $W_N$. Indeed, since $u_{N+1}\in G_{\tilde b_{N+1}}\le G_{\tilde w_N}$, we have
\[
\rho_nW_N=q_Nu_{N+1}q_N^{-1}q_N\tilde w_N=q_Nu_{N+1}\tilde w_N=q_N\tilde w_N=W_N .
\] Moreover $B_N=q_N\tilde b_N$ by Lemma~\ref{lem:adjacency}(\ref{lem:adjacency_1}), so $\rho_nB_N=q_Nu_{N+1}\tilde b_N$. Hence $\rho_nB_N=B_N$ if and only if $u_{N+1}\in H_{\tilde b_N}$. But $1\neq u_{N+1}\in H_{\tilde b_{N+1}}$ and $H_{\tilde b_N}\cap H_{\tilde b_{N+1}}=1$ by Lemma~\ref{lem:consectutive_blacks}. So $\rho_nB_N\neq B_N$, and the whole sequence is an edge path without backtracking from $W_0=\tilde w_0$ to $\rho_nW_0=\rho_n\tilde w_0$, with $4N$ edges. A path without backtracking in a tree is a geodesic.
\end{proof}

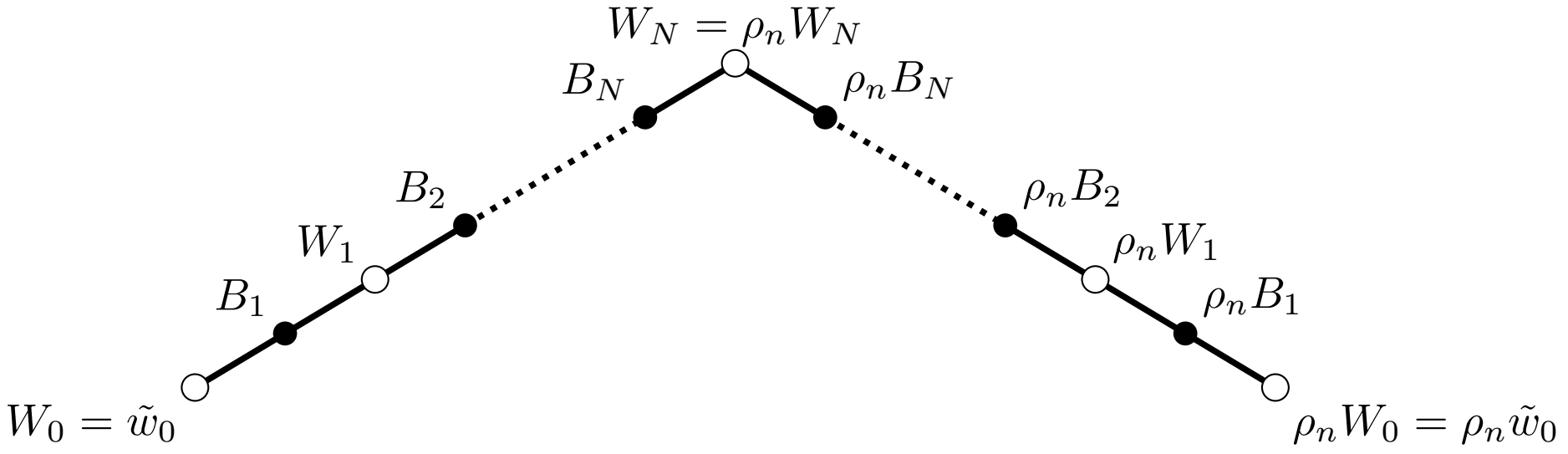
\begin{figure}[htb]
\centering
\begin{tikzpicture}[scale=0.95]
\coordinate (W0) at (0,0);  \coordinate (B1) at (0.75,0.45); \coordinate (W1) at (1.5,0.9);
\coordinate (B2) at (2.25,1.35); \coordinate (BN) at (3.75,2.25); \coordinate (WN) at (4.5,2.7);
\coordinate (rBN) at (5.25,2.25); \coordinate (rB2) at (6.75,1.35); \coordinate (rW1) at (7.5,0.9);
\coordinate (rB1) at (8.25,0.45); \coordinate (rW0) at (9,0);
\draw[hor] (W0)--(B1)--(W1)--(B2);
\draw[cont] (B2)--(BN);
\draw[hor] (BN)--(WN)--(rBN);
\draw[cont] (rBN)--(rB2);
\draw[hor] (rB2)--(rW1)--(rB1)--(rW0);
\foreach \v in {W0,W1,WN,rW1,rW0} \node[wfib] at (\v) {};
\foreach \v in {B1,B2,BN,rBN,rB2,rB1} \node[blk] at (\v) {};
\node[font=\small,below left] at (W0) {$W_0=\tilde w_0$};
\node[font=\small,above left] at (B1) {$B_1$};
\node[font=\small,above left] at (W1) {$W_1$};
\node[font=\small,above left] at (B2) {$B_2$};
\node[font=\small,above left] at (BN) {$B_N$};
\node[font=\small,above] at (WN) {$W_N=\rho_nW_N$};
\node[font=\small,above right] at (rBN) {$\rho_nB_N$};
\node[font=\small,above right] at (rB2) {$\rho_nB_2$};
\node[font=\small,above right] at (rW1) {$\rho_nW_1$};
\node[font=\small,above right] at (rB1) {$\rho_nB_1$};
\node[font=\small,below right] at (rW0) {$\rho_nW_0=\rho_n\tilde w_0$};
\end{tikzpicture}
\caption{The geodesic of Lemma~\ref{lem:sequence_geodesic} from $\tilde w_0$ to $\rho_n\tilde w_0$, with $W_j=q_j\tilde w_j$ and $B_j=q_{j-1}\tilde b_j$. It turns at $W_N$, which is fixed by $\rho_n$.}
\label{fig:geodesic}
\end{figure}

Recall from Notation~\ref{not:points_xl} that $x_0\in\mathcal H(\tilde w_0)\subset\Omega_{T_i}$. Then $\rho_nx_0\in\mathcal H(\rho_n\tilde w_0)\subset\Omega_{T_i}$.

\begin{lem} \label{lem:supremum}
    Let $\gamma:[0,1]\to\Omega_{T_i}$ be a path from $x_0$ to $\rho_nx_0$. For $1\le j\le N$ let
$$\tau_j:=\sup\{\,s\in[0,1]:\ \gamma(s)\in\mathcal H(B_j)\,\},\quad y_j:=\gamma(\tau_j).$$

Then: 
\begin{enumerate}
    \item \label{lem:supremum_1} the sets in the definition of $\tau_j$ are nonempty, and $y_j\in\mathcal H(B_j)$,
    \item \label{lem:supremum_2} $\tau_1<\tau_2<\dots<\tau_N$,
    \item \label{lem:supremum_3} $\displaystyle\ell_d(\gamma)\ \ge\ \sum_{j=1}^{N-1}d(y_j,y_{j+1})$.
\end{enumerate}
\end{lem}

\begin{proof}
    Let $\bar\gamma:=\Pi\circ\gamma:[0,1]\to T_i$ be a continuous path from $\Pi(x_0)=W_0$ to $\Pi(\rho_nx_0)=\rho_nW_0$. By Lemma~\ref{lem:basisalmostfiber}(\ref{lem:basisalmostfiber2}), for every vertex $\tilde v$ of $T_H$ and every $s\in[0,1]$,
\begin{equation}\label{eq:gamma-in-piece}
\gamma(s)\in\mathcal{H}(\tilde v)\iff\bar\gamma(s)=\tilde v .
\end{equation}

(\ref{lem:supremum_1}) The set $\bar\gamma([0,1])$ is connected and contains $W_0$ and $\rho_nW_0$, so it contains the geodesic of Lemma~\ref{lem:sequence_geodesic}, in particular it contains each $B_j$, so by (\ref{eq:gamma-in-piece}) the path $\gamma$ meets $\mathcal H(B_j)$. The set $\{s:\bar\gamma(s)=B_j\}$ is closed, so the supremum $\tau_j$ is attained, and $y_j=\gamma(\tau_j)\in\mathcal H(B_j)$.

(\ref{lem:supremum_2}) Fix $1\le j\le N-1$. The set $\bar\gamma([\tau_j,1])$ is connected and contains $\bar\gamma(\tau_j)=B_j$ and $\bar\gamma(1)=\rho_nW_0$. So it contains the geodesic from $B_j$ to $\rho_nW_0$, which is a subpath of the geodesic of Lemma~\ref{lem:sequence_geodesic} and therefore passes through $B_{j+1}$. Thus there is $s'\in[\tau_j,1]$ with $\bar\gamma(s')=B_{j+1}$. Since $\bar\gamma(\tau_j)=B_j\neq B_{j+1}$, we have $s'>\tau_j$. By the definition of $\tau_{j+1}$ as a supremum, $\tau_{j+1}\ge s'>\tau_j$.

(\ref{lem:supremum_3}) By (\ref{lem:supremum_2}), the intervals $[\tau_j,\tau_{j+1}]$, $1\le j\le N-1$, have disjoint interiors. Length is additive, and the length of a path is at least the distance between its endpoints:

$$\ell_d(\gamma)\ \ge\ \sum_{j=1}^{N-1}\ell_d\bigl(\gamma|_{[\tau_j,\tau_{j+1}]}\bigr)\ \ge\ \sum_{j=1}^{N-1}d\bigl(\gamma(\tau_j),\gamma(\tau_{j+1})\bigr)=\sum_{j=1}^{N-1}d(y_j,y_{j+1}).$$
\end{proof}

The points $y_j$ and $y_{j+1}$ lie on the black lines $\mathcal H(B_j)$ and $\mathcal H(B_{j+1})$ (each $\mathcal H(B_j)$ is a line, being the convex hull of an orbit of the infinite cyclic group $H_{B_j}$ acting by translations), and both $B_j$ and $B_{j+1}$ are adjacent to the white vertex $W_j$. The next lemma says that this configuration is the $q_j$-translate of a fixed reference configuration at $\tilde w_j$.

\begin{lem} \label{lem:WjBj}
    For $1\le j\le N-1$:
    \begin{enumerate}
        \item \label{lem:WjBj_1} $\mathcal H(W_j)=q_j\,\mathcal H(\tilde w_j)$, $\ \mathcal H(B_j)=q_j\,\mathcal H(\tilde b_j)$ and $\ \mathcal H(B_{j+1})=q_j\,\mathcal H(\tilde b_{j+1})$,
        \item \label{lem:WjBj_2} the edge from $B_j$ to $W_j$ is $q_j\tilde e_j^+$, and the edge from $B_{j+1}$ to $W_j$ is $q_j\tilde e_{j+1}^-$,
        \item \label{lem:WjBj_3} the corresponding attaching maps satisfy

$$\iota_{q_j\tilde e_j^+}\circ q_j=q_j\circ\iota_{\tilde e_j^+},\quad \iota_{q_j\tilde e_{j+1}^-}\circ q_j=q_j\circ\iota_{\tilde e_{j+1}^-},$$
\item \label{lem:WjBj_4} on the other side of $B_j$: $\mathcal H(B_j)=q_{j-1}\mathcal H(\tilde b_j)$ and the edge from $B_j$ to $W_{j-1}$ is $q_{j-1}\tilde e_j^-$.
    \end{enumerate}
\end{lem}

\begin{proof}
(\ref{lem:WjBj_1}) and (\ref{lem:WjBj_4}) follow from Lemma~\ref{lem:adjacency}(\ref{lem:adjacency_1}),(\ref{lem:adjacency_3}) and Proposition~\ref{prop:horizontal pieces}(\ref{prop:horizontal pieces_1}). For $B_{j+1}$, note that $B_{j+1}=q_j\tilde b_{j+1}$ by definition. (\ref{lem:WjBj_2}) is Lemma~\ref{lem:adjacency}(\ref{lem:adjacency_4}), applied to $j$ and to $j+1$. (\ref{lem:WjBj_3}) is the equivariance $g\circ\iota_{\tilde e}=\iota_{g\tilde e}\circ g$ of (\ref{D1}).
\end{proof}

\section{Lower bound on distortion of horizontal components}\label{sec:overlap}
Recall from Definition~\ref{defn:horizontal_component} that
the connected components $T_i$ of $\widehat Y$ (see Definition~\ref{defn:hatYandOmegaTi}) are called the \textit{horizontal components} of $H$. A white vertex $\tilde w$ of $T_i$ is an \textit{essential piece} of $T_i$ if it has at least two neighbours in $T_i$. We say that $T_i$ \textit{consists of one piece} if $T_i$ has at most one essential piece.

In this section, we are going to prove the following theorem.

\begin{thm} \label{thm:lower_bound}
   Let $T_i$ be a component of $\widehat Y$  and $H_i=\mathrm{Stab}_H(T_i)$. Suppose $H_i$ contains an element $h$ acting loxodromically on $T_i$. Then:
    \begin{enumerate}
        \item  \label{thm:lower_bound_1}  $\Delta_{H_i}^G(x)\succeq x^2$,
        \item   \label{thm:lower_bound_2}  if $\omega_i(h)\neq1$, then $\Delta_{H_i}^G(x)\succeq 2^{x}$, that is, the distortion is at least exponential.
    \end{enumerate}
\end{thm}

\medskip
\noindent\textit{Idea of the proof.} Using the key identity of Lemma~\ref{lem:periodic_key_identity} and the exponents of Lemma~\ref{lem:HN19}, we constructed in Section~\ref{sec:rho} elements $\rho_n\in H_i$ whose length in $G$ is at most linear in $n$ (see Proposition~\ref{prop:Ambient_upper_bound}). We show that their length in $H_i$ is at least quadratic in $n$, and at least exponential when $\omega_i(h)\neq1$. This is done in the subspace $\Omega_{T_i}$ of the model space, using the bounded overlap of axes in the trees of the white vertex spaces.

For a white vertex $\tilde w\in T$ over $w$, write $X_{\tilde w}$ for the tree factor of $\widetilde M^{\tilde w}=X_{\tilde w}\times\mathbb R$. It is a copy of $\widetilde X_w$, with every edge of length $1$. Let $\pi_X:\widetilde M^{\tilde w}\to X_{\tilde w}$ be the projection. 

Every $g\in G$ maps $\widetilde M^{\tilde w}$ to $\widetilde M^{g\tilde w}$ preserving the product structure, so it induces an isometry of trees

\begin{equation}\label{eq:bar_g}
\bar g:X_{\tilde w}\to X_{g\tilde w},\quad \pi_X\circ g=\bar g\circ\pi_X .
\end{equation}

For $g\in G_{\tilde w}$ this is the map $\bar g=p_w(g)$ of (\ref{D4}). In particular $G_{\tilde w}$ acts on $X_{\tilde w}$ through $p_w$, and this action of $F_w$ is free, because it is the deck action on a universal cover.

Recall the standard facts about an isometry $p$ of a tree $X$ that fixes no point. There is a unique line $\mathrm{Ax}(p)\subset X$, the axis, on which $p$ acts as a translation of length $\|p\|>0$. Moreover $\mathrm{Ax}(p^{k})=\mathrm{Ax}(p)$ for $k\neq0$, and $\mathrm{Ax}(gpg^{-1})=g\,\mathrm{Ax}(p)$, $\|gpg^{-1}\|=\|p\|$ for any isometry $g$.

\begin{lem} \label{lem:axisintersection}
    Let $F$ be a finitely generated group acting freely by isometries on a simplicial tree $X$, and let $p,p'\in F$ be nontrivial. If
$$\operatorname{diam}\bigl(\mathrm{Ax}(p)\cap\mathrm{Ax}(p')\bigr)\ \ge\ \|p\|+\|p'\| ,$$
then $\langle p\rangle\cap\langle p'\rangle\neq1$.
Equivalently, if $\langle p\rangle\cap\langle p'\rangle=1$, then $\operatorname{diam}(\mathrm{Ax}(p)\cap\mathrm{Ax}(p'))<\|p\|+\|p'\|$.
\end{lem} 
\begin{proof}
Let $a:=\|p\|$ and $a':=\|p'\|$. Replacing $p'$ by $p'^{-1}$ if necessary, we may assume that $p$ and $p'$ translate $\mathrm{Ax}(p)\cap\mathrm{Ax}(p')$ in the same direction. Let $\sigma\colon[0,a+a']\to\mathrm{Ax}(p)\cap\mathrm{Ax}(p')$ be an isometric parametrization in this direction. Then $pp'\sigma(0)=\sigma(a+a')=p'p\,\sigma(0)$, so the commutator of $p$ and $p'$ fixes $\sigma(0)$. Since the action is free, $p$ and $p'$ commute. Hence they lie in a common cyclic subgroup $\langle c\rangle$ of the free group $F$, say $p=c^k$ and $p'=c^{k'}$ with $k,k'\neq0$, and $p^{k'}=p'^{k}$ is a nontrivial element of $\langle p\rangle\cap\langle p'\rangle$.
\end{proof}

\begin{lem} \label{lem:boundedgates}
    Let $L,L'$ be two geodesic lines in a tree $X$, and let $D_0\ge0$ with
$\operatorname{diam}(L\cap L')\le D_0$
(the intersection may be empty). Then there are points $o\in L$ and $o'\in L'$ such that for all $p\in L$ and $p'\in L'$,
$$d(p,p')\ \ge\ d(p,o)+d(o',p')-2D_0 .$$
We call $o$ and $o'$ gates of the pair $(L,L')$.
\end{lem} 

\begin{figure}[htb]
\centering
\begin{tikzpicture}[scale=0.95,>=stealth]
\begin{scope}
\draw[line width=7pt,gray!25,line cap=round,line join=round] (-2.4,1.2)--(-0.8,0)--(0.8,0)--(2.4,1.2);
\draw[gray] (-1.6,0.6)--(-1.9,0.05); \draw[gray,dotted] (-1.9,0.05)--(-2.05,-0.25);
\draw[gray] (1.6,-0.6)--(1.9,-0.05); \draw[gray,dotted] (1.9,-0.05)--(2.05,0.25);
\draw[gray] (-0.4,0)--(-0.6,-0.7); \draw[gray,dotted] (-0.6,-0.7)--(-0.68,-1.0);
\draw[thick] (-3.2,1.8)--(-0.8,0)--(0.8,0)--(3.2,-1.8);
\draw[thick,dotted] (-3.2,1.8)--(-3.5,2.03); \draw[thick,dotted] (3.2,-1.8)--(3.5,-2.03);
\draw[thick,dashed] (-3.2,-1.8)--(-0.8,0);
\draw[thick,dashed] (0.8,0)--(3.2,1.8);
\draw[thick,dotted] (-3.2,-1.8)--(-3.5,-2.03); \draw[thick,dotted] (3.2,1.8)--(3.5,2.03);
\draw[line width=2.6pt] (-0.8,0)--(0.8,0);
\draw[->,thick] (1.95,-1.12)--(2.75,-1.72); \node[font=\scriptsize,above right] at (2.35,-1.42) {$\bar u_j$};
\draw[->,thick] (-2.95,-1.86)--(-2.15,-1.26); \node[font=\scriptsize,above left] at (-2.55,-1.56) {$\bar u_{j+1}$};
\node[font=\small,anchor=east] at (-3.5,2.03) {$\ell_j^+$};
\node[font=\small,anchor=east] at (-3.5,-2.03) {$\ell_{j+1}^-$};
\fill (0,0) circle (1.9pt);
\node[font=\scriptsize,above=3pt] at (0,0) {$o_j^+=o_{j+1}^-$};
\fill (-2.4,1.2) circle (1.9pt) node[above right,font=\small] {$p$};
\fill (2.4,1.2) circle (1.9pt) node[above left,font=\small] {$p'$};
\draw[decorate,decoration={brace,mirror,amplitude=4pt}] (-0.8,-0.2)--(0.8,-0.2);
\node[font=\scriptsize,align=center] at (0,-2.05) {$\operatorname{diam}\bigl(\ell_j^+\cap\ell_{j+1}^-\bigr)<D_0$};
\draw[gray,->] (0,-1.85)--(0,-0.45);
\node[font=\small] at (0,2.45) {$X_{\tilde w_j}$};
\end{scope}
\draw[->,thick] (3.9,0.2)--(5.3,0.2) node[midway,above,font=\small] {$\bar h$};
\begin{scope}[shift={(8.2,0)},scale=0.7]
\draw[thick] (-3.2,1.8)--(-0.8,0)--(0.8,0)--(3.2,-1.8);
\draw[thick,dashed] (-3.2,-1.8)--(-0.8,0);
\draw[thick,dashed] (0.8,0)--(3.2,1.8);
\draw[line width=2.2pt] (-0.8,0)--(0.8,0);
\node[font=\small,anchor=east] at (-3.2,1.8) {$\ell_{j+m}^+$};
\node[font=\small,anchor=east] at (-3.2,-1.8) {$\ell_{j+m+1}^-$};
\fill (0,0) circle (2.2pt);
\node[font=\scriptsize,above=3pt] at (0,0) {$o_{j+m}^\pm=\bar h(o_j^\pm)$};
\node[font=\small] at (0,3.5) {$X_{\tilde w_{j+m}}$};
\end{scope}
\end{tikzpicture}
\caption{Lemma~\ref{lem:bounded_axis}. Left: the axes $\ell_j^+$ (solid) and $\ell_{j+1}^-$ (dashed) in $X_{\tilde w_j}$ overlap in a segment of diameter less than $D_0$, which contains the gates. The geodesic from $p$ to $p'$ (shaded) passes through the overlap. If the two axes are disjoint, the gates are the endpoints of the bridge between them. Right: $\bar h$ maps this configuration to the one of index $j+m$.}
\label{fig:lines-in-tree}
\end{figure}
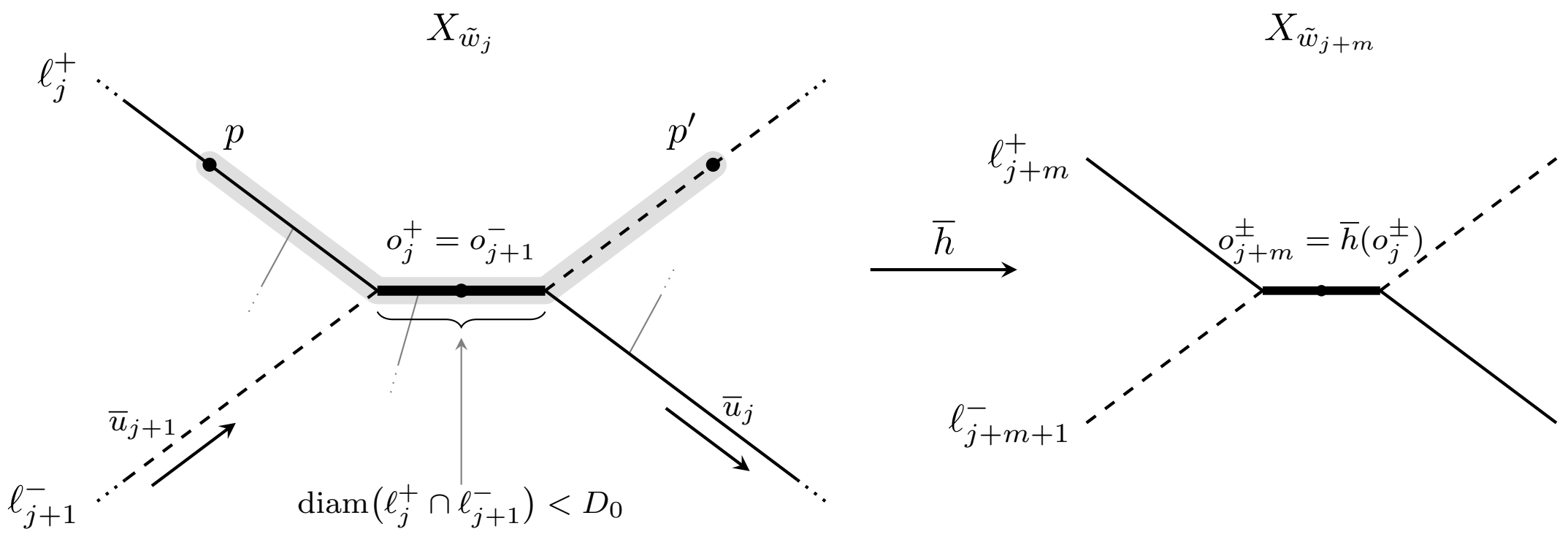

\begin{proof}
If $L\cap L'=\emptyset$, let $[o,o']$ be the bridge between $L$ and $L'$, with $o\in L$ and $o'\in L'$. Every geodesic from a point $p\in L$ to a point $p'\in L'$ passes through $o$ and $o'$, so $d(p,p')\ge d(p,o)+d(o',p')$. Otherwise $I:=L\cap L'$ is a point or a compact segment, and we let $o=o'$ be any point of $I$. For $p\in L$ and $p'\in L'$, let $\pi(p),\pi(p')\in I$ be the points of $I$ closest to $p$ and $p'$. The path $[p,\pi(p)]\cup[\pi(p),\pi(p')]\cup[\pi(p'),p']$ has no backtracking, because the two segments leaving $I$ along $L$ and $L'$ start in different directions (a common initial segment would lie in $L\cap L'=I$). So it is the geodesic from $p$ to $p'$, and
$$d(p,p')\ \ge\ d(p,\pi(p))+d(\pi(p'),p')\ \ge\ d(p,o)+d(o',p')-2D_0 ,$$
since $d(\pi(p),o),d(\pi(p'),o')\le\operatorname{diam}I\le D_0$ (see Figure~\ref{fig:lines-in-tree}).
\end{proof}

Recall from Section~\ref{subsec:axis_notation} that $u_j$ and $u_{j+1}$ both lie in $H_{\tilde w_j}$, so $\bar u_j=p_w(u_j)$ and $\bar u_{j+1}=p_w(u_{j+1})$ both act on the tree $X_{\tilde w_j}$. They are nontrivial, because $u_j,u_{j+1}\neq1$ and $p_w$ is injective on $H_{\tilde w_j}$ ($\tilde w_j$ is of fiber type). Since $F_w$ acts freely on $X_{\tilde w_j}$, they fix no point of $X_{\tilde w_j}$.

\begin{defn}\label{defn:axes_lj}
    For each $j\in\mathbb Z$, the element $u_j$ acts on two trees, $X_{\tilde w_{j-1}}$ and $X_{\tilde w_j}$. Define
$$\ell_j^-:=\mathrm{Ax}(\bar u_j)\subset X_{\tilde w_{j-1}},\quad \ell_j^+:=\mathrm{Ax}(\bar u_j)\subset X_{\tilde w_j},$$
with translation lengths
$$\lambda_j^-:=\|\bar u_j\|\ \text{on }X_{\tilde w_{j-1}},\quad \lambda_j^+:=\|\bar u_j\|\ \text{on }X_{\tilde w_j}.$$

So the tree $X_{\tilde w_j}$ contains the two lines $\ell_j^+$ (from the black vertex before $\tilde w_j$) and $\ell_{j+1}^-$ (from the black vertex after it).
\end{defn}

Recall that for $g\in G$, the map $\bar g$ denotes the isometry of trees induced by $g$, see (\ref{eq:bar_g}). In particular, $\bar h$ is the isometry induced by the element $h$ fixed in Convention~\ref{conv:fix_h}. We have the following basic properties of the lines $\ell_j^\pm$. We refer the reader to Figure~\ref{fig:lines-in-tree} for an illustration of the following lemma.

\begin{lem} \label{lem:bounded_axis}
    \begin{enumerate}
        \item \label{lem:bounded_axis_1} $\bar h(\ell_j^\pm)=\ell_{j+m}^\pm$ and $\lambda_{j+m}^\pm=\lambda_j^\pm$ for all $j$.

        \item  \label{lem:bounded_axis_2} Let $D_0:=\max_{1\le j\le m}\bigl(\lambda_j^++\lambda_{j+1}^-\bigr)$. Then for all $j\in\mathbb Z$,

$$\operatorname{diam}\bigl(\ell_j^+\cap\ell_{j+1}^-\bigr)<D_0 .$$

\item \label{lem:bounded_axis_3}  There are points $o_j^+\in\ell_j^+$ and $o_{j+1}^-\in\ell_{j+1}^-$, for all $j\in\mathbb Z$, with $\bar h(o_j^\pm)=o_{j+m}^\pm$, such that for all $p\in\ell_j^+$ and $p'\in\ell_{j+1}^-$,

$$d_X(p,p')\ \ge\ d_X(p,o_j^+)+d_X(o_{j+1}^-,p')-2D_0 .$$
    \end{enumerate}
\end{lem}

\begin{proof}
    (\ref{lem:bounded_axis_1}): By Section~\ref{subsec:axis_notation}, $u_{j+m}=hu_jh^{-1}$, hence $\bar u_{j+m}=\bar h\,\bar u_j\,\bar h^{-1}$ as isometries of $X_{\tilde w_{j+m}}$ and of $X_{\tilde w_{j+m-1}}$. Conjugation carries axes to axes and preserves translation length.

(\ref{lem:bounded_axis_2}): By Corollary~\ref{cor:trivialintersectionofpw}, $\langle\bar u_j\rangle\cap\langle\bar u_{j+1}\rangle=1$ in $F_w$. The group $F_w$ is free and acts freely on $X_{\tilde w_j}$. So Lemma~\ref{lem:axisintersection} gives $$\operatorname{diam}(\ell_j^+\cap\ell_{j+1}^-)<\lambda_j^++\lambda_{j+1}^-$$ By (\ref{lem:bounded_axis_1}) the right side is $m$-periodic in $j$, so it is at most $D_0$ for every $j$.

(\ref{lem:bounded_axis_3}): For $1\le j\le m$, apply Lemma~\ref{lem:boundedgates} to the lines $L=\ell_j^+$ and $L'=\ell_{j+1}^-$ in $X_{\tilde w_j}$, with the constant $D_0$ of (\ref{lem:bounded_axis_2}). This gives gates $o_j^+$ and $o_{j+1}^-$. For other $j$, define the gates by $o_{j+m}^+:=\bar h(o_j^+)$ and $o_{j+1+m}^-:=\bar h(o_{j+1}^-)$. Since $\bar h:X_{\tilde w_j}\to X_{\tilde w_{j+m}}$ is an isometry carrying $\ell_j^+,\ell_{j+1}^-$ to $\ell_{j+m}^+,\ell_{j+1+m}^-$ by (\ref{lem:bounded_axis_1}), the inequality is transported to all $j$.
\end{proof}

Each black vertex $\tilde b_j$ thus has two gates attached to it: $o_j^-\in\ell_j^-\subset X_{\tilde w_{j-1}}$, on the side of the previous white vertex, and $o_j^+\in\ell_j^+\subset X_{\tilde w_j}$, on the side of the next one.

\begin{cor} \label{cor:linesaxis}
    Let $1\le j\le N-1$. In the tree $X_{W_j}$, consider the lines $\bar q_j(\ell_j^+)$ and $\bar q_j(\ell_{j+1}^-)$. For all $p\in\bar q_j(\ell_j^+)$ and $p'\in\bar q_j(\ell_{j+1}^-)$,

$$d_X(p,p')\ \ge\ d_X\bigl(p,\ \bar q_j(o_j^+)\bigr)+d_X\bigl(\bar q_j(o_{j+1}^-),\ p'\bigr)-2D_0 .$$
\end{cor}

\begin{proof}
    Apply the isometry $\bar q_j^{-1}\colon X_{W_j}\to X_{\tilde w_j}$ and use Lemma~\ref{lem:bounded_axis}(\ref{lem:bounded_axis_3}).
\end{proof}

\subsection{Gates on the black lines}\label{subsec:gates_black_lines}

Let $\tilde e=(\tilde b,\tilde w)$ be an edge of $T$ over $e=(b,w)$, and write $\tilde b=g\bar b$ with $g\in G$. We now introduce coordinates on the plane $\widetilde M^{\tilde b}$ which are compatible with the coordinates on the group $G_{\tilde b}$ of Definition~\ref{defn:coordinates}. The basis we use is the basis $\{r_b,t_b\}$ of $G_b\cong\mathbb Z^2$ which is given by Proposition~\ref{prop:split-over-tori}.

First consider $\bar b$. The plane $\widetilde M^{\bar b}=\widetilde X_b\times\mathbb R$ is the universal cover of the torus $M^b=X_b\times S_b$, where both circles $X_b$ and $S_b$ have length $1$. The generators $r_b$ and $t_b$ of $G_b$ are represented by the loops $r_e$ and $t_e$, which go once around $X_b$ and once around $S_b$ (see Section~\ref{subsection:topological_representative}, item~(\ref{item:multitwist-decomp})). So $r_b$ and $t_b$ act on $\widetilde M^{\bar b}$ as translations by $1$ along the first and the second factor. We fix a base point in $\widetilde M^{\bar b}$ and identify $\widetilde M^{\bar b}$ with $\mathbb R^2$, so that the base point is the origin, the first factor is the $X$-axis and the second factor is the $Y$-axis. In these coordinates, $r_b$ acts by $(X,Y)\mapsto(X+1,Y)$ and $t_b$ acts by $(X,Y)\mapsto(X,Y+1)$. Hence the element $r_b^{X_0}t_b^{Y_0}$, that is, the element with coordinates $(X_0,Y_0)$, acts by $(X,Y)\mapsto(X+X_0,Y+Y_0)$.

We transport these coordinates to $\widetilde M^{\tilde b}=g\widetilde M^{\bar b}$ by $g$. Since $g$ is determined by $\tilde b$ up to right multiplication by $G_b$, which acts by translations, the coordinates on $\widetilde M^{\tilde b}$ are well defined up to a translation, which does not affect differences of points. Since $x\in G_{\tilde b}$ acts on $\widetilde M^{\tilde b}$ as $g^{-1}xg\in G_b$ acts on $\widetilde M^{\bar b}$, an element of $G_{\tilde b}$ with coordinates $(X_0,Y_0)$ acts on $\widetilde M^{\tilde b}$ by the translation $$(X,Y)\mapsto(X+X_0,Y+Y_0).$$ 

In these coordinates, $$d_{\tilde b}((X,Y),(X',Y'))=|X-X'|+|Y-Y'|.$$ For $z,z'\in\widetilde M^{\tilde b}$ we write $|z-z'|$ for the Euclidean distance in these coordinates, and we use it only to measure distances along lines. Every $g'\in G$ maps $\widetilde M^{\tilde b}$ to $\widetilde M^{g'\tilde b}$ by a translation in these coordinates, so $$|g'z-g'z'|=|z-z'|.$$

Recall that $p_e$ is an immersed circuit in $X_w$, and let $|p_e|\ge1$ be its number of edges. By Proposition~\ref{prop:split-over-tori}(\ref{prop:split-over-tori4}), the map $\iota_e^+$ sends the circle $X_b$ of length $1$, subdivided into $|p_e|$ edges, combinatorially onto $p_e$, hence with constant speed $|p_e|$. By item~(\ref{item:multitwist-decomp}) in Section~\ref{subsection:topological_representative}, the $X_w$-component of $I_e$ sends the point with coordinates $(X,Y)$ to the point of $p_e$ with parameter $X+k(e)Y$ modulo $1$. Lifting to universal covers, there are a line $\ell_{\tilde e}\subset X_{\tilde w}$ that is a lift of the immersed circuit $p_e$, an arc-length parametrization $\eta_{\tilde e}:\mathbb R\to\ell_{\tilde e}$ and a constant $c_{\tilde e}$ such that
\begin{equation}\label{eq:affine-attaching}
\pi_X\circ\iota_{\tilde e}(X,Y)=\eta_{\tilde e}\bigl(|p_e|\cdot(X+k(e)Y)+c_{\tilde e}\bigr).
\end{equation}
Similarly, the $S_w$-component of $I_e$ sends the point with coordinates $(X,Y)$ to the point of $S_w$ with parameter $Y$ modulo $1$. Lifting, there is a constant $c'_{\tilde e}$ such that
\begin{equation}\label{eq:vertical-attaching}
\text{the $\mathbb R$-component of $\iota_{\tilde e}(X,Y)$ is } Y+c'_{\tilde e}.
\end{equation}
This agrees with (\ref{D4}), because an element of $G_{\tilde b}$ with coordinates $(X_0,Y_0)$ is sent to an element $g\in G_{\tilde w}$ with $q_w(g)=Y_0$.

For each $j\in\mathbb Z$, let $L_j:=\mathcal H(\tilde b_j)\subset\widetilde M^{\tilde b_j}$. Let $(a_j,\gamma_j)$ be the coordinates of $u_j$ where $u_j$ is defined in Section~\ref{subsec:axis_notation} and $$\theta_j:=\sqrt{a_j^2+\gamma_j^2}>0.$$ Then $L_j$ is a line in direction $(a_j,\gamma_j)$, and $u_j$ translates it by Euclidean length $\theta_j$. For $z,z'\in L_j$ we write $|z-z'|$ for their Euclidean distance.

For any $g\in G$, the translate $gL_j$ is a line in $\widetilde M^{g\tilde b_j}$, and $|gz-gz'|=|z-z'|$.
By periodicity, $hL_j=L_{j+m}$ and $\theta_{j+m}=\theta_j$ (since $u_{j+m}=hu_jh^{-1}$ has the same coordinates as $u_j$, see the paragraph preceding (\ref{eq:chi-equivariant})). Let $$\theta:=\min_{1\le j\le m}\theta_j>0$$

We now study the images of the lines $L_j$ in the trees of the two neighbouring white vertices.

 For each $j$, let
$$\phi_j^+:=\pi_X\circ\iota_{\tilde e_j^+}\big|_{L_j}:L_j\to X_{\tilde w_j},\quad \phi_j^-:=\pi_X\circ\iota_{\tilde e_j^-}\big|_{L_j}:L_j\to X_{\tilde w_{j-1}} .$$

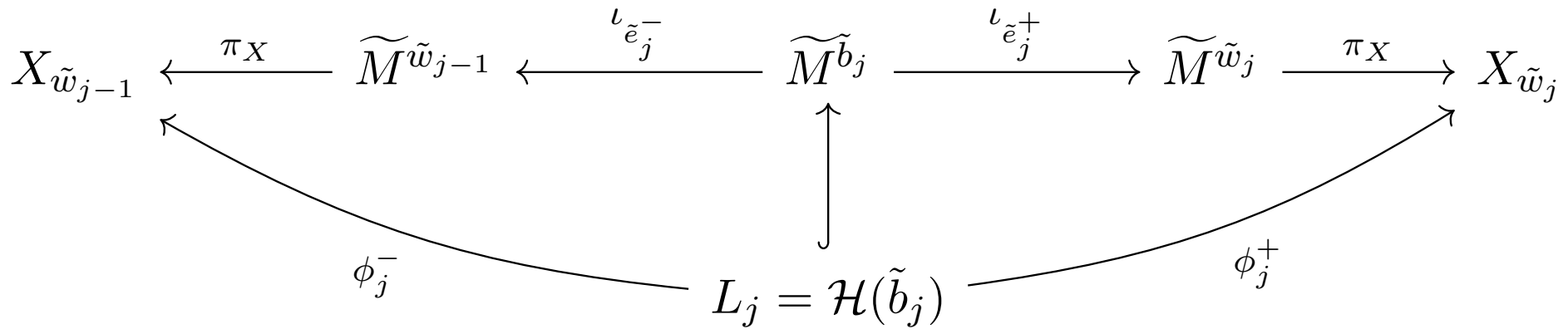
\[
\begin{tikzcd}[column sep=large,row sep=large]
X_{\tilde w_{j-1}} & \widetilde M^{\tilde w_{j-1}} \arrow[l,"\pi_X"'] & \widetilde M^{\tilde b_j} \arrow[l,"\iota_{\tilde e_j^-}"'] \arrow[r,"\iota_{\tilde e_j^+}"] & \widetilde M^{\tilde w_j} \arrow[r,"\pi_X"] & X_{\tilde w_j}\\
 & & L_j = \mathcal H(\tilde b_j) \arrow[u,hook] \arrow[ull,bend left=12,"\phi_j^-"] \arrow[urr,bend right=12,"\phi_j^+"'] & &
\end{tikzcd}
\]

We will see in Lemma~\ref{lem:comparephij} that $\phi_{j}^{\pm}$ are bijection. Thus, the following definition makes sense.
\begin{defn} \label{defn:gates}
    For $j\in\mathbb Z$, we pull the gates of Lemma~\ref{lem:bounded_axis}(\ref{lem:bounded_axis_3}) back to the black line: 
$$\hat o_j^-:=(\phi_j^-)^{-1}(o_j^-)\in L_j,\quad \hat o_j^+:=(\phi_j^+)^{-1}(o_j^+)\in L_j .$$

For $1\le j\le N$, define the \textit{gates on $\mathcal H(B_j)$} by

$$G_j^-:=q_{j-1}\,\hat o_j^-,\quad G_j^+:=q_j\,\hat o_j^+ .$$
\end{defn}

\begin{figure}[htb]
\centering
\begin{tikzpicture}[scale=0.85,>=stealth]
\begin{scope}[shift={(0,3.2)}]
\draw[thick] (-2.6,1.0)--(-0.9,0)--(0.3,0)--(2.0,1.0);
\draw[thick,gray!60!black,dash pattern=on 4pt off 2pt] (-2.6,-1.0)--(-0.9,0)--(0.3,0)--(2.0,-1.0);
\draw[line width=2.4pt] (-0.9,0)--(0.3,0);
\fill (-0.3,0) circle (2pt) node[above,font=\small] {$o_j^-$};
\node[font=\small,right] at (2.0,1.0) {$\ell_{j-1}^+$};
\node[font=\small,right] at (2.0,-1.0) {$\ell_j^-$};
\node[font=\small] at (-0.3,1.45) {$X_{\tilde w_{j-1}}$};
\end{scope}
\begin{scope}[shift={(8.4,3.2)}]
\draw[thick,gray!60!black,dash pattern=on 4pt off 2pt] (-2.6,-1.0)--(-0.9,0)--(0.3,0)--(2.0,-1.0);
\draw[thick] (-2.6,1.0)--(-0.9,0)--(0.3,0)--(2.0,1.0);
\draw[line width=2.4pt] (-0.9,0)--(0.3,0);
\fill (-0.3,0) circle (2pt) node[above,font=\small] {$o_j^+$};
\node[font=\small,left] at (-2.6,-1.0) {$\ell_j^+$};
\node[font=\small,right] at (2.0,1.0) {$\ell_{j+1}^-$};
\node[font=\small] at (-0.3,1.45) {$X_{\tilde w_j}$};
\end{scope}
\draw[thick] (0.5,0)--(7.9,0);
\node[font=\small,right] at (7.9,0) {$L_j=\mathcal H(\tilde b_j)$};
\fill (3.4,0) circle (2pt) node[below,font=\small] {$\hat o_j^-$};
\fill (5.0,0) circle (2pt) node[below,font=\small] {$\hat o_j^+$};
\draw[->,dashed] (3.3,0.15) .. controls (2.6,1.4) and (0.6,1.6) .. (-0.25,3.0);
\draw[->,dashed] (5.1,0.15) .. controls (5.8,1.4) and (7.4,1.6) .. (8.1,3.0);
\node[font=\small] at (2.75,1.55) {$\phi_j^-$};
\node[font=\small] at (5.65,1.55) {$\phi_j^+$};
\begin{scope}[shift={(0,-2.1)}]
\draw[thick] (-1.6,0)--(9.6,0);
\node[font=\small,right] at (9.6,0) {$\mathcal H(B_j)$};
\fill (-0.6,0) circle (2pt) node[above,font=\small] {$G_j^-=q_{j-1}\hat o_j^-$};
\fill (8.6,0) circle (2pt) node[above,font=\small] {$G_j^+=q_j\hat o_j^+$};
\draw[decorate,decoration={brace,mirror,amplitude=6pt}] (-0.6,-0.2)--(8.6,-0.2) node[midway,below=8pt,font=\small] {length $\ \ge\ \theta|\nu_j|-\Gamma_0$};
\end{scope}
\end{tikzpicture}
\caption{The gates $o_j^\pm$ in the trees $X_{\tilde w_{j-1}}$ and $X_{\tilde w_j}$, their preimages $\hat o_j^\pm$ on the line $L_j$, and the gates $G_j^\pm$ on $\mathcal H(B_j)$ (bottom), which are far apart (see Lemma~\ref{lem:lowerbound_Gj}).}
\label{fig:gates_black}
\end{figure}
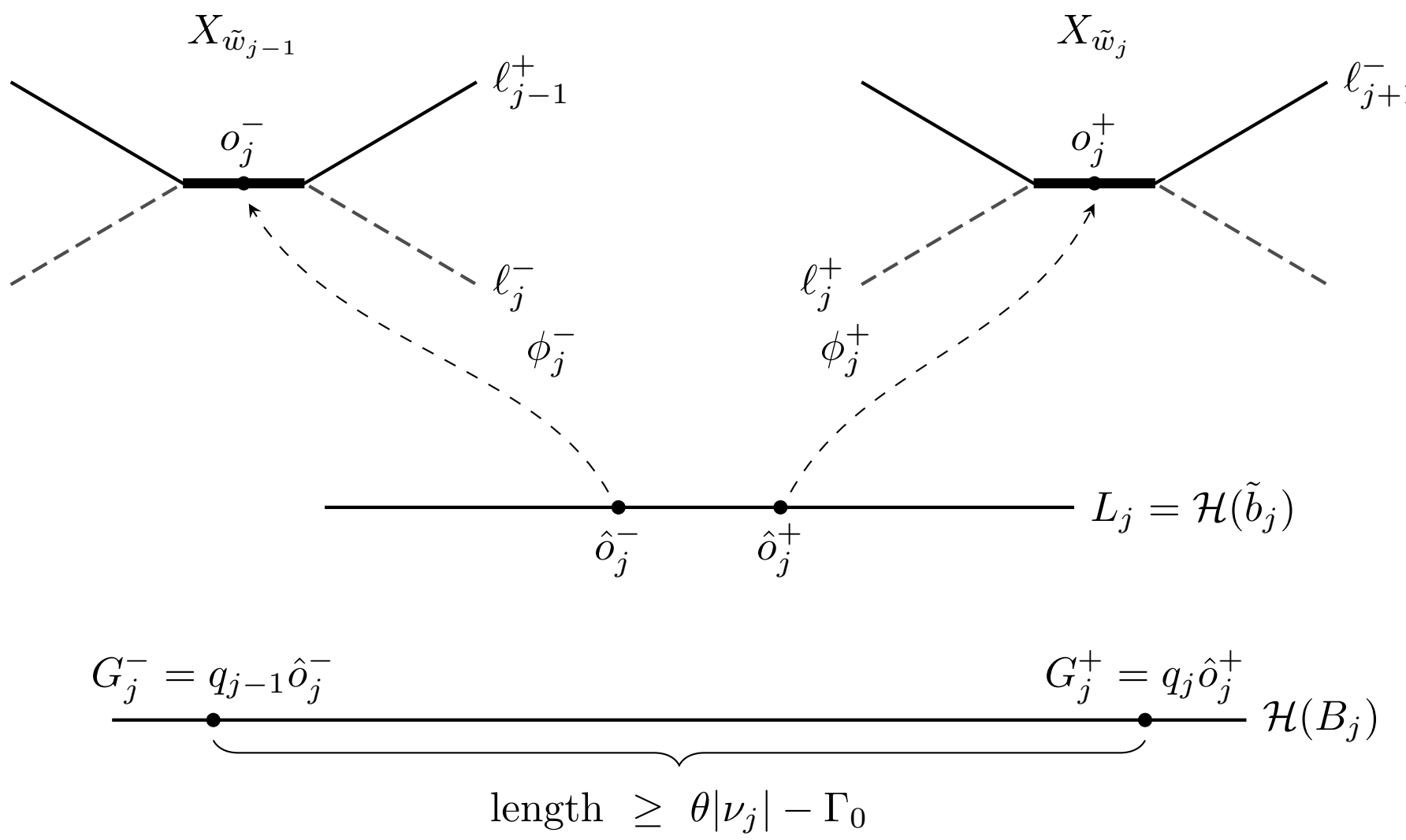

\begin{rem}\label{rem:gates_on_line}
    Both gates $G_{j}^{-}, G_{j}^{+}$ lie on $\mathcal H(B_j)=q_{j-1}L_j=q_jL_j$ (see Lemma~\ref{lem:adjacency}(\ref{lem:adjacency_1}) and Proposition~\ref{prop:horizontal pieces}(\ref{prop:horizontal pieces_1})). The gate $G_j^-$ is the reference gate as seen from $W_{j-1}$, and $G_j^+$ is the reference gate as seen from $W_j$ (see Figure~\ref{fig:gates_black}).
\end{rem}

Since $\phi_{j+m}^\pm\circ h=\bar h\circ\phi_j^\pm$ and $o_{j+m}^\pm=\bar h(o_j^\pm)$, we get $\hat o_{j+m}^\pm=h\,\hat o_j^\pm$, so $|\hat o_j^+-\hat o_j^-|$ is $m$-periodic. Let

$$\Gamma_0:=\max_{1\le j\le m}|\hat o_j^+-\hat o_j^-| .$$
\begin{lem} \label{lem:lowerbound_Gj}
    For $1\le j\le N$, we have
$$|G_j^+-G_j^-|\ \ge\ \theta\,|\nu_j|-\Gamma_0 .$$
\end{lem}

\begin{proof}
Since $q_j=q_{j-1}u_j^{\nu_j}$ and $q_{j-1}$ preserves distances along the line,

$$|G_j^+-G_j^-|=\bigl|q_{j-1}u_j^{\nu_j}\hat o_j^+-q_{j-1}\hat o_j^-\bigr|=\bigl|u_j^{\nu_j}\hat o_j^+-\hat o_j^-\bigr| .$$

The element $u_j^{\nu_j}$ translates $L_j$ by $|\nu_j|\theta_j$. By the triangle inequality on the line,

$$\bigl|u_j^{\nu_j}\hat o_j^+-\hat o_j^-\bigr|\ \ge\ \bigl|u_j^{\nu_j}\hat o_j^+-\hat o_j^+\bigr|-|\hat o_j^+-\hat o_j^-|\ \ge\ \theta_j|\nu_j|-\Gamma_0\ \ge\ \theta|\nu_j|-\Gamma_0 $$
\end{proof}

\begin{lem} \label{lem:comparephij}  $\phi_j^\pm$ is a bijection from $L_j$ onto the axis $\ell_j^\pm$ of Definition~\ref{defn:axes_lj}.  Moreover,
for all $z,z'\in L_j$, we have
\begin{equation} \label{eq:important}
d_X\bigl(\phi_j^\pm(z),\phi_j^\pm(z')\bigr)=\frac{\lambda_j^\pm}{\theta_j}\,|z-z'| .
\end{equation}
\end{lem}
\begin{proof}
    We treat $\phi_j^+$. The other case is identical. Write $\tilde e=\tilde e_j^+$, over $e$, and $\chi=\chi_{\tilde e}$.

 Let $z,z'\in L_j$. Since $L_j$ has direction $(a_j,\gamma_j)$, we have $z'-z=\tfrac{s}{\theta_j}(a_j,\gamma_j)$ with $|s|=|z-z'|$. By (\ref{eq:affine-attaching}), the parameters of $\phi_j^+(z)$ and $\phi_j^+(z')$ on $\ell_{\tilde e}$ differ by

$$|p_e|\cdot\frac{s}{\theta_j}\bigl(a_j+k(e)\gamma_j\bigr)=|p_e|\,\chi(u_j)\,\frac{s}{\theta_j}.$$

As $\eta_{\tilde e}$ is an arc-length parametrization,

\begin{equation}\label{eq:phi-param}
d_X\bigl(\phi_j^+(z),\phi_j^+(z')\bigr)=\frac{|p_e|\,|\chi(u_j)|}{\theta_j}\,|z-z'| .
\end{equation}

 Since $\chi(u_j)\neq0$ (see Lemma~\ref{lem:chi_nonzero}), the map $\phi_j^+$ is affine with nonzero rate in the parameter of $\ell_{\tilde e}$, so it is a bijection of the line $L_j$ onto the line $\ell_{\tilde e}$. By equivariance of $\iota_{\tilde e}$ under $u_j$ (\ref{D5}) and $\pi_X\circ u_j=\bar u_j\circ\pi_X$, we have $\phi_j^+(u_jz)=\bar u_j\,\phi_j^+(z)$. So the line $\ell_{\tilde e}$ is invariant under the fixed-point-free isometry $\bar u_j$, and an invariant line of such an isometry is its axis. Hence $\ell_{\tilde e}=\ell_j^+$.

 Take $z'=u_jz$, so $|z-z'|=\theta_j$. By (\ref{eq:phi-param}), $\bar u_j$ moves the point $\phi_j^+(z)$ of its axis by $|p_e|\,|\chi(u_j)|$. The displacement of a point on the axis is the translation length, so $\lambda_j^+=|p_e|\,|\chi(u_j)|$. Substituting into (\ref{eq:phi-param}) gives (\ref{eq:important}).
\end{proof}

\subsection{Proof of Theorem~\ref{thm:lower_bound}}

Let
$$\varsigma:=\min_{1\le j\le m}\min\Bigl\{\frac{\lambda_j^+}{\theta_j},\ \frac{\lambda_j^-}{\theta_j}\Bigr\}>0 .$$
Since $\lambda_j^\pm$ (see Lemma~\ref{lem:bounded_axis}(\ref{lem:bounded_axis_1})) and $\theta_j$ are $m$-periodic in $j$, we have $\lambda_j^\pm/\theta_j\ge\varsigma$ for all $j$. Let $K$ be the constant of (\ref{D3}), let $D_0$ be the constant of Lemma~\ref{lem:bounded_axis}(\ref{lem:bounded_axis_2}), and set
$$c_\ast:=\varsigma/K,\qquad C_\ast:=2D_0/K+4K,$$
and
\begin{equation}\label{eq:choice-mu}
\mu:=\max\Bigl\{1,\ \frac{4(c_\ast\Gamma_0+C_\ast)}{c_\ast\theta}\Bigr\}.
\end{equation}
The constants $\varsigma$, $c_\ast$, $C_\ast$, $\theta$ and $\Gamma_0$ depend only on $K$, $D_0$, $h$, the generators $u_j$ and the gates. They do not depend on $\mu$, on the sequence $(s_j)$, or on $n$. From now on, the sequence $(s_j)$ of Definition~\ref{defn:q_j} is chosen with this $\mu$.

\begin{prop} \label{prop:lowerboundofrhoninHi} Let $\Sigma(n):=\sum_{j=1}^{nm}|s_j|$. There are constants $c_3>0$ and $C_3\ge0$ such that for all $n\ge4$,
$$|\rho_n|_{H_i}\ \ge\ c_3\,\Sigma(n-1)-C_3 .$$
\end{prop}

\begin{proof}
Let $n\ge4$, let $N:=nm$, and let $\gamma$ be a path in $\Omega_{T_i}$ from $x_0$ to $\rho_nx_0$. Let $y_1,\dots,y_N$ be the points of Lemma~\ref{lem:supremum}. Thus $y_j\in\mathcal H(B_j)$, and
\begin{equation}\label{eq:path-yj}
\ell_d(\gamma)\ \ge\ \sum_{j=1}^{N-1}d(y_j,y_{j+1})
\end{equation}
by Lemma~\ref{lem:supremum}(\ref{lem:supremum_1}),(\ref{lem:supremum_3}). We now bound each distance $d(y_j,y_{j+1})$ from below in terms of the distances from $y_j$ and $y_{j+1}$ to the gates $G_j^+$ and $G_{j+1}^-$ (see Definition~\ref{defn:gates}).

\smallskip
\noindent\textit{Claim 1.} For $1\le j\le N-1$ and all points $y\in\mathcal H(B_j)$ and $y'\in\mathcal H(B_{j+1})$,
$$d(y,y')\ \ge\ c_\ast\bigl(|y-G_j^+|+|y'-G_{j+1}^-|\bigr)-C_\ast .$$

\smallskip
\noindent\textit{Proof of Claim 1.} Let $\tilde e:=q_j\tilde e_j^+$ and $\tilde e':=q_j\tilde e_{j+1}^-$ be the edges joining $B_j$ and $B_{j+1}$ to $W_j$ (see Lemma~\ref{lem:WjBj}(\ref{lem:WjBj_2})), and let
$$\hat y:=\iota_{\tilde e}(y),\quad \hat y':=\iota_{\tilde e'}(y')\quad\text{in }\widetilde M^{W_j}.$$
The vertical segment $\{y\}\times[-1,1]$ of the cylinder of $\tilde e$ joins $y$ to $\hat y$, so $d(y,\hat y)\le2\le 2K$ by (\ref{eq:vertical-segment}). Likewise $d(y',\hat y')\le2K$. Hence
$$d(y,y')\ \ge\ d(\hat y,\hat y')-4K .$$
Both $\hat y,\hat y'$ lie in $\widetilde M^{W_j}$. By (\ref{D3}), and since $d_{W_j}$ is the $\ell^1$-product metric on $X_{W_j}\times\mathbb R$ (see Section~\ref{subsec:metrics}),
$$d(\hat y,\hat y')\ \ge\ \tfrac1K\,d_{W_j}(\hat y,\hat y')\ \ge\ \tfrac1K\,d_X\bigl(\pi_X\hat y,\ \pi_X\hat y'\bigr).$$
Write $y=q_jz$ with $z\in L_j$, and $y'=q_jz'$ with $z'\in L_{j+1}$, which is possible by Lemma~\ref{lem:WjBj}(\ref{lem:WjBj_1}). By Lemma~\ref{lem:WjBj}(\ref{lem:WjBj_3}) and $\pi_X\circ q_j=\bar q_j\circ\pi_X$, we have
$$\pi_X\hat y=\pi_X\bigl(q_j\,\iota_{\tilde e_j^+}(z)\bigr)=\bar q_j\,\phi_j^+(z)\in\bar q_j(\ell_j^+),\quad \pi_X\hat y'=\bar q_j\,\phi_{j+1}^-(z')\in\bar q_j(\ell_{j+1}^-).$$
By Corollary~\ref{cor:linesaxis}, and since $\bar q_j$ is an isometry,
$$d_X(\pi_X\hat y,\pi_X\hat y')\ \ge\ d_X\bigl(\phi_j^+(z),\,o_j^+\bigr)+d_X\bigl(o_{j+1}^-,\,\phi_{j+1}^-(z')\bigr)-2D_0 .$$
By Definition~\ref{defn:gates}, $o_j^+=\phi_j^+(\hat o_j^+)$ and $o_{j+1}^-=\phi_{j+1}^-(\hat o_{j+1}^-)$. Hence Lemma~\ref{lem:comparephij} gives
$$d_X\bigl(\phi_j^+(z),o_j^+\bigr)=\frac{\lambda_j^+}{\theta_j}|z-\hat o_j^+|\ \ge\ \varsigma\,|z-\hat o_j^+|\quad\text{and}\quad d_X\bigl(o_{j+1}^-,\phi_{j+1}^-(z')\bigr)\ \ge\ \varsigma\,|z'-\hat o_{j+1}^-| .$$
Since $q_j$ preserves distances along the lines, $|z-\hat o_j^+|=|y-G_j^+|$ and $|z'-\hat o_{j+1}^-|=|y'-G_{j+1}^-|$. Combining these inequalities, we obtain
$$d(y,y')\ \ge\ \frac{\varsigma}{K}\bigl(|y-G_j^+|+|y'-G_{j+1}^-|\bigr)-\frac{2D_0}{K}-4K\ =\ c_\ast\bigl(|y-G_j^+|+|y'-G_{j+1}^-|\bigr)-C_\ast ,$$
which proves Claim~1.

\smallskip
\noindent\textit{Claim 2.} We have
$$\ell_d(\gamma)\ \ge\ c_\ast\theta\sum_{j=2}^{N-1}|\nu_j|-(c_\ast\Gamma_0+C_\ast)\,N .$$

\smallskip
\noindent\textit{Proof of Claim 2.} By (\ref{eq:path-yj}) and Claim~1 applied to $y=y_j$ and $y'=y_{j+1}$,
$$\ell_d(\gamma)\ \ge\ c_\ast\sum_{j=1}^{N-1}\bigl(|y_j-G_j^+|+|y_{j+1}-G_{j+1}^-|\bigr)-C_\ast(N-1).$$
Regrouping the sum and using the triangle inequality on the line $\mathcal H(B_j)$, we get
$$\sum_{j=1}^{N-1}\bigl(|y_j-G_j^+|+|y_{j+1}-G_{j+1}^-|\bigr)\ \ge\ \sum_{j=2}^{N-1}\bigl(|y_j-G_j^-|+|y_j-G_j^+|\bigr)\ \ge\ \sum_{j=2}^{N-1}|G_j^+-G_j^-| .$$
By Lemma~\ref{lem:lowerbound_Gj},
$$\sum_{j=2}^{N-1}|G_j^+-G_j^-|\ \ge\ \theta\sum_{j=2}^{N-1}|\nu_j|-\Gamma_0 N .$$
Putting these inequalities together proves Claim~2.

\smallskip
\noindent\textit{Claim 3.} We have
$$\ell_d(\gamma)\ \ge\ \frac{c_\ast\theta}{2}\sum_{j=2}^{N-1}|\nu_j| .$$

\smallskip
\noindent\textit{Proof of Claim 3.} Since $|\nu_j|=|\kappa_j||s_j|\ge\mu$ (see Lemma~\ref{lem:HN19}(\ref{lem:HN19_1})), we have $\sum_{j=2}^{N-1}|\nu_j|\ge(N-2)\mu$. Since $n\ge4$, we have $N\ge4$, hence $N\le2(N-2)$. By the choice (\ref{eq:choice-mu}) of $\mu$,
$$(c_\ast\Gamma_0+C_\ast)N\ \le\ 2(c_\ast\Gamma_0+C_\ast)(N-2)\ \le\ \frac{2(c_\ast\Gamma_0+C_\ast)}{\mu}\sum_{j=2}^{N-1}|\nu_j|\ \le\ \frac{c_\ast\theta}{2}\sum_{j=2}^{N-1}|\nu_j| .$$
Subtracting this from the bound of Claim~2 proves Claim~3.

\smallskip
We now conclude the proof. Since $m\ge1$, we have $(n-1)m\le nm-1=N-1$. Using $|\nu_j|\ge|s_j|$, we get
$$\sum_{j=2}^{N-1}|\nu_j|\ \ge\ \sum_{j=2}^{(n-1)m}|s_j|=\Sigma(n-1)-|s_1| .$$
Since $\gamma$ was an arbitrary path in $\Omega_{T_i}$ from $x_0$ to $\rho_nx_0$, Claim~3 gives
$$\widehat d_{\Omega_{T_i}}(x_0,\rho_nx_0)\ \ge\ \frac{c_\ast\theta}{2}\bigl(\Sigma(n-1)-|s_1|\bigr).$$
By Corollary~\ref{cor:Hiundistorted}, $H_i$ has a finite generating set $S_i$, with respect to which $|\cdot|_{H_i}$ is computed. Let $\lambda_i:=\max\bigl\{1,\ \max_{s\in S_i}\widehat d_{\Omega_{T_i}}(x_0,sx_0)\bigr\}$, which is finite by Remark~\ref{rem:remarkonHi}. Since $H_i$ acts on $(\Omega_{T_i},\widehat d_{\Omega_{T_i}})$ by isometries, the triangle inequality gives $\widehat d_{\Omega_{T_i}}(x_0,gx_0)\le\lambda_i|g|_{H_i}$ for all $g\in H_i$. Hence
$$|\rho_n|_{H_i}\ \ge\ \frac{c_\ast\theta}{2\lambda_i}\,\Sigma(n-1)-\frac{c_\ast\theta}{2\lambda_i}\,|s_1| .$$
Take $c_3:=c_\ast\theta/(2\lambda_i)$ and $C_3:=c_3|s_1|$.
\end{proof}

We are now ready for the proof of Theorem~\ref{thm:lower_bound}.

\begin{proof}[ Proof of Theorem~\ref{thm:lower_bound}]
    Since $h^{-1}$ is also loxodromic and $\omega_i(h^{-1})=\omega_i(h)^{-1}$, we may replace $h$ by $h^{-1}$ and assume $\omega_i(h)\le1$, as in Convention~\ref{conv:fix_h}. Then all constructions of Sections~\ref{sec:labelling}--\ref{sec:overlap} apply to $T_i$ and $h$, and $\Omega=\omega_i(h)^{-1}\ge1$ is given by (\ref{eq:Omega}). Moreover, $\Omega>1$ if and only if $\omega_i(h)\neq1$.
Let $\rho_n\in H_i$ be the elements of Definition~\ref{defn:q_j}, with $\mu$ as in (\ref{eq:choice-mu}). By Proposition~\ref{prop:Ambient_upper_bound} and Proposition~\ref{prop:lowerboundofrhoninHi}, for all $n\ge4$,
$$|\rho_n|_G\le C_1n+C_2\quad\text{and}\quad|\rho_n|_{H_i}\ge c_3\Sigma(n-1)-C_3 .$$
By the definition of distortion,
$\Delta_{H_i}^G(C_1n+C_2)\ \ge\ c_3\,\Sigma(n-1)-C_3\quad(n\ge4).$
Since $\Sigma(n-1)\succeq(n-1)^2+\Omega^{\,n-1}$, and $(n-1)^2+\Omega^{n-1}\simeq n^2+\Omega^n$, and $\Delta_{H_i}^G$ is non-decreasing, a linear change of variable gives
$\Delta_{H_i}^G(x)\ \succeq\ x^2+\Omega^{\,x}.$
This gives (\ref{thm:lower_bound_1}) in all cases. If $\omega_i(h)\neq1$, then $\Omega>1$ and $x^2+\Omega^x\simeq2^x$ (see Section~\ref{subsec:distortion}), which gives (\ref{thm:lower_bound_2}).
\end{proof}

\subsection{Characterization for horizontal components}

\begin{thm}\label{thm:dichotomy}
Let $T_i$ be a horizontal component of $H$ and $H_i=\mathrm{Stab}_H(T_i)$. The following are equivalent:
\begin{enumerate}
\item\label{dich1} no element of $H_i$ acts loxodromically on $T_i$,
\item\label{dich2} $H_i$ fixes a vertex of $T_i$,
\item\label{dich3} $T_i$ is bounded,
\item\label{dich4} there is a white vertex $\tilde w\in T_i$ such that $H_i=H_{\tilde w}$ and $T_i$ is the combinatorial ball of radius $2$ about $\tilde w$ in $T_i$, and moreover, for every black neighbour $\tilde b$ of $\tilde w$ in $T_i$ and every neighbour $\tilde w'\neq\tilde w$ of $\tilde b$ in $T_i$, we have $H_{\tilde w'}=H_{\tilde b}\cong\mathbb Z$,
\item\label{dich5} $H_i$ is undistorted in $G$.
\end{enumerate}
If these conditions fail, then $\Delta^G_{H_i}(x)\succeq x^2$, and $\Delta^G_{H_i}(x)\succeq 2^x$ whenever the dilation homomorphism $\omega_i$ is nontrivial.
\end{thm}

\begin{proof}
(\ref{dich1})$\Rightarrow$(\ref{dich2}): $H_i$ is finitely generated (see Corollary~\ref{cor:Hiundistorted}) and acts on $T_i$ without inversions and without loxodromic elements, so it fixes a vertex \cite[Corollary~7.3]{Bas93}.

(\ref{dich2})$\Rightarrow$(\ref{dich4}): Let $H_i$ fix $\tilde v\in T_i$. Since $H_{\tilde v}\le H_i$ (see Lemma~\ref{lem:basifactsofTi}(\ref{lem:basifactsofTi_1})), we get $H_i=H_{\tilde v}$. If $\tilde v$ is black, then for a neighbour $\tilde w$ of $\tilde v$ in $T_i$ we have $H_{\tilde v}\le H_{\tilde w}\le H_i=H_{\tilde v}$ by (\ref{D5}), so we may replace $\tilde v$ by the white vertex $\tilde w$. So let $\tilde w:=\tilde v$ be white. Let $\tilde w-\tilde b-\tilde w'$ be a path in $T_i$. Then $H_{\tilde w'}\le H_i=H_{\tilde w}$ fixes $\tilde w$ and $\tilde w'$, hence fixes $\tilde b$. Thus $H_{\tilde w'}=H_{\tilde b}$, which is infinite cyclic by Lemma~\ref{lem:vertices_of_Ti}(\ref{lem:vertices_of_Ti_1}). If $\tilde b''\neq\tilde b$ were another black neighbour of $\tilde w'$ in $T_i$, then $H_{\tilde b''}\le H_{\tilde w'}=H_{\tilde b}$, so $H_{\tilde b''}=H_{\tilde b}\cap H_{\tilde b''}=1$ by Lemma~\ref{lem:consectutive_blacks}, contradicting Lemma~\ref{lem:vertices_of_Ti}(\ref{lem:vertices_of_Ti_1}). Hence $T_i$ is the combinatorial ball of radius $2$ about $\tilde w$ in $T_i$.

(\ref{dich4})$\Rightarrow$(\ref{dich3}) is clear, and (\ref{dich3})$\Rightarrow$(\ref{dich1}) holds because a loxodromic element of $H_i$ has its axis in $T_i$, and an axis is unbounded.

(\ref{dich4})$\Rightarrow$(\ref{dich5}): $H_i=H_{\tilde w}$ is a finitely generated subgroup of $G_{\tilde w}\cong F_w\times\mathbb Z$, hence undistorted in $G_{\tilde w}$ \cite[Lemma~8.6]{Lin25}, and $G_{\tilde w}$ is undistorted in $G$ \cite[Lemma~4.1]{Mut24}.

(\ref{dich5})$\Rightarrow$(\ref{dich1}) and the last statement: if some $h\in H_i$ is loxodromic on $T_i$, then $\Delta^G_{H_i}(x)\succeq x^2$ by Theorem~\ref{thm:lower_bound}, so $H_i$ is distorted. If $\omega_i\neq1$, pick $h\in H_i$ with $\omega_i(h)\neq1$. Then $h$ is not elliptic, because $\omega_i$ is trivial on elliptic elements (see Proposition~\ref{prop:dilation_basis}(\ref{prop:dilation_basis_3})), so $h$ is loxodromic and Theorem~\ref{thm:lower_bound}(\ref{thm:lower_bound_2}) gives $\Delta^G_{H_i}(x)\succeq 2^x$.
\end{proof}

\begin{lem}\label{lem:one_piece}
For a horizontal component $T_i$ the following are equivalent:
\begin{enumerate}
\item\label{one_piece_1} $T_i$ consists of one piece,
\item\label{one_piece_2} for every black vertex $\tilde b$ of $T_i$, all but at most one of the neighbours $\tilde w$ of $\tilde b$ in $T_i$ satisfy $H_{\tilde w}=H_{\tilde b}$,
\item\label{one_piece_3} $T_i$ is bounded.
\end{enumerate}
\end{lem}

\begin{proof}
We first note: if $\tilde w\in T_i$ is white and $H_{\tilde w}=H_{\tilde b}$ for some neighbour $\tilde b$ of $\tilde w$ in $T_i$, then $\tilde w$ is a leaf of $T_i$. Indeed, if $\tilde b''\neq\tilde b$ were another black neighbour of $\tilde w$ in $T_i$, then $H_{\tilde b''}\le H_{\tilde w}=H_{\tilde b}$ and $H_{\tilde b}\cap H_{\tilde b''}=1$ (see Lemma~\ref{lem:consectutive_blacks}), so $H_{\tilde b''}=1$, contradicting Lemma~\ref{lem:vertices_of_Ti}(\ref{lem:vertices_of_Ti_1}).

(\ref{one_piece_3})$\Rightarrow$(\ref{one_piece_1}) and (\ref{one_piece_3})$\Rightarrow$(\ref{one_piece_2}): By Theorem~\ref{thm:dichotomy}(\ref{dich4}), $T_i$ is the ball of radius $2$ about a white vertex $\tilde w$, every black vertex of $T_i$ is adjacent to $\tilde w$, and every white vertex $\tilde w'\neq\tilde w$ of $T_i$ is a leaf with $H_{\tilde w'}=H_{\tilde b}$ for its unique neighbour $\tilde b$.

(\ref{one_piece_1})$\Rightarrow$(\ref{one_piece_3}): If $T_i$ is unbounded, it contains a geodesic $\tilde w_0-\tilde b_1-\tilde w_1-\tilde b_2-\tilde w_2-\tilde b_3-\tilde w_3$, and $\tilde w_1,\tilde w_2$ are two distinct essential pieces.

(\ref{one_piece_2})$\Rightarrow$(\ref{one_piece_1}): Suppose that $\tilde w_1\neq\tilde w_2$ are essential pieces, and let $\tilde b$ be the black vertex following $\tilde w_1$ on the geodesic from $\tilde w_1$ to $\tilde w_2$, and let $\tilde u$ be the white vertex following $\tilde b$ on this geodesic. By (\ref{one_piece_2}) and the first paragraph, at most one neighbour of $\tilde b$ in $T_i$ is not a leaf of $T_i$, so $\tilde w_1$ or $\tilde u$ is a leaf of $T_i$. Since $\tilde w_1$ is essential, $\tilde u$ is a leaf. A leaf of a tree is the endpoint of every geodesic through it, so $\tilde u=\tilde w_2$, which is essential, a contradiction.
\end{proof}

\section{Reduction to the horizontal components}\label{sec:undistortion}

In this section we show that the distortion of $H$ is controlled by the distortions of the stabilizers of its horizontal components (see Theorem~\ref{thm:reduction}). Together with the lower bounds of Section~\ref{sec:overlap} and the upper bounds of Section~\ref{sec:upper}, this gives the main theorem (see Theorem~\ref{thm:main}). 

 Throughout this section, $M_\varphi$ is the topological representative of Section~\ref{subsection:lineartopological}.

Let $\tilde w$ be a white vertex of $T_H$. We use the tree $X_{\tilde w}$ and the projection $\pi_X\colon\widetilde M^{\tilde w}\to X_{\tilde w}$ of Section~\ref{sec:overlap}, and we define
\[
\widehat C_{\tilde w}:=\pi_X\bigl(\mathcal H(\tilde w)\bigr)\subset X_{\tilde w}.
\]
For the representatives $\tilde w_i$ of Section~\ref{subsec:horizontalpiece_white} this is the tree $\widehat C_{\tilde w_i}$, and $\widehat C_{h\tilde w}=\bar h(\widehat C_{\tilde w})$ for $h\in H$. Hence $\widehat C_{\tilde w}$ is a subtree and a subgraph of $X_{\tilde w}$, and $H_{\tilde w}$ acts on it, through $g\mapsto\bar g$, with finitely many orbits of edges. We write $P_{\tilde w}$ for the image of $H_{\tilde w}$ in $\mathrm{Isom}(X_{\tilde w})$.

For an edge $\tilde e=(\tilde b,\tilde w)$ of $T$, let $$\ell_{\tilde e}=\pi_X\iota_{\tilde e}(\widetilde M^{\tilde b})\subset X_{\tilde w}$$ be the line of (\ref{eq:affine-attaching}). For $z\in\widetilde M^{\tilde b}$ with coordinates $(X,Y)$ we write $\chi_{\tilde e}(z):=X+k(e)Y$, so that by (\ref{eq:affine-attaching})
\begin{equation}\label{eq:chi-distance}
d_X\bigl(\pi_X\iota_{\tilde e}(z),\pi_X\iota_{\tilde e}(z')\bigr)=|p_e|\,\bigl|\chi_{\tilde e}(z)-\chi_{\tilde e}(z')\bigr|\qquad(z,z'\in\widetilde M^{\tilde b}).
\end{equation}

In particular an element $x\in G_{\tilde b}$ acts on $\ell_{\tilde e}$ as the translation by $|p_e|\chi_{\tilde e}(x)$. 

Since $G_{\tilde e}=G_{\tilde b}$ contains $\langle t_{\tilde w}\rangle=\ker p_w$, the group $p_w(G_{\tilde e})$ is the infinite cyclic group $\langle\bar p_e\rangle$ stabilizing $\ell_{\tilde e}$, it is a maximal cyclic subgroup of the free group $p_w(G_{\tilde w})$ (see Proposition~\ref{prop:split-over-tori}(\ref{prop:split-over-tori4})), and $G_{\tilde e}=p_w^{-1}(\langle\bar p_e\rangle)$. 

Finally, for a subtree $Z$ of a tree we write $\mathrm{proj}_Z$ for the closest point projection, which is $1$-Lipschitz, and we define
\[
Q_{\tilde e}:=\mathrm{proj}_{\ell_{\tilde e}}\bigl(\widehat C_{\tilde w}\bigr)\subset\ell_{\tilde e}.
\]

\begin{lem}\label{lem:tree_gates}
There are constants $D_0,D_0'\ge0$, depending only on $M_\varphi$ and $H$, such that the following hold for every white vertex $\tilde w$ of $T_H$ and all edges $\tilde e=(\tilde b,\tilde w)$, $\tilde e'=(\tilde b',\tilde w)$ of $T_H$.
\begin{enumerate}
\item\label{tree_gates_1} $\pi_X\iota_{\tilde e}(\mathcal H(\tilde b))\subset\ell_{\tilde e}\cap\widehat C_{\tilde w}$. In particular, $\ell_{\tilde e}\cap\widehat C_{\tilde w}\neq\emptyset$ and $Q_{\tilde e}=\ell_{\tilde e}\cap\widehat C_{\tilde w}$.
\item\label{tree_gates_2} If $\chi_{\tilde e}(H_{\tilde b})\neq0$, then $\ell_{\tilde e}\subset\widehat C_{\tilde w}$, so $Q_{\tilde e}=\ell_{\tilde e}$. This holds if $\tilde b$ is of type $(H)$, of type (2), or of type $(V)$ and $\tilde w$ is of fiber type.
\item\label{tree_gates_3} If $\chi_{\tilde e}(H_{\tilde b})=0$, then $\operatorname{diam}Q_{\tilde e}\le D_0$. This holds if $\tilde b$ is of type (0), or of type $(V)$ and $\tilde w$ is of compact-fiber type.
\item\label{tree_gates_4} If $\tilde e\neq\tilde e'$, then $$\operatorname{diam}\mathrm{proj}_{\ell_{\tilde e}}(\ell_{\tilde e'})\le D_0' \quad \text{and} \quad \mathrm{proj}_{\ell_{\tilde e}}(\ell_{\tilde e'})\subset\mathcal N_{D_0'}(Q_{\tilde e})$$
\end{enumerate}
\end{lem} 

\begin{proof}
(\ref{tree_gates_1}): By Proposition~\ref{prop:horizontal pieces}(\ref{prop:horizontal pieces_3}), we have $\iota_{\tilde e}(\mathcal H(\tilde b))\subset\mathcal H(\tilde w)$, and $\mathcal H(\tilde b)\neq\emptyset$. The claim follows since the projection of a subtree onto a line which it meets is their intersection.

(\ref{tree_gates_2}): Let $u\in H_{\tilde b}$ with $\chi_{\tilde e}(u)\neq0$. The set $\pi_X\iota_{\tilde e}(\mathcal H(\tilde b))$ is nonempty and invariant under $\bar u$, which translates $\ell_{\tilde e}$ nontrivially, so the convex hull of this set is $\ell_{\tilde e}$.

Since $\widehat C_{\tilde w}$ is convex, (\ref{tree_gates_1}) gives $\ell_{\tilde e}\subset\widehat C_{\tilde w}$. For type $(H)$, and for type $(V)$ with $\tilde w$ of fiber type, we have $\chi_{\tilde e}(H_{\tilde b})\neq0$ by Lemma~\ref{lem:chi_nonzero}. For type (2), $H_{\tilde b}$ has finite index in $G_{\tilde b}$ and $\chi_{\tilde e}\neq0$ on $G_{\tilde b}$. (Alternatively, for type $(V)$ with $\tilde w$ of fiber type, $H_{\tilde b}$ contains a power $t_{\tilde w''}^a$, $a\neq0$, where $\tilde w''$ is the compact-fiber neighbour of $\tilde b$, and $\chi_{\tilde e}(t_{\tilde w''})=k(e)-k(e'')\neq0$ by Lemma~\ref{lem:slope_vertical_extended}(\ref{lem:slope_vertical_3}).)

(\ref{tree_gates_3}): Since $H_{\tilde w}\cap G_{\tilde e}=H_{\tilde b}$ and $G_{\tilde e}=p_w^{-1}(\langle\bar p_e\rangle)$, we have $P_{\tilde w}\cap\langle\bar p_e\rangle=p_w(H_{\tilde b})$. Each $x\in H_{\tilde b}$ acts on $\ell_{\tilde e}$ as the translation by $|p_e|\chi_{\tilde e}(x)=0$, so $p_w(x)$ fixes a point of $X_{\tilde w}$. As $p_w(G_{\tilde w})$ acts freely, $p_w(x)=1$. Hence $P_{\tilde w}\cap\langle\bar p_e\rangle=1$.

By (\ref{tree_gates_1}), $Q_{\tilde e}=\ell_{\tilde e}\cap\widehat C_{\tilde w}$ is a segment, consisting of edges of $X_{\tilde w}$. Let $n_E$ be the number of $P_{\tilde w}$-orbits of edges of $\widehat C_{\tilde w}$, and suppose $Q_{\tilde e}$ contains consecutive edges $a_1,\dots,a_N$, oriented in the direction of translation of $\bar p_e$, with $N>2n_E|p_e|$.  There are $i<j$ such that $a_j=ga_i$ as oriented edges for some $g\in P_{\tilde w}$, and $j-i=m|p_e|$ with $m\ge1$. Then also $\bar p_e^{\,m}a_i=a_j$, so $g^{-1}\bar p_e^{\,m}$ fixes $a_i$, and freeness gives $g=\bar p_e^{\,m}\in P_{\tilde w}\cap\langle\bar p_e\rangle\setminus\{1\}$, a contradiction. Hence $$\operatorname{diam}Q_{\tilde e}\le2n_E|p_e|+1$$ Since $H$ has finitely many orbits of vertices of $T_H$ and $\Delta$ has finitely many edges, this bound is uniform. Let $D_0$ be the maximum.

For type (0) we have $H_{\tilde b}=1$. For type $(V)$ with $\tilde w$ of compact-fiber type, $H_{\tilde b}$ is commensurable with $\langle t_{\tilde w}\rangle$, and $\chi_{\tilde e}(t_{\tilde w})=0$ by Lemma~\ref{lem:slope_vertical_extended}(\ref{lem:slope_vertical_2}).

(\ref{tree_gates_4}): The lines $\ell_{\tilde e}$ and $\ell_{\tilde e'}$ are the axes of $\bar p_e$ and $\bar p_{e'}$, with translation lengths $|p_e|$ and $|p_{e'}|$. Suppose $\operatorname{diam}(\ell_{\tilde e}\cap\ell_{\tilde e'})\ge|p_e|+|p_{e'}|$. By Lemma~\ref{lem:axisintersection}, $\langle\bar p_e\rangle\cap\langle\bar p_{e'}\rangle\neq1$. Since both are maximal cyclic subgroups of a free group, $\langle\bar p_e\rangle=\langle\bar p_{e'}\rangle$, hence $G_{\tilde e}=p_w^{-1}(\langle\bar p_e\rangle)=G_{\tilde e'}$. Thus $G_{\tilde b}=G_{\tilde b'}$ with $\tilde b\neq\tilde b'$, so $G_{\tilde b}\cap G_{\tilde b'}\cong\mathbb Z^2$ is not cyclic, contradicting Proposition~\ref{prop:split-over-tori}(\ref{prop:split-over-tori9}).

Let $D_0':=2\max_e|p_e|$. The projection of a line $\ell'$ onto a line $\ell$ in a tree is a point if $\ell\cap\ell'=\emptyset$, and is $\ell\cap\ell'$ otherwise. Hence $\operatorname{diam}\mathrm{proj}_{\ell_{\tilde e}}(\ell_{\tilde e'})\le D_0'$. By (\ref{tree_gates_1}) there is $z\in\ell_{\tilde e'}\cap\widehat C_{\tilde w}$, and $\mathrm{proj}_{\ell_{\tilde e}}(z)\in\mathrm{proj}_{\ell_{\tilde e}}(\ell_{\tilde e'})\cap Q_{\tilde e}$. Since $\mathrm{proj}_{\ell_{\tilde e}}(\ell_{\tilde e'})$ has diameter at most $D_0'$ and meets $Q_{\tilde e}$, it lies in $\mathcal N_{D_0'}(Q_{\tilde e})$.
\end{proof}

\subsection{Subgroup distortion is determined by the horizontal components}

\begin{defn}[Superadditive closure]\label{defn:superadditive}
For a non-decreasing function $f\colon\mathbb N\to\mathbb N$, its \textit{superadditive closure} is
\[
\overline f(n):=\max\bigl\{f(n_1)+\dots+f(n_\ell)\ \big|\ \ell\ge1,\ n_1+\dots+n_\ell=n,\ n_s\ge1\bigr\}.
\]
Here $n_1,\dots,n_\ell\ge1$. The function $\overline f$ is non-decreasing, $\overline f\ge f$, and $\sum_sf(n_s)\le\overline f(\sum_sn_s)$ for all $n_s\ge1$. We set $\overline f(0):=f(0)$.  We extend functions on $\mathbb N$ to $[0,\infty)$ by $f(x):=f(\lceil x\rceil)$.
\end{defn}

For a horizontal component $T_i$, let
\[
\tilde\delta_i(n):=\sup\bigl\{\widehat d_{\Omega_{T_i}}(a,b)\ \big|\ a,b\in\Omega_{T_i},\ d(a,b)\le n\bigr\}.
\]
By $H$-equivariance, $\tilde\delta_i$ depends only on the $H$-orbit of $T_i$.

\begin{lem}\label{lem:component_distortion}
For every horizontal component $T_i$ there is $C_1\ge1$ such that $\tilde\delta_i(n)\le C_1\,\Delta^G_{H_i}(C_1n+C_1)+C_1$ for all $n$.
\end{lem}

\begin{proof}
Fix $x_0\in\Omega_{T_i}$. By Remark~\ref{rem:remarkonHi} there is $r_0$ with $H_i\cdot B(x_0,r_0)=\Omega_{T_i}$, where the ball is taken with respect to $\widehat d_{\Omega_{T_i}}$, and, as in the proof of Proposition~\ref{prop:lowerboundofrhoninHi}, there is $\lambda_i\ge1$ with $\widehat d_{\Omega_{T_i}}(x_0,gx_0)\le\lambda_i|g|_{H_i}$ for $g\in H_i$. Let $a,b\in\Omega_{T_i}$ with $d(a,b)\le n$, and choose $g,h\in H_i$ with $\widehat d_{\Omega_{T_i}}(a,gx_0)\le r_0$ and $\widehat d_{\Omega_{T_i}}(b,hx_0)\le r_0$. Since $d\le\widehat d_{\Omega_{T_i}}$, we get $d(x_0,g^{-1}hx_0)\le n+2r_0$, so $|g^{-1}h|_G\le\lambda_G(n+2r_0+\epsilon_G)$ by Lemma~\ref{lem:distortionuniversal}, and therefore $|g^{-1}h|_{H_i}\le\Delta^G_{H_i}(\lambda_G(n+2r_0+\epsilon_G))$. Hence $\widehat d_{\Omega_{T_i}}(a,b)\le2r_0+\lambda_i\Delta^G_{H_i}(\lambda_G(n+2r_0+\epsilon_G))$.
\end{proof}

The following result reduces the problem of the distortion of $H$ in $G$ to the problem of the distortion of the stabilizers of its horizontal components. 

\begin{thm} \label{thm:reduction}
Let $T_1,\dots,T_r$ be representatives of the $H$-orbits of horizontal components of $H$ (there are finitely many such orbits, since every horizontal component contains an edge of $T_H$ and $\Lambda_H=H\backslash T_H$ is finite), and let
\[
f(n):=\max\bigl\{n,\ \Delta^G_{H_1}(n),\dots,\Delta^G_{H_r}(n)\bigr\}
\]
(so $f(n)=n$ if $\widehat Y=\emptyset$). Then
\[
f\ \preceq\ \Delta^G_H\ \preceq\ \overline f .
\]
\end{thm}

To prove Theorem~\ref{thm:reduction}, we need two technical lemmas, Lemma~\ref{lem:ordered_path} and Lemma~\ref{lem:local_gates}. We postpone their proofs and first prove Theorem~\ref{thm:reduction}. 

\begin{lem}\label{lem:ordered_path}
There is a constant $K_1\ge1$ with the following property. Let $v_0,v_1,\dots,v_k$ be a geodesic in $T$ with $k\ge2$ and $v_0,v_k$ white, and let $x\in\widetilde M^{v_0}$, $y\in\widetilde M^{v_k}$. For $j$ odd, the vertex $v_j$ is black, and, as in Section~\ref{subsec:axis_notation}, we let $\tilde e_j^-:=(v_j,v_{j-1})$ and $\tilde e_j^+:=(v_j,v_{j+1})$ be the edges by which the geodesic enters and leaves $v_j$, with attaching maps $\iota_j^\pm:=\iota_{\tilde e_j^\pm}\colon\widetilde M^{v_j}\to\widetilde M^{v_{j\pm1}}$. Then there are points $z_j\in\widetilde M^{v_j}$ for $j$ odd, and paths $\gamma'_j\subset\widetilde M^{v_j}$ for $j$ even, such that:
\begin{enumerate}
\item\label{ordered_path_1} $\gamma'_j$ joins $\iota_{j-1}^+(z_{j-1})$ to $\iota_{j+1}^-(z_{j+1})$, where $\iota_{-1}^+(z_{-1}):=x$ and $\iota_{k+1}^-(z_{k+1}):=y$,
\item\label{ordered_path_2} $$\sum_{j\ \mathrm{even}}|\gamma'_j|\le K_1\,d(x,y)+K_1$$ where $|\cdot|=\ell_d$, and $2k\le d(x,y)$.
\end{enumerate}
\end{lem}

In the setting of Lemma~\ref{lem:ordered_path}, assume moreover that $v_0,\dots,v_k$ is a geodesic in $T_H$, $x\in\mathcal H(v_0)$ and $y\in\mathcal H(v_k)$. We set $|\gamma'_{-1}|=|\gamma'_{k+1}|:=0$.

\begin{lem}[Local gates]\label{lem:local_gates}
There is a constant $C_0\ge1$, depending only on $M_\varphi$ and $H$, with the following properties.
\begin{enumerate}
\item\label{local_gates_1} For $j$ odd and $\tilde e=\tilde e_j^\pm$,
\[
d_X\bigl(\pi_X\iota_j^\pm(z_j),\,Q_{\tilde e}\bigr)\le K|\gamma'_{j\pm1}|+D_0' .
\]
\item\label{local_gates_2} Let $j$ be odd and let $v_j$ be of type (0), with $\mathcal H(v_j)=\{\ast_j\}$. Then
\[
d(z_j,\ast_j)\le C_0\bigl(|\gamma'_{j-1}|+|\gamma'_{j+1}|\bigr)+C_0 .
\]
\item\label{local_gates_3} Let $j$ be odd, and suppose that $v_j$ is of type (2), or of type $(V)$ with its compact-fiber neighbour equal to $v_{j-1}$ or $v_{j+1}$. Then there is $\zeta_j\in\mathcal H(v_j)$ with
\[
d(z_j,\zeta_j)\le C_0\bigl(|\gamma'_{j-1}|+|\gamma'_{j+1}|\bigr)+C_0 .
\]
\end{enumerate}
\end{lem}

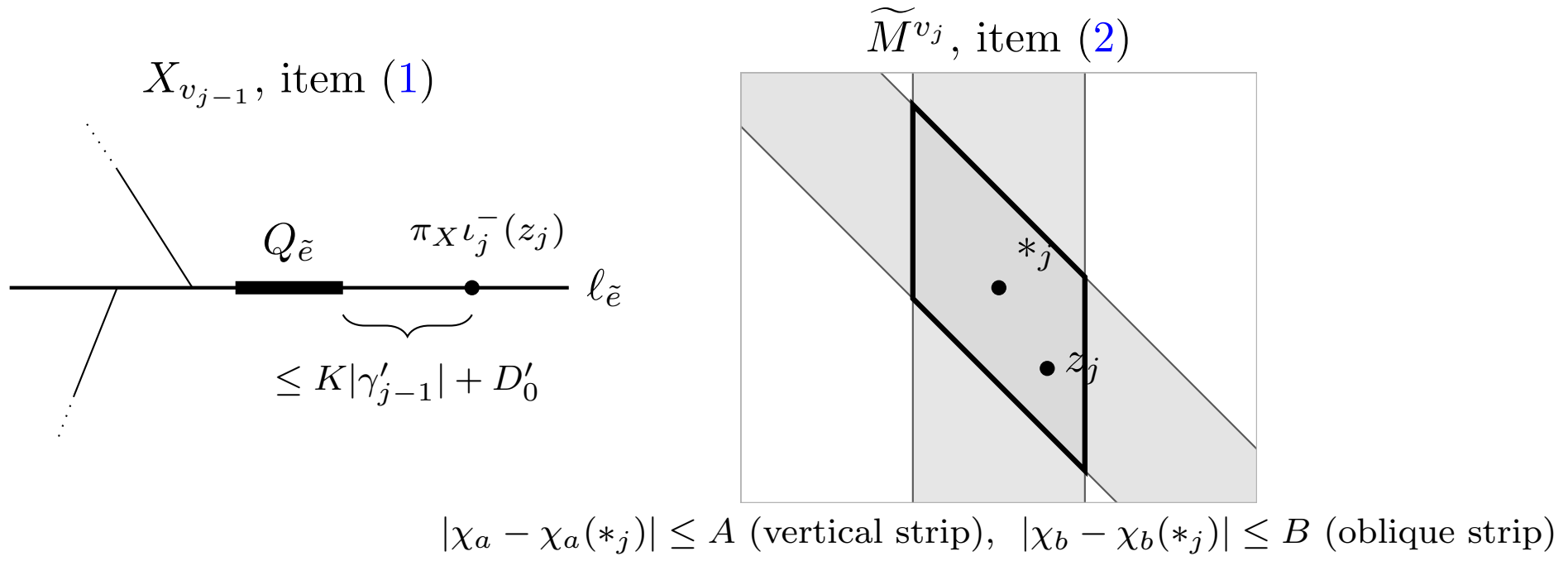
\begin{figure}[htb]
\centering
\begin{tikzpicture}[scale=0.86,>=stealth]
\begin{scope}
\draw[thick] (-2.6,0)--(2.6,0) node[right,font=\small] {$\ell_{\tilde e}$};
\draw (-0.9,0)--(-1.6,1.1); \draw[dotted] (-1.6,1.1)--(-1.9,1.55);
\draw (-1.6,0)--(-2.0,-1.0); \draw[dotted] (-2.0,-1.0)--(-2.15,-1.4);
\draw[line width=3pt] (-0.5,0)--(0.5,0);
\node[font=\small,above] at (0,0.08) {$Q_{\tilde e}$};
\fill (1.7,0) circle (2pt);
\node[font=\scriptsize,above=2pt] at (1.85,0.05) {$\pi_X\iota_j^{-}(z_j)$};
\draw[decorate,decoration={brace,mirror,amplitude=5pt}] (0.5,-0.25)--(1.7,-0.25) node[midway,below=7pt,font=\scriptsize,align=center] {$\le K|\gamma'_{j-1}|+D_0'$};
\node[font=\small] at (0,1.9) {$X_{v_{j-1}}$, item (\ref{local_gates_1})};
\end{scope}
\begin{scope}[shift={(6.6,0)}]
\clip (-2.4,-2.0) rectangle (2.4,2.0);
\fill[gray!20] (-0.8,-2.0) rectangle (0.8,2.0);
\fill[gray!35,opacity=0.6] (-2.9,2.0)--(1.1,-2.0)--(2.9,-2.0)--(-1.1,2.0)--cycle;
\draw[gray!70!black] (-0.8,-2.0)--(-0.8,2.0); \draw[gray!70!black] (0.8,-2.0)--(0.8,2.0);
\draw[gray!70!black] (-1.1,2.0)--(2.9,-2.0); \draw[gray!70!black] (-2.9,2.0)--(1.1,-2.0);
\draw[line width=1.2pt] (0.8,0.1)--(0.8,-1.7)--(-0.8,-0.1)--(-0.8,1.7)--cycle;
\fill (0,0) circle (2pt) node[above right,font=\small] {$\ast_j$};
\fill (0.45,-0.75) circle (2pt) node[right,font=\small] {$z_j$};
\draw[gray!60] (-2.4,-2.0) rectangle (2.4,2.0);
\end{scope}
\node[font=\small] at (6.6,2.35) {$\widetilde M^{v_j}$, item (\ref{local_gates_2})};
\node[font=\scriptsize] at (6.6,-2.3) {$|\chi_a-\chi_a(\ast_j)|\le A$ (vertical strip), \ $|\chi_b-\chi_b(\ast_j)|\le B$ (oblique strip)};
\end{tikzpicture}
\caption{Lemma~\ref{lem:local_gates}. Left: item (\ref{local_gates_1}) for $\tilde e=\tilde e_j^-$ in the tree $X_{v_{j-1}}$. Right: item (\ref{local_gates_2}), where $z_j$ lies in the intersection of two transverse strips around $\ast_j$.}
\label{fig:local_gates}
\end{figure}

Figure~\ref{fig:local_gates} illustrates items (\ref{local_gates_1}) and (\ref{local_gates_2}) of Lemma~\ref{lem:local_gates}.

\begin{proof}[Proof of Theorem~\ref{thm:reduction}]
The lower bound follows from Proposition~\ref{prop:distortion-Hi-H}, the definition of $\preceq$, and the inequality $n\preceq\Delta^G_H(n)$.

Fix $x_0$ in the piece $\mathcal H(\tilde w_0)$ of a white vertex $\tilde w_0\in T_H$. Let $h\in H$ and let $x:=x_0$, $y:=hx_0\in\mathcal H(h\tilde w_0)$. Let $v_0=\tilde w_0,v_1,\dots,v_k=h\tilde w_0$ be the geodesic in $T_H$. If $k=0$, then $$\widehat d_{\Omega_H}(x,y)\le A\,d(x,y)$$ by Proposition~\ref{prop:horizontal pieces}(\ref{prop:horizontal pieces_4}). So assume $k\ge2$, and let $z_j$, $\gamma'_j$ be given by Lemma~\ref{lem:ordered_path}.

Call a black vertex $v_j$ of the geodesic a \textit{gate} if at least one of the edges $\tilde e_j^-$, $\tilde e_j^+$ is not horizontal, and let $j_1<\dots<j_m$ be the indices of the gates. 

By Lemma~\ref{lem:black_types}, a gate is of type $(0)$, of type $(2)$, or of type $(V)$ with its compact-fiber neighbour on the geodesic. For each gate $v_j$ choose a point $\zeta_j\in\mathcal H(v_j)$ as follows: 
\begin{itemize}
    \item $\zeta_j:=\ast_j$ if $v_j$ is of type $(0)$.
    \item Otherwise, $\zeta_j$ as in Lemma~\ref{lem:local_gates}(\ref{local_gates_3}).
\end{itemize} 

In both cases, it follows from Lemma~\ref{lem:local_gates} that
\begin{equation}\label{eq:gate-error}
d(z_j,\zeta_j)\le E_j:=C_0\bigl(|\gamma'_{j-1}|+|\gamma'_{j+1}|\bigr)+C_0 .
\end{equation}
The gates cut the geodesic into segments 
$$S_0=[v_0,v_{j_1-1}] , S_1=[v_{j_1+1},v_{j_2-1}], \dots S_m=[v_{j_m+1},v_k]$$ each starting and ending at a white vertex. If there are no gates, that is, if $m=0$, then there is a single segment $S_0=[v_0,v_k]$. 

Recall from Definition~\ref{defn:horizontal_black} that $\mathcal{B}_H$ is the set of horizontal black vertices of $T_H$, which are of type $(V)$ or $(H)$.
Every black vertex $v_j$ inside a segment has two horizontal edges $\tilde e_j^-$ and $\tilde e_j^+$, so $v_j\in\mathcal B_H$, and these two edges and their endpoints lie in $\widehat Y$, in the same component. Moreover, two consecutive black vertices of a segment share a white neighbour. Hence each segment $S=[v_a,v_b]$ either consists of a single white vertex, or is contained in a single horizontal component $T_{i(S)}$. We refer the reader to Figure~\ref{fig:cut} for an illustration.

\begin{figure}[htb]
\centering
\begin{tikzpicture}[scale=0.82,>=stealth]
\foreach \i in {0,...,14} \coordinate (v\i) at (\i,0);
\draw[hor] (v0)--(v4);
\draw[nhor] (v4)--(v6);
\draw[hor] (v6)--(v7);
\draw[nhor] (v7)--(v11);
\draw[hor] (v11)--(v14);
\foreach \i in {0,2,4,6,12,14} \node[wfib] at (v\i) {};
\foreach \i in {8,10} \node[wcpt] at (v\i) {};
\foreach \i in {1,3,5,7,9,11,13} \node[blk] at (v\i) {};
\foreach \i in {0,...,14} \node[font=\scriptsize,above=5pt] at (v\i) {$v_{\i}$};
\foreach \i/\t in {1/H,3/H,5/0,7/V,9/2,11/V,13/H} \node[typ,below=6pt] at (v\i) {(\t)};
\foreach \i in {5,7,9,11} \draw[gray,thick] ($(v\i)+(-0.2,-0.2)$) rectangle ($(v\i)+(0.2,0.2)$);
\draw[decorate,decoration={brace,mirror,amplitude=5pt}] (-0.2,-0.95)--(4.2,-0.95) node[midway,below=6pt,font=\small] {$S_0\subset T_{i(S_0)}$};
\draw[decorate,decoration={brace,mirror,amplitude=4pt}] (5.8,-0.95)--(6.2,-0.95) node[midway,below=6pt,font=\small] {$S_1$};
\draw[decorate,decoration={brace,mirror,amplitude=4pt}] (7.8,-0.95)--(8.2,-0.95) node[midway,below=6pt,font=\small] {$S_2$};
\draw[decorate,decoration={brace,mirror,amplitude=4pt}] (9.8,-0.95)--(10.2,-0.95) node[midway,below=6pt,font=\small] {$S_3$};
\draw[decorate,decoration={brace,mirror,amplitude=5pt}] (11.8,-0.95)--(14.2,-0.95) node[midway,below=6pt,font=\small] {$S_4\subset T_{i(S_4)}$};
\end{tikzpicture}
\caption{The proof of Theorem~\ref{thm:reduction}. Thick edges are horizontal, and the gates are framed. The gates cut the geodesic into segments, each of which is a single white vertex or lies in one horizontal component.}
\label{fig:cut}
\end{figure}
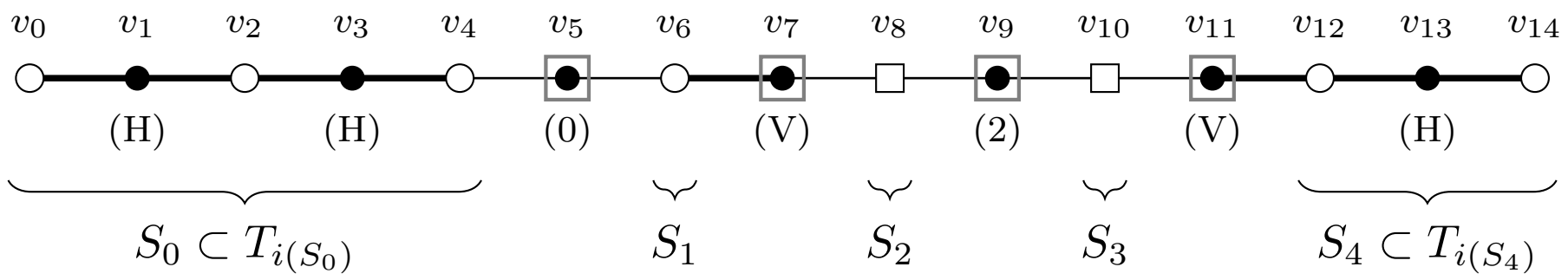

For a segment $S=[v_a,v_b]$ let $a_S:=\iota_{a-1}^+(\zeta_{a-1})$ if $a>0$ and $a_S:=x$ if $a=0$. Similarly, let $b_S:=\iota_{b+1}^-(\zeta_{b+1})$ if $b<k$ and $b_S:=y$ if $b=k$.

By Proposition~\ref{prop:horizontal pieces}(\ref{prop:horizontal pieces_3}), $a_S\in\mathcal H(v_a)$ and $b_S\in\mathcal H(v_b)$.

Let $\sigma_S$ be a path from $a_S$ to $b_S$ chosen as follows: 
\begin{itemize}
    \item if $S$ is a single white vertex,  we choose $\sigma_S$ to be the path of Proposition~\ref{prop:horizontal pieces}(\ref{prop:horizontal pieces_4}) in $\mathcal H(v_a)$, so $|\sigma_S|\le A\,d(a_S,b_S)$,

    \item otherwise, we choose $\sigma_S$ to be a path in $\Omega_{T_{i(S)}}$ with $$|\sigma_S|\le\widehat d_{\Omega_{T_{i(S)}}}(a_S,b_S)+1\le\tilde\delta_{i(S)}(d(a_S,b_S))+1$$
\end{itemize}
 Consecutive segments $S,S'$ separated by the gate $v_j$ are joined in $\Omega_H$ by the vertical segments $t\mapsto\Phi_{\tilde e}(\zeta_j,t)$ in the strips $\mathcal Z(\tilde e_j^-)$ and $\mathcal Z(\tilde e_j^+)$: they go from $b_S=\iota_j^-(\zeta_j)$ to $\zeta_j$ and from $\zeta_j$ to $a_{S'}=\iota_j^+(\zeta_j)$, and have total length at most $4$ by (\ref{eq:vertical-segment}). Hence
\begin{equation}\label{eq:dhat-sum}
\widehat d_{\Omega_H}(x,y)\le\sum_{s=0}^m|\sigma_{S_s}|+4m .
\end{equation}

\textit{Claim:} Define $n_S : = d(a_S, b_S)$. There exists a constant $C_2$ depending only on $M_\varphi$ and $H$ such that 
\[
\sum_{s =0}^{m} n_{S_s} \le C_2\, d(x,y) + C_2
\]

Indeed, for a segment $S=[v_a,v_b]$, follow the path of Lemma~\ref{lem:ordered_path} from $\iota_{a-1}^+(z_{a-1})$ to $\iota_{b+1}^-(z_{b+1})$. It has length at most $\sum_{a\le j\le b,\ j\ \mathrm{even}}|\gamma'_j|+2(b-a)$. If $a=0$, then $a_S=x=\iota_{-1}^+(z_{-1})$, so there is no error at the left end. If $a>0$, then by (\ref{eq:vertical-segment}) and (\ref{eq:gate-error}), we have

\begin{align*}
 d(a_S,\iota_{a-1}^+(z_{a-1})) &\le d(\iota_{a-1}^+(\zeta_{a-1}),\zeta_{a-1})+d(\zeta_{a-1},z_{a-1})+d(z_{a-1},\iota_{a-1}^+(z_{a-1})) \\~\\
 &\le E_{a-1}+4   
\end{align*}
and similarly at the right end, where there is no error if $b=k$. Therefore
\[
n_S\le\sum_{a\le j\le b,\ j\ \mathrm{even}}|\gamma'_j|+2(b-a)+E_{a-1}+E_{b+1}+8 ,
\]
with $E_{-1}=E_{k+1}:=0$. Each $|\gamma'_j|$ occurs in at most two of the numbers $E_{j'}$, and each $E_{j'}$ occurs for at most two segments. Using Lemma~\ref{lem:ordered_path}(\ref{ordered_path_2}) and $m\le k\le d(x,y)/2$, we get
\begin{equation}\label{eq:sum-nS}
\sum_{s=0}^mn_{S_s}\le(1+4C_0)\sum_{j\ \mathrm{even}}|\gamma'_j|+(4C_0+10)(k+1)\le C_2\,d(x,y)+C_2
\end{equation}
for a constant $C_2$ depending only on $M_\varphi$ and $H$.
The claim is proved.

Choose $C_1\ge\max\{1,A\}$, and enlarge it so that it dominates the constants given by Lemma~\ref{lem:component_distortion} for the representatives $T_1,\dots,T_r$, if there are any. Since $\tilde\delta_i$ depends only on the $H$-orbit of $T_i$, and $f\ge\Delta^G_{H_i}$ for $1\le i\le r$, we get $$A\,n\le C_1f(C_1n+C_1)+C_1 \quad \text{and} \quad \tilde\delta_i(n)\le C_1f(C_1n+C_1)+C_1$$ for every horizontal component $T_i$ and all $n\ge0$, after enlarging $C_1$. Since $f$ is non-decreasing, we have $$f(C_1n_S+C_1)\le f(\lceil C_1n_S+C_1\rceil)$$ and $$\sum_s\lceil C_1n_{S_s}+C_1\rceil\le C_1\sum_sn_{S_s}+(C_1+1)(m+1).$$ By (\ref{eq:dhat-sum}), (\ref{eq:sum-nS}) and Definition~\ref{defn:superadditive},
\[
\widehat d_{\Omega_H}(x_0,hx_0)\le C_1\,\overline f\bigl(C_3\,d(x_0,hx_0)+C_3\bigr)+C_3\,d(x_0,hx_0)+C_3
\]
for a constant $C_3$. Finally, $H$ acts on $(\Omega_H,\widehat d_{\Omega_H})$ properly and coboundedly (see Lemma~\ref{lem:basisalmostfiber}), so by Lemma~\ref{lem:svarc-milnor} there are $\lambda_H,\epsilon_H$ with $$|h|_H\le\lambda_H\widehat d_{\Omega_H}(x_0,hx_0)+\epsilon_H$$ and $$d(x_0,hx_0)\le\lambda_G|h|_G+\epsilon_G$$ by Lemma~\ref{lem:distortionuniversal}. Combining these inequalities gives $\Delta^G_H\preceq\overline f$.
\end{proof}

\begin{proof}[Proof of Lemma~\ref{lem:ordered_path}]
Since $\widetilde\pi$ is $1$-Lipschitz and the edges of $T$ have length $2$, we have $2k=d_T(v_0,v_k)\le d(x,y)$. Let $\gamma\colon[0,1]\to\widetilde M_\varphi$ be a geodesic from $x$ to $y$.

Let $\tau_0:=0$, and for $1\le j\le k$ let $\tau_j$ be the first time $\tau\ge\tau_{j-1}$ with $\widetilde\pi\gamma(\tau)=v_j$. It exists because every path in $T$ from $v_{j-1}$ to $v_k$ passes through $v_j$. Let $\tau_{k+1}:=1$, $q_j:=\gamma(\tau_j)$ for $0\le j\le k+1$ (so $q_0=x$, $q_j\in\widetilde M^{v_j}$ for $j\le k$, and $q_{k+1}=y$). Let $l_j:=\ell_d(\gamma|_{[\tau_{j-1},\tau_j]})$. Then
$$l_j:=\ell_d(\gamma|_{[\tau_{j-1},\tau_j]})\ge d(q_{j-1},q_j) \, \text{for} \, 1\le j\le k+1 \quad \text{and} \quad \sum_jl_j=d(x,y)$$

\begin{itemize}
    \item For $j$ odd let $z_j:=q_j$.

    \item For $j$ even with $0<j<k$, let $\gamma'_j$ be the concatenation of a $d_{v_j}$-geodesic from $\iota_{j-1}^+(z_{j-1})$ to $q_j$ and a $d_{v_j}$-geodesic from $q_j$ to $\iota_{j+1}^-(z_{j+1})$.

    \item  Let $\gamma'_0$ be a $d_{v_0}$-geodesic from $x$ to $\iota_1^-(z_1)$, and let $\gamma'_k$ be the concatenation of $d_{v_k}$-geodesics from $\iota_{k-1}^+(z_{k-1})$ to $q_k$ and from $q_k$ to $y$.
\end{itemize}
  By (\ref{eq:Klength}), (\ref{eq:K}) and (\ref{eq:vertical-segment}) of Section~\ref{subsec:metrics}, a $d_{v_j}$-geodesic from $\iota_{j-1}^+(z_{j-1})$ to $q_j$ has $d$-length at most $K^2d(\iota_{j-1}^+(z_{j-1}),q_j)\le K^2(2+l_j)$. Similarly, the other parts have $d$-length at most $K^2(2+l_{j+1})$, and the last part of $\gamma'_k$ has $d$-length at most $K^2l_{k+1}$. Summing over $j$ even gives \[
\sum|\gamma'_j|\le K^2d(x,y)+2K^2(k+1)\le 2K^2d(x,y)+2K^2 \le K_1 d(x,y) + K_1\] with $K_1 : = 2K^2$.
\end{proof}

\begin{proof}[Proof of Lemma~\ref{lem:local_gates}]
(\ref{local_gates_1}): We treat $\tilde e=\tilde e_j^-$, which joins $v_j$ to the white vertex $v_{j-1}$. The case $\tilde e_j^+$ is the same with $\gamma'_{j+1}$ reversed and $y$ in place of $x$. The path $\gamma'_{j-1}$ ends at $\iota_j^-(z_j)$, and $\pi_X\iota_j^-(z_j)$ lies on $\ell_{\tilde e}$. It starts either at $x\in\mathcal H(v_0)$ if $j=1$, whose projection $\pi_X(x)$ lies in $\widehat C_{v_0}$, or at $\iota_{j-2}^+(z_{j-2})$ if $j\ge3$, whose projection lies on $\ell_{\tilde e_{j-2}^+}$. In both cases $\mathrm{proj}_{\ell_{\tilde e}}\pi_X$ of the starting point lies in $\mathcal N_{D_0'}(Q_{\tilde e})$, by definition of $Q_{\tilde e}$ and Lemma~\ref{lem:tree_gates}(\ref{tree_gates_4}). The map $\mathrm{proj}_{\ell_{\tilde e}}\circ\pi_X$ is $1$-Lipschitz for the $\ell^1$-metric $d_{v_{j-1}}$, and $\ell_{d_{v_{j-1}}}(\gamma'_{j-1})\le K|\gamma'_{j-1}|$ by (\ref{eq:Klength}). The claim follows.

(\ref{local_gates_2}): Let $\tilde e_a:=\tilde e_j^-$ and $\tilde e_b:=\tilde e_j^+$, which join $v_j$ to $\tilde w_a:=v_{j-1}$ and $\tilde w_b:=v_{j+1}$, let $e_a,e_b$ be their images in $\Delta$, and write $\chi_a:=\chi_{\tilde e_a}$, $\chi_b:=\chi_{\tilde e_b}$. Since $v_j$ is of type (0), both $\tilde w_a$ and $\tilde w_b$ are of fiber type and $H_{v_j}=1$, so $\operatorname{diam}Q_{\tilde e_a},\operatorname{diam}Q_{\tilde e_b}\le D_0$ by Lemma~\ref{lem:tree_gates}(\ref{tree_gates_3}), and $\pi_X\iota_{\tilde e_a}(\ast_j)\in Q_{\tilde e_a}$, $\pi_X\iota_{\tilde e_b}(\ast_j)\in Q_{\tilde e_b}$ by Lemma~\ref{lem:tree_gates}(\ref{tree_gates_1}). By (\ref{local_gates_1}) and (\ref{eq:chi-distance}),
\[
A:=|\chi_a(z_j)-\chi_a(\ast_j)|\le K|\gamma'_{j-1}|+D_0'+D_0 \] and 

\[B:=|\chi_b(z_j)-\chi_b(\ast_j)|\le K|\gamma'_{j+1}|+D_0'+D_0 .
\]
Write $z_j-\ast_j=(\Delta X,\Delta Y)$ in the coordinates of $\widetilde M^{v_j}$. Then $\Delta X+k(e_a)\Delta Y=\pm A$ and $\Delta X+k(e_b)\Delta Y=\pm B$. Since $|k(e_a)-k(e_b)|\ge1$, we get $|\Delta Y|\le A+B$ and $|\Delta X|\le A+|k(e_a)|(A+B)$. In these coordinates $d_{v_j}$ is the $\ell^1$-metric, so by (\ref{eq:K}),
\[
d(z_j,\ast_j)\le K(2+\kappa_{\max})(A+B),\qquad\kappa_{\max}:=\max_e|k(e)| .
\]

(\ref{local_gates_3}): We are going to define $\zeta_j$ as follows:
\begin{itemize}
    \item If $v_j$ is of type (2), then $\mathcal H(v_j)$, the convex hull of an orbit of a lattice of $\mathbb R^2$, is all of $\widetilde M^{v_j}$, and we define $\zeta_j$ to be the point $z_j$.

    \item Let $v_j$ be of type $(V)$, let $\tilde w:=v_{j\pm1}$ be its compact-fiber neighbour on the geodesic, and let $\tilde e:=\tilde e_j^\pm$ be the edge between them, over $e$. By Lemma~\ref{lem:black_types}, a generator $u$ of $H_{v_j}$ has a power in $\langle t_{\tilde w}\rangle$, so its coordinates are proportional to $(-k(e),1)$ (see Lemma~\ref{lem:slope_vertical}). Hence $\mathcal H(v_j)$ is a line in direction $(-k(e),1)$. It contains exactly one point $\zeta_j$ with the same $Y$-coordinate as $z_j$, and $\chi_{\tilde e}$ is constant on it.
\end{itemize}  
In the latter case, we have
$z_j-\zeta_j=(\Delta X,0)$ with $|\Delta X|=|\chi_{\tilde e}(z_j)-\chi_{\tilde e}(\zeta_j)|$. 

By Lemma~\ref{lem:tree_gates}(\ref{tree_gates_1}),(\ref{tree_gates_3}), we have $\pi_X\iota_{\tilde e}(\zeta_j)\in Q_{\tilde e}$ and $\operatorname{diam}Q_{\tilde e}\le D_0$. Thus (\ref{local_gates_1}) and (\ref{eq:chi-distance}) imply $$|\Delta X|\le K|\gamma'_{j\pm1}|+D_0'+D_0$$ and $$d(z_j,\zeta_j)\le K|\Delta X|$$

Taking $C_0:=K(2+\kappa_{\max})\max\{K,\,2(D_0+D_0')\}$, all three estimates hold.
\end{proof}

\section{Upper bounds and the trichotomy}\label{sec:upper}

In this section, we prove the upper bounds and the main theorem (see Theorem~\ref{thm:main}). The lower bounds were proved in Section~\ref{sec:overlap}.

Throughout this section, $T_i$ is a horizontal component of $H$, and $H_i=\mathrm{Stab}_H(T_i)$. We use the constants $R$ and $L$ of Section~\ref{subsec:horizontalpiece_white}, $A$ of Proposition~\ref{prop:horizontal pieces} and $K$ of (\ref{D3}), and we let $\kappa_{\max}:=\max_e|k(e)|$, where $e$ runs over the edges of $\Delta$. For a white vertex $\tilde w$ of $T_H$, we use the tree $X_{\tilde w}$, the projection $\pi_X$ and the subtree $\widehat C_{\tilde w}=\pi_X(\mathcal H(\tilde w))$ of Sections~\ref{sec:overlap} and~\ref{sec:undistortion}. A \textit{vertical line} of $\widetilde M^{\tilde w}=X_{\tilde w}\times\mathbb R$ is a fibre $\{\tilde x\}\times\mathbb R$ of $\pi_X$.

\subsection{Heights in white vertex spaces}

Let $\tilde w$ be a white vertex of $T_i$. By Lemma~\ref{lem:vertices_of_Ti}, the vertex $\tilde w$ is of fiber type, so $\mathcal H(\tilde w)=g\operatorname{Band}_R(\tilde w_l)$ for some representative $\tilde w_l$ of Section~\ref{subsec:horizontalpiece_white} and some $g\in H$. We write $\tilde w_l$ and $f_l$ for the white representatives $\tilde w_i$ and the functions $f_i$ of Section~\ref{subsec:horizontalpiece_white}, since $i$ now indexes the component $T_i$.

\begin{rem}\label{rem:height_formula}
\begin{enumerate}
    \item Let $\tilde x\in\widehat C_{\tilde w_l}$. By Definition~\ref{defn:band}, the vertical line $\{\tilde x\}\times\mathbb R$ meets $\operatorname{Band}_R(\tilde w_l)$ in the segment $\{\tilde x\}\times[f_l(\tilde x)-R,f_l(\tilde x)+R]$. Since $d_{\tilde w_l}$ restricts to the usual distance on vertical lines, the $d_{\tilde w_l}$-distance from $(\tilde x,s)$ to this segment is
\begin{equation}\label{eq:height_formula}
\max\bigl\{0,\, |s-f_l(\tilde x)|-R\bigr\}.
\end{equation}

\item The element $g$ is a deck transformation, so it maps $\widetilde M^{\tilde w_l}$ to $\widetilde M^{\tilde w}$ isometrically for the product metrics and maps vertical lines to vertical lines. Hence, by item (1), $\mathcal H(\tilde w)$ meets the vertical line over each point of $\widehat C_{\tilde w}$ in a segment of length $2R$.
\end{enumerate}

\end{rem}

\begin{defn} \label{defn:height}
Let $\tilde w$ be a white vertex of $T_i$, and let $p\in\widetilde M^{\tilde w}$ with $\pi_X(p)\in \widehat C_{\tilde w}=\pi_X(\mathcal H(\tilde w))$. The \textit{height} $\mathrm{ht}_{\tilde w}(p)$ of $p$ is the $d_{\tilde w}$-distance from $p$ to the segment (see Remark~\ref{rem:height_formula}) in which $\mathcal H(\tilde w)$ meets the vertical line through $p$.
\end{defn}

\begin{lem}\label{lem:height_white}
Let $\tilde w$ be a white vertex of $T_i$, and let $p,p'\in\widetilde M^{\tilde w}$ with $\pi_X(p),\pi_X(p')\in\widehat C_{\tilde w}$. Then
\begin{enumerate}
\item\label{height_white_1} $\mathrm{ht}_{\tilde w}(p)=0$ if and only if $p\in\mathcal H(\tilde w)$,
\item\label{height_white_2} $|\mathrm{ht}_{\tilde w}(p)-\mathrm{ht}_{\tilde w}(p')|\le(1+L)\,d_{\tilde w}(p,p')$.
\end{enumerate}
\end{lem}

\begin{proof}
Write $\mathcal H(\tilde w)=g\operatorname{Band}_R(\tilde w_l)$ as above. Since $g$ preserves the product metrics and vertical lines, we have $\mathrm{ht}_{\tilde w}(gq)=\mathrm{ht}_{\tilde w_l}(q)$ for all $q$, and $\widehat C_{\tilde w}=\bar g(\widehat C_{\tilde w_l})$, where $\bar g\colon X_{\tilde w_l}\to X_{\tilde w}$ is the isometry of trees induced by $g$ (see (\ref{eq:bar_g})). So we may assume that $\tilde w=\tilde w_l$. In the space $\widetilde M^{\tilde w_l}=\widetilde X_{w_l}\times\mathbb R$, Remark~\ref{rem:height_formula} gives
\[
\mathrm{ht}_{\tilde w_l}(\tilde x,s)=\max\bigl\{0, |s-f_l(\tilde x)|-R\bigr\}\qquad(\tilde x\in\widehat C_{\tilde w_l}).
\]
This gives (\ref{height_white_1}).

The function $f_l$ is $L$-Lipschitz on the subtree $\widehat C_{\tilde w_l}$, which is convex in $\widetilde X_{w_l}$ (see Section~\ref{subsec:horizontalpiece_white}), and $d_{\tilde w_l}$ is the $\ell^1$-product metric of Section~\ref{subsec:metrics}. Hence $(\tilde x,s)\mapsto|s-f_l(\tilde x)|$ is $(1+L)$-Lipschitz, and so is $\mathrm{ht}_{\tilde w_l}$.
\end{proof}

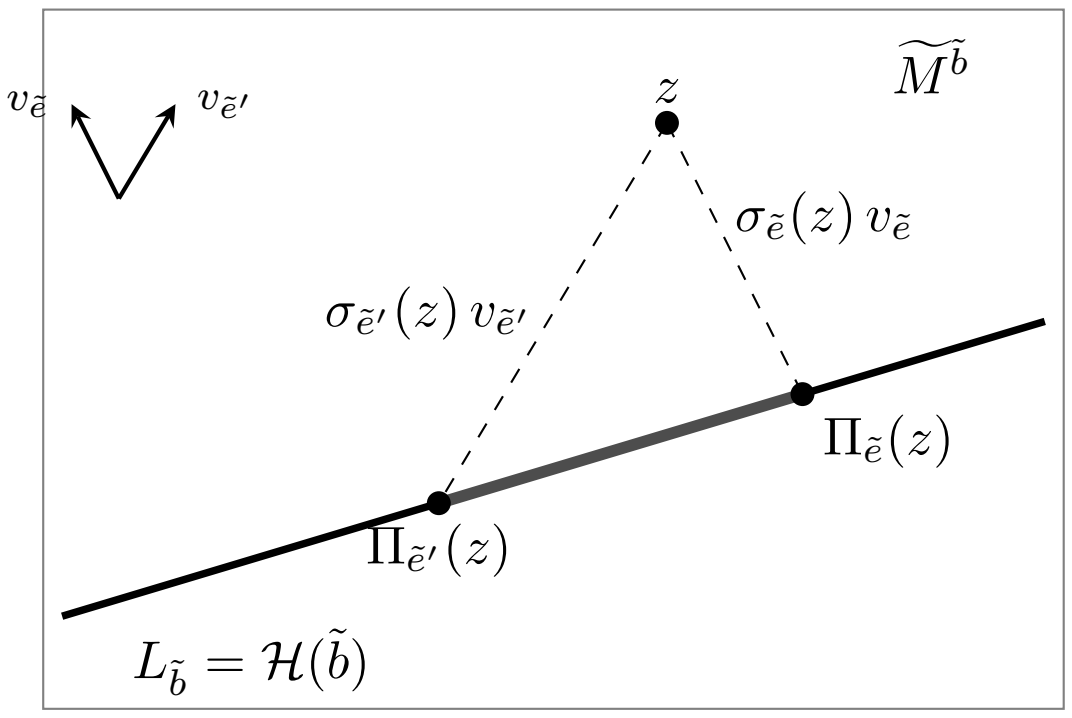
\begin{figure}[htb]
\centering
\begin{tikzpicture}[scale=1.25,>=stealth]
\draw[gray] (-0.8,-0.7) rectangle (4.6,3.0);
\draw[line width=1.4pt] (-0.7,-0.21)--(4.5,1.35);
\node[font=\small] at (0.3,-0.45) {$L_{\tilde b}=\mathcal H(\tilde b)$};
\coordinate (z) at (2.5,2.4);
\coordinate (P1) at (3.217,0.965);
\coordinate (P2) at (1.293,0.388);
\draw[dashed] (z)--(P1);
\draw[dashed] (z)--(P2);
\draw[line width=2.2pt,gray!60!black] (P2)--(P1);
\fill (z) circle (1.8pt) node[above,font=\small] {$z$};
\fill (P1) circle (1.8pt) node[below right,font=\small] {$\Pi_{\tilde e}(z)$};
\fill (P2) circle (1.8pt) node[below,font=\small] {$\Pi_{\tilde e'}(z)$};
\node[font=\small,right] at ($(z)!0.35!(P1)$) {$\sigma_{\tilde e}(z)\,v_{\tilde e}$};
\node[font=\small,left] at ($(z)!0.5!(P2)$) {$\sigma_{\tilde e'}(z)\,v_{\tilde e'}$};
\draw[->,thick] (-0.4,2.0)--(-0.65,2.5) node[left,font=\scriptsize] {$v_{\tilde e}$};
\draw[->,thick] (-0.4,2.0)--(-0.1,2.5) node[right,font=\scriptsize] {$v_{\tilde e'}$};
\node[font=\small] at (3.9,2.7) {$\widetilde M^{\tilde b}$};
\end{tikzpicture}
\caption{Lemma~\ref{lem:height_black} in the plane $\widetilde M^{\tilde b}$. The feet $\Pi_{\tilde e}(z)$ and $\Pi_{\tilde e'}(z)$ of $z$ on $L_{\tilde b}$ in the fibre directions $v_{\tilde e}$ and $v_{\tilde e'}$ of the two white neighbours.}
\label{fig:heights}
\end{figure}

Now let $\tilde b$ be a black vertex of $T_i$, and let $u$ be a generator of $H_{\tilde b}\cong\mathbb Z$, with coordinates $(a,c)$ as in Definition~\ref{defn:coordinates}. We use the coordinates $(X,Y)$ on the plane $\widetilde M^{\tilde b}$ of Section~\ref{subsec:gates_black_lines}, in which $u$ acts as the translation by the vector $(a,c)$. Since $H_{\tilde b}=\langle u\rangle$, the piece $L_{\tilde b}:=\mathcal H(\tilde b)$ is the line through the base point $\ast=(X_\ast,Y_\ast)$ in direction $(a,c)$. Let
\[
\ell_{\tilde b}(X,Y):=a\,(Y-Y_\ast)-c\,(X-X_\ast),
\]
so that $L_{\tilde b}=\{\ell_{\tilde b}=0\}$. For an edge $\tilde e=(\tilde b,\tilde w)$ of $T_i$ over $e$, let $v_{\tilde e}:=(-k(e),1)$ be the coordinate vector of $t_{\tilde w}$ (see Lemma~\ref{lem:slope_vertical}). Then $$\ell_{\tilde b}(z+v_{\tilde e})-\ell_{\tilde b}(z)=a+k(e)c=\chi_{\tilde e}(u)$$ which is nonzero by Lemma~\ref{lem:chi_nonzero}. For $z\in\widetilde M^{\tilde b}$ we define
\[
\sigma_{\tilde e}(z):=\frac{\ell_{\tilde b}(z)}{\chi_{\tilde e}(u)}\in\mathbb R\qquad\text{and}\qquad \Pi_{\tilde e}(z):=z-\sigma_{\tilde e}(z)\,v_{\tilde e} .
\]
Thus $\Pi_{\tilde e}(z)$ is the point where the line through $z$ in the direction of the fibres of $\tilde w$ meets $L_{\tilde b}$ (see Figure~\ref{fig:heights}).

\begin{lem}\label{lem:height_black}
Let $\tilde b$ be a black vertex of $T_i$, let $\tilde e=(\tilde b,\tilde w)$ and $\tilde e'=(\tilde b,\tilde w')$ be edges of $T_i$ at $\tilde b$, and let $z\in\widetilde M^{\tilde b}$. Then:
\begin{enumerate}
\item\label{height_black_1} $\Pi_{\tilde e}(z)\in L_{\tilde b}$, the points $\iota_{\tilde e}(z)$ and $\iota_{\tilde e}(\Pi_{\tilde e}(z))$ lie on the same vertical line of $\widetilde M^{\tilde w}$, at $d_{\tilde w}$-distance $|\sigma_{\tilde e}(z)|$, and $\pi_X\iota_{\tilde e}(z)\in\widehat C_{\tilde w}$,
\item\label{height_black_2} $|\sigma_{\tilde e}(z)|-2R\le\mathrm{ht}_{\tilde w}(\iota_{\tilde e}(z))\le|\sigma_{\tilde e}(z)|$,
\item\label{height_black_3} $|\sigma_{\tilde e'}(z)|=\mathfrak{sl}(\tilde e\to\tilde e')^{-1}\,|\sigma_{\tilde e}(z)|$, and consequently
\[
\mathrm{ht}_{\tilde w'}(\iota_{\tilde e'}(z))\ \le\ \mathfrak{sl}(\tilde e\to\tilde e')^{-1}\bigl(\mathrm{ht}_{\tilde w}(\iota_{\tilde e}(z))+2R\bigr),
\]
\item\label{height_black_4} $d_{\tilde b}\bigl(\Pi_{\tilde e}(z),\Pi_{\tilde e'}(z)\bigr)\le(1+\kappa_{\max})\bigl(|\sigma_{\tilde e}(z)|+|\sigma_{\tilde e'}(z)|\bigr)$.
\end{enumerate}
\end{lem}

\begin{proof}
(\ref{height_black_1}): By definition of $\ell_{\tilde b}$, we have $$\ell_{\tilde b}(\Pi_{\tilde e}(z))=\ell_{\tilde b}(z)-\sigma_{\tilde e}(z)\chi_{\tilde e}(u)=0$$ so $\Pi_{\tilde e}(z)\in L_{\tilde b}$.

By (\ref{eq:affine-attaching}) and (\ref{eq:vertical-attaching}), the $X_{\tilde w}$-component of $\iota_{\tilde e}(X,Y)$ depends only on $X+k(e)Y$, and its $\mathbb R$-component is $Y+c'_{\tilde e}$ for a constant $c'_{\tilde e}$. Moving $z$ by $s\,v_{\tilde e}$ with $s\in\mathbb R$ does not change $X+k(e)Y$ and increases $Y$ by $s$. Hence $\iota_{\tilde e}(z)$ is obtained from $\iota_{\tilde e}(\Pi_{\tilde e}(z))$ by a vertical translation by $\sigma_{\tilde e}(z)$. Finally $\iota_{\tilde e}(\Pi_{\tilde e}(z))\in\iota_{\tilde e}(\mathcal H(\tilde b))\subset\mathcal H(\tilde w)$ by Proposition~\ref{prop:horizontal pieces}(\ref{prop:horizontal pieces_3}), so $$\pi_X\iota_{\tilde e}(z)=\pi_X\iota_{\tilde e}(\Pi_{\tilde e}(z))\in\widehat C_{\tilde w}$$

(\ref{height_black_2}): The point $\iota_{\tilde e}(\Pi_{\tilde e}(z))$ lies in the segment of length $2R$ in which $\mathcal H(\tilde w)$ meets the vertical line through $\iota_{\tilde e}(z)$. By (\ref{height_black_1}), the distance from $\iota_{\tilde e}(z)$ to this point is $|\sigma_{\tilde e}(z)|$, which gives both inequalities.

(\ref{height_black_3}): By definition $$\sigma_{\tilde e'}(z)\chi_{\tilde e'}(u)=\ell_{\tilde b}(z)=\sigma_{\tilde e}(z)\chi_{\tilde e}(u)$$ and $$\mathfrak{sl}(\tilde e\to\tilde e')=|\chi_{\tilde e'}(u)/\chi_{\tilde e}(u)|$$ (see Definition~\ref{defn:slope}). The inequality follows from (\ref{height_black_2}) applied to $\tilde e'$ and to $\tilde e$.

(\ref{height_black_4}): $\Pi_{\tilde e}(z)-\Pi_{\tilde e'}(z)=\sigma_{\tilde e'}(z)v_{\tilde e'}-\sigma_{\tilde e}(z)v_{\tilde e}$, and $d_{\tilde b}$ is the $\ell^1$-metric in our coordinates, in which $v_{\tilde e}$ and $v_{\tilde e'}$ have $\ell^1$-norm at most $1+\kappa_{\max}$.
\end{proof}

\subsection{Upper bounds for the stabilizers}

By Lemma~\ref{lem:slope}(\ref{lem:slope_3}), the slope $\mathfrak{sl}(\tilde e\to\tilde e')$ depends only on the $H_i$-orbit of the triple $(\tilde b,\tilde e,\tilde e')$. The tree $T_i$ has finitely many $H_i$-orbits of vertices (see Lemma~\ref{lem:basifactsofTi}(\ref{lem:basifactsofTi_2})), and every black vertex of $T$ has finitely many edges, one over each edge of $\Delta$ at its image (see Section~\ref{sec:slopes}). Hence the slopes along $T_i$ take finitely many values, and
\[
\Theta_i:=\max\bigl\{\mathfrak{sl}(\tilde e\to\tilde e')\ \big|\ \tilde e,\tilde e'\ \text{edges of}\ T_i\ \text{at a black vertex}\bigr\}
\]
is a real number with $\Theta_i\ge1$ (by Lemma~\ref{lem:slope}(\ref{lem:slope_2}), the maximum of the inverse slopes is also $\Theta_i$).

\begin{lem}\label{lem:bounded_slope_products}
If the dilation homomorphism $\omega_i$ is trivial, then
\[
M_i:=\sup\bigl\{\mathfrak S(\tilde w,\tilde w')\ \big|\ \tilde w,\tilde w'\ \text{white vertices of}\ T_i\bigr\}<\infty 
\] where $\mathfrak S (\tilde w, \tilde w')$ is the slope-product from $\tilde w$ to $\tilde w'$ given by Definition~\ref{defn:slope_product}.
\end{lem}

\begin{proof}
Let $\tilde w^1,\dots,\tilde w^r$ be representatives of the $H_i$-orbits of white vertices of $T_i$, and fix a white vertex $\tilde w_0$ of $T_i$. Let $\tilde w=g\tilde w^a$ and $\tilde w'=g'\tilde w^b$ with $g,g'\in H_i$. By Lemma~\ref{lem:basisPi}, 
\begin{align*}
 \mathfrak S(\tilde w,\tilde w') &=\mathfrak S(g\tilde w^a,g\tilde w_0)\,\mathfrak S(g\tilde w_0,g'\tilde w_0)\,\mathfrak S(g'\tilde w_0,g'\tilde w^b) \\~\\
 &=\mathfrak S(\tilde w^a,\tilde w_0)\,\omega_i(g^{-1}g')\,\mathfrak S(\tilde w_0,\tilde w^b),   
\end{align*}
where we used Lemma~\ref{lem:basisPi}(\ref{lem:basisPi_3}) and Proposition~\ref{prop:dilation_basis}(\ref{prop:dilation_basis_2}) to get $\mathfrak S(g\tilde w_0,g'\tilde w_0)=\mathfrak S(\tilde w_0,g^{-1}g'\tilde w_0)=\omega_i(g^{-1}g')$. If $\omega_i$ is trivial, this takes finitely many values.
\end{proof}

\begin{prop}\label{prop:main_estimate}
There is a constant $C_5\ge1$, depending only on $M_\varphi$, $H$ and $T_i$, with the following property. Let $\tilde w,\tilde w'$ be white vertices of $T_i$, let $x\in\mathcal H(\tilde w)$ and $y\in\mathcal H(\tilde w')$, and let $d:=d(x,y)$. Let $v_0=\tilde w,v_1,\dots,v_k=\tilde w'$ be the geodesic in $T_i$, and let
\[
\mathfrak m:=\max\bigl\{\mathfrak S(v_a,v_b)^{-1}\ \big|\ 0\le a\le b\le k,\ a,b\ \text{even}\bigr\}\ \ge1 .
\]
Then
\[
\widehat d_{\Omega_{T_i}}(x,y)\ \le\ C_5\,\mathfrak m\,(d+1)^2 .
\]
Moreover $\mathfrak m\le\Theta_i^{\,d/4}$, and $\mathfrak m\le M_i$ if $\omega_i$ is trivial.
\end{prop}

\begin{proof}
The last sentence holds because $\mathfrak S(v_a,v_b)^{-1}$ is a product of at most $k/2$ inverse slopes, each at most $\Theta_i$, because $2k\le d$ (see Lemma~\ref{lem:ordered_path}(\ref{ordered_path_2})), and because $\mathfrak S(v_a,v_b)^{-1}=\mathfrak S(v_b,v_a)\le M_i$ (see Lemma~\ref{lem:basisPi}(\ref{lem:basisPi_1})). Here $M_i<\infty$ by Lemma~\ref{lem:bounded_slope_products}.

If $k=0$, then Proposition~\ref{prop:horizontal pieces}(\ref{prop:horizontal pieces_4}) gives a path in $\mathcal H(\tilde w)\subset\Omega_{T_i}$ from $x$ to $y$ of length at most $A\,d$. So for the rest of the proof, we  assume that $k\ge2$.

The geodesic $v_0,\dots,v_k$ of $T$ lies in the subtree $T_i$, so its vertices of odd index are black vertices of $T_i$, its vertices of even index are white vertices of $T_i$, and its edges are edges of $T_i$. Let $z_j\in\widetilde M^{v_j}$ ($j$ odd) and $\gamma'_j\subset\widetilde M^{v_j}$ ($j$ even) be given by Lemma~\ref{lem:ordered_path}, so that
\begin{equation}\label{eq:sum_gamma}
\sum_{j\ \mathrm{even}}|\gamma'_j|\le K_1d+K_1\qquad\text{and}\qquad 2k\le d .
\end{equation}

 For $j$ odd, let $\tilde e_j^\pm$ and $\iota_j^\pm$ be as in Lemma~\ref{lem:ordered_path}, and let
\begin{itemize}
    \item $\mathrm{ht}^{\mathrm{in}}_j:=\mathrm{ht}_{v_{j-1}}\bigl(\iota_j^-(z_j)\bigr)$

    \item $\mathrm{ht}^{\mathrm{out}}_j:=\mathrm{ht}_{v_{j+1}}\bigl(\iota_j^+(z_j)\bigr)$

    \item $ \mathfrak s_j:=\mathfrak{sl}\bigl(\tilde e_j^-\to\tilde e_j^+\bigr).$
\end{itemize}
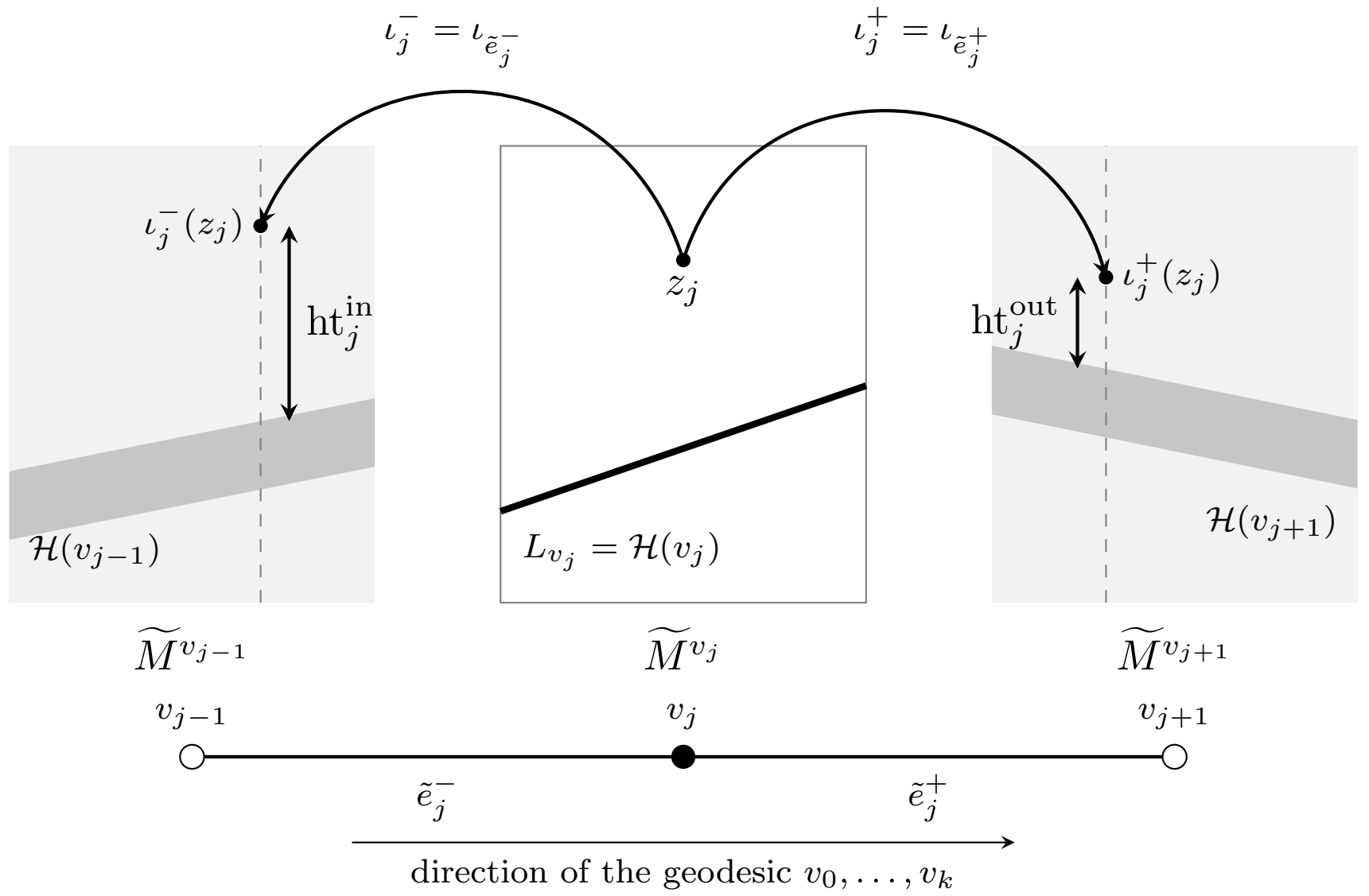
\begin{figure}[htb]
\centering
\begin{tikzpicture}[scale=0.95,>=stealth]
\fill[gray!10] (0,0) rectangle (3.2,4);
\fill[gray!45] (0,0.55)--(3.2,1.19)--(3.2,1.79)--(0,1.15)--cycle;
\node[font=\scriptsize] at (0.75,0.45) {$\mathcal H(v_{j-1})$};
\draw[gray,dashed] (2.2,0)--(2.2,4);
\coordinate (A) at (2.2,3.3);
\fill (A) circle (1.8pt) node[left,font=\scriptsize] {$\iota_j^-(z_j)$};
\draw[<->,thick] (2.45,1.59)--(2.45,3.3) node[midway,right,font=\small] {$\mathrm{ht}^{\mathrm{in}}_j$};
\node[font=\small] at (1.6,-0.4) {$\widetilde M^{v_{j-1}}$};
\draw[gray] (4.3,0) rectangle (7.5,4);
\draw[line width=1.6pt] (4.3,0.8)--(7.5,1.9);
\node[font=\scriptsize,anchor=west] at (4.35,0.45) {$L_{v_j}=\mathcal H(v_j)$};
\coordinate (Z) at (5.9,3.0);
\fill (Z) circle (1.8pt) node[below,font=\small] {$z_j$};
\node[font=\small] at (5.9,-0.4) {$\widetilde M^{v_j}$};
\fill[gray!10] (8.6,0) rectangle (11.8,4);
\fill[gray!45] (8.6,1.65)--(11.8,1.0)--(11.8,1.6)--(8.6,2.25)--cycle;
\node[font=\scriptsize] at (11.05,0.7) {$\mathcal H(v_{j+1})$};
\draw[gray,dashed] (9.6,0)--(9.6,4);
\coordinate (B) at (9.6,2.85);
\fill (B) circle (1.8pt) node[right,font=\scriptsize] {$\iota_j^+(z_j)$};
\draw[<->,thick] (9.35,2.05)--(9.35,2.85) node[midway,left,font=\small] {$\mathrm{ht}^{\mathrm{out}}_j$};
\node[font=\small] at (10.2,-0.4) {$\widetilde M^{v_{j+1}}$};
\draw[->,thick] (Z) .. controls (5.3,4.9) and (2.9,4.9) .. (A);
\node[font=\scriptsize] at (3.9,4.95) {$\iota_j^-=\iota_{\tilde e_j^-}$};
\draw[->,thick] (Z) .. controls (6.5,4.9) and (9.1,4.6) .. (B);
\node[font=\scriptsize] at (8.0,4.95) {$\iota_j^+=\iota_{\tilde e_j^+}$};
\begin{scope}[shift={(0,-1.35)}]
\draw[thick] (1.6,0)--(5.9,0) node[midway,below,font=\scriptsize] {$\tilde e_j^-$};
\draw[thick] (5.9,0)--(10.2,0) node[midway,below,font=\scriptsize] {$\tilde e_j^+$};
\fill[white,draw=black] (1.6,0) circle (3pt); \fill (5.9,0) circle (3pt); \fill[white,draw=black] (10.2,0) circle (3pt);
\node[font=\scriptsize,above] at (1.6,0.1) {$v_{j-1}$};
\node[font=\scriptsize,above] at (5.9,0.1) {$v_j$};
\node[font=\scriptsize,above] at (10.2,0.1) {$v_{j+1}$};
\draw[->] (3.0,-0.75)--(8.8,-0.75) node[midway,below,font=\scriptsize] {direction of the geodesic $v_0,\dots,v_k$};
\end{scope}
\end{tikzpicture}
\caption{The heights $\mathrm{ht}^{\mathrm{in}}_j$ and $\mathrm{ht}^{\mathrm{out}}_j$ of the points $\iota_j^-(z_j)$ and $\iota_j^+(z_j)$ in the two white vertex spaces adjacent to the black vertex $v_j$.}
\label{fig:ht_in_out}
\end{figure}

These heights are defined by Lemma~\ref{lem:height_black}(\ref{height_black_1}).  We refer the reader to Figure~\ref{fig:ht_in_out} for an illustration of $\mathrm{ht}^{\mathrm{out}}_j$ and $\mathrm{ht}^{\mathrm{in}}_j$.

We claim that
\begin{equation}\label{eq:height_white_step}
\mathrm{ht}^{\mathrm{in}}_1\le(1+L)K|\gamma'_0|,\quad \mathrm{ht}^{\mathrm{in}}_{j+1}\le\mathrm{ht}^{\mathrm{out}}_{j-1}+(1+L)K|\gamma'_j|\quad(j\ \text{even},\ 2\le j\le k-2),
\end{equation} and
\begin{equation}\label{eq:height_black_step}
\mathrm{ht}^{\mathrm{out}}_j\le \mathfrak s_j^{-1}(\mathrm{ht}^{\mathrm{in}}_j+2R)\quad(j\ \text{odd}).
\end{equation}
Indeed, let $p:=\iota_{j-1}^+(z_{j-1})$ if $j\ge2$ and $p:=x$ if $j=0$, and $p':=\iota_{j+1}^-(z_{j+1})$. These are the endpoints of the path $\gamma'_j\subset\widetilde M^{v_j}$ (see Lemma~\ref{lem:ordered_path}(\ref{ordered_path_1})). By definition, $\mathrm{ht}_{v_j}(p)=\mathrm{ht}^{\mathrm{out}}_{j-1}$ if $j\ge2$, and $\mathrm{ht}_{v_j}(p')=\mathrm{ht}^{\mathrm{in}}_{j+1}$. If $j=0$, then $\mathrm{ht}_{v_0}(x)=0$ by Lemma~\ref{lem:height_white}(\ref{height_white_1}), since $x\in\mathcal H(v_0)$. The projections $\pi_X(p),\pi_X(p')$ lie in $\widehat C_{v_j}$ by Lemma~\ref{lem:height_black}(\ref{height_black_1}) (and since $x\in\mathcal H(v_0)$). By (\ref{eq:Klength}), $d_{v_j}(p,p')\le\ell_{d_{v_j}}(\gamma'_j)\le K|\gamma'_j|$. Hence Lemma~\ref{lem:height_white}(\ref{height_white_2}) gives $\mathrm{ht}_{v_j}(p')\le\mathrm{ht}_{v_j}(p)+(1+L)K|\gamma'_j|$, which is (\ref{eq:height_white_step}).

Inequality (\ref{eq:height_black_step}) is Lemma~\ref{lem:height_black}(\ref{height_black_3}) at $\tilde b=v_j$. The claim is confirmed.

 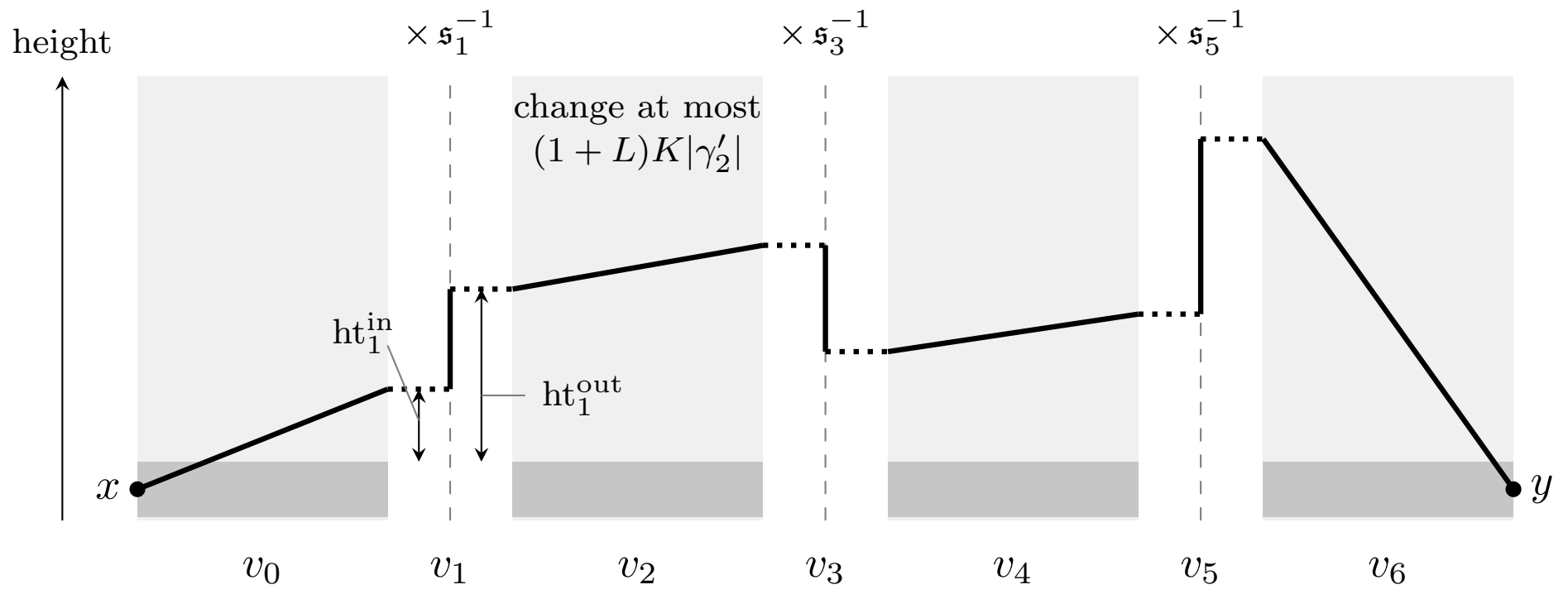
\begin{figure}[htb]
\centering
\begin{tikzpicture}[scale=1,>=stealth]
\foreach \x/\lab in {0/0,3/2,6/4,9/6}{
  \fill[gray!12] (\x,-0.25) rectangle (\x+2,3.3);
  \fill[gray!45] (\x,-0.22) rectangle (\x+2,0.22);
  \node[font=\small] at (\x+1,-0.65) {$v_{\lab}$};
}
\foreach \x/\lab in {2.5/1,5.5/3,8.5/5}{
  \draw[gray,dashed] (\x,-0.25)--(\x,3.3);
  \node[font=\small] at (\x,-0.65) {$v_{\lab}$};
}
\draw[line width=1.2pt] (0,0)--(2,0.8);
\draw[line width=1.2pt,dotted] (2,0.8)--(2.5,0.8);
\draw[line width=1.2pt] (2.5,0.8)--(2.5,1.6);
\draw[line width=1.2pt,dotted] (2.5,1.6)--(3,1.6);
\draw[line width=1.2pt] (3,1.6)--(5,1.95);
\draw[line width=1.2pt,dotted] (5,1.95)--(5.5,1.95);
\draw[line width=1.2pt] (5.5,1.95)--(5.5,1.1);
\draw[line width=1.2pt,dotted] (5.5,1.1)--(6,1.1);
\draw[line width=1.2pt] (6,1.1)--(8,1.4);
\draw[line width=1.2pt,dotted] (8,1.4)--(8.5,1.4);
\draw[line width=1.2pt] (8.5,1.4)--(8.5,2.8);
\draw[line width=1.2pt,dotted] (8.5,2.8)--(9,2.8);
\draw[line width=1.2pt] (9,2.8)--(11,0);
\fill (0,0) circle (1.8pt) node[left,font=\small] {$x$};
\fill (11,0) circle (1.8pt) node[right,font=\small] {$y$};
\draw[<->] (2.25,0.22)--(2.25,0.8);
\node[font=\scriptsize,anchor=east] at (2.2,1.25) {$\mathrm{ht}^{\mathrm{in}}_1$};
\draw[gray] (2.25,0.55)--(2.0,1.15);
\draw[<->] (2.75,0.22)--(2.75,1.6);
\node[font=\scriptsize,anchor=west] at (3.1,0.75) {$\mathrm{ht}^{\mathrm{out}}_1$};
\draw[gray] (2.75,0.75)--(3.1,0.75);
\node[font=\scriptsize,align=center] at (2.5,3.65) {$\times\,\mathfrak s_1^{-1}$};
\node[font=\scriptsize,align=center] at (5.5,3.65) {$\times\,\mathfrak s_3^{-1}$};
\node[font=\scriptsize,align=center] at (8.5,3.65) {$\times\,\mathfrak s_5^{-1}$};
\node[font=\scriptsize,align=center] at (4,2.85) {change at most\\ $(1+L)K|\gamma'_2|$};
\draw[->] (-0.6,-0.25)--(-0.6,3.3) node[above,font=\scriptsize] {height};
\end{tikzpicture}
\caption{The heights along the geodesic $v_0,\dots,v_6$. Inside a white vertex space the height changes by at most $(1+L)K|\gamma'_j|$ (see (\ref{eq:height_white_step})). At a black vertex $v_j$ it is multiplied by at most $\mathfrak s_j^{-1}$, up to an additive error (see (\ref{eq:height_black_step})).}
\label{fig:height_profile}
\end{figure}

 For $0\le a\le b\le k$ with $a,b$ even, the black vertices of the geodesic $[v_a,v_b]$ are the $v_l$ with $a<l<b$ odd, entered through $\tilde e_l^-$ and left through $\tilde e_l^+$. Hence 
\[
\mathfrak S(v_a,v_b)=\prod_{a<l<b,\ l\ \mathrm{odd}}\mathfrak s_l \](see Definition~\ref{defn:slope_product}). Iterating inequalities (\ref{eq:height_white_step}) and (\ref{eq:height_black_step}), we get for $j$ odd
\begin{align*}
    \mathrm{ht}^{\mathrm{out}}_j &\le  \sum_{\substack{i<j\\ i \,\mathrm{even}}}(1+L)K|\gamma'_i|\,\mathfrak S(v_i,v_{j+1})^{-1} + \sum_{\substack{i\le j\\ i\ \mathrm{odd}}}2R\,\mathfrak S(v_{i-1},v_{j+1})^{-1} \\~\\
    &\le\ \mathfrak m\Bigl((1+L)K\sum_{i\ \mathrm{even}}|\gamma'_i|+2R(k+1)\Bigr).
\end{align*}

Denote by $\Xi$ the right-hand side of the last inequality, so that $\mathrm{ht}^{\mathrm{out}}_j\le\Xi$ for every odd $j$. Since $\mathfrak m\ge1$, we have $(1+L)K|\gamma'_i|\le\Xi$ for every even $i$. Hence (\ref{eq:height_white_step}) gives $\mathrm{ht}^{\mathrm{in}}_1\le\Xi$ and
\[
\mathrm{ht}^{\mathrm{in}}_j\le\mathrm{ht}^{\mathrm{out}}_{j-2}+(1+L)K|\gamma'_{j-1}|\le2\,\Xi\qquad(j\ \text{odd},\ j\ge3).
\]
By (\ref{eq:sum_gamma}), there is a constant $C_6$, depending only on $M_\varphi$, $H$ and $T_i$, such that
\begin{equation}\label{eq:height_bound}
\mathrm{ht}^{\mathrm{in}}_j,\ \mathrm{ht}^{\mathrm{out}}_j\ \le\ C_6\,\mathfrak m\,(d+1)\qquad\text{for all odd }j .
\end{equation}

 For $j$ odd let
\[
p_j^{\mathrm{in}}:=\Pi_{\tilde e_j^-}(z_j)\in L_{v_j}\qquad\text{and}\qquad p_j^{\mathrm{out}}:=\Pi_{\tilde e_j^+}(z_j)\in L_{v_j}.
\]

For $j$ even let $a_j:=\iota_{j-1}^+(p_{j-1}^{\mathrm{out}})$ if $j\ge2$ and $a_0:=x$, and let $b_j:=\iota_{j+1}^-(p_{j+1}^{\mathrm{in}})$ if $j\le k-2$ and $b_k:=y$. By Proposition~\ref{prop:horizontal pieces}(\ref{prop:horizontal pieces_3}), $a_j,b_j\in\mathcal H(v_j)$ (see Figure~\ref{fig:model_path}). We concatenate the following paths, in the order of the geodesic:
\begin{itemize}
\item for $j$ even, a path in $\mathcal H(v_j)$ from $a_j$ to $b_j$ of length at most $A\,d(a_j,b_j)$ (see Proposition~\ref{prop:horizontal pieces}(\ref{prop:horizontal pieces_4})),
\item for $j$ odd, the vertical segment in the strip $\mathcal Z(\tilde e_j^-)$ from $b_{j-1}=\iota_j^-(p_j^{\mathrm{in}})$ to $p_j^{\mathrm{in}}$, then the straight segment of the line $L_{v_j}$ from $p_j^{\mathrm{in}}$ to $p_j^{\mathrm{out}}$, then the vertical segment in the strip $\mathcal Z(\tilde e_j^+)$ from $p_j^{\mathrm{out}}$ to $a_{j+1}=\iota_j^+(p_j^{\mathrm{out}})$.
\end{itemize}

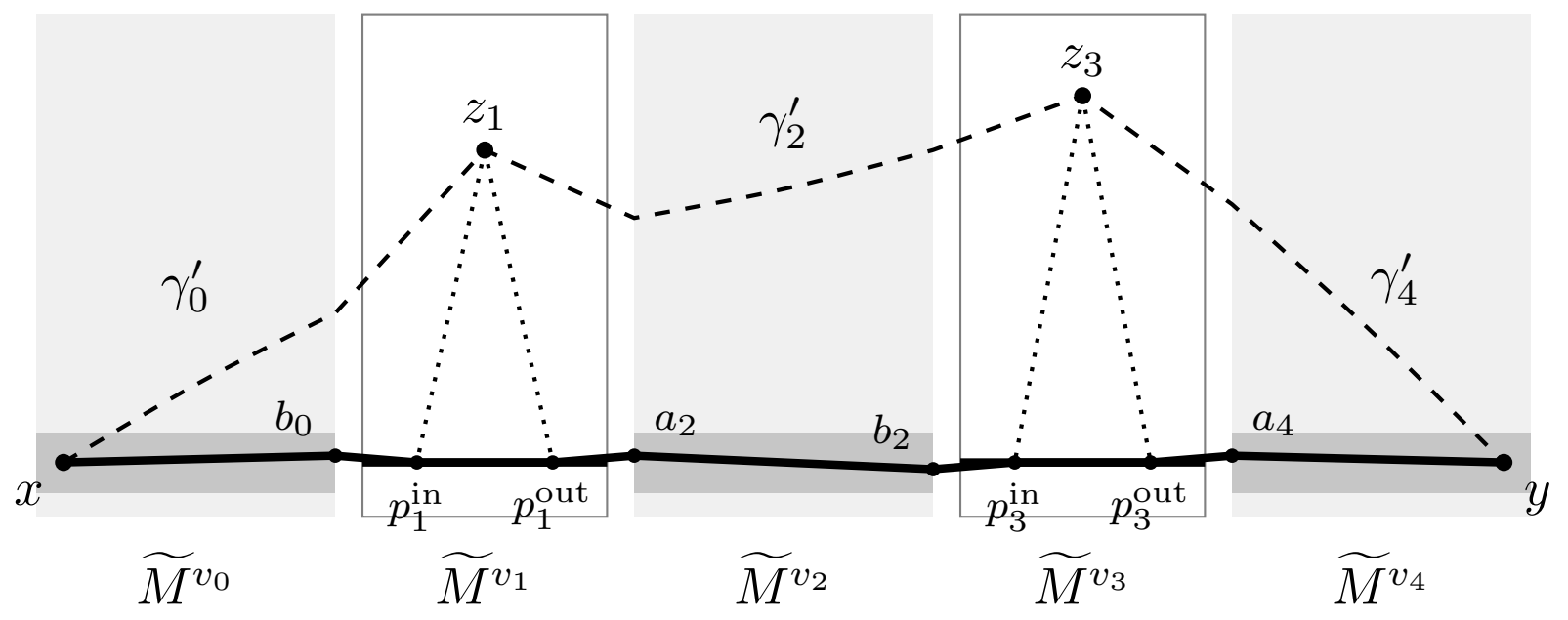
\begin{figure}[htb]
\centering
\begin{tikzpicture}[scale=0.92,>=stealth]
\foreach \x/\lab in {0/0,4.4/2,8.8/4}{
  \fill[gray!12] (\x,-0.4) rectangle (\x+2.2,3.3);
  \fill[gray!45] (\x,-0.22) rectangle (\x+2.2,0.22);
  \node[font=\small] at (\x+1.1,-0.85) {$\widetilde M^{v_{\lab}}$};
}
\foreach \x/\lab in {3.3/1,7.7/3}{
  \draw[gray] (\x-0.9,-0.4) rectangle (\x+0.9,3.3);
  \draw[line width=1.6pt] (\x-0.9,0)--(\x+0.9,0);
  \node[font=\small] at (\x,-0.85) {$\widetilde M^{v_{\lab}}$};
}
\draw[dashed,thick] (0.2,0) .. controls (1.2,0.6) .. (2.2,1.1)
   -- (3.3,2.3) -- (4.4,1.8) .. controls (5.5,2.0) .. (6.6,2.3)
   -- (7.7,2.7) -- (8.8,1.9) .. controls (9.8,1.0) .. (10.8,0);
\fill (3.3,2.3) circle (1.8pt) node[above,font=\small] {$z_1$};
\fill (7.7,2.7) circle (1.8pt) node[above,font=\small] {$z_3$};
\node[font=\small] at (1.1,1.3) {$\gamma'_0$};
\node[font=\small] at (5.5,2.5) {$\gamma'_2$};
\node[font=\small] at (10,1.35) {$\gamma'_4$};
\foreach \x/\z/\a/\b in {3.3/2.3/2.8/3.8,7.7/2.7/7.2/8.2}{
  \draw[dotted,thick] (\x,\z)--(\a,0);
  \draw[dotted,thick] (\x,\z)--(\b,0);
}
\draw[line width=1.6pt] (0.2,0)--(2.2,0.05) -- (2.8,0) -- (3.8,0) -- (4.4,0.05) -- (6.6,-0.05) -- (7.2,0) -- (8.2,0) -- (8.8,0.05) -- (10.8,0);
\fill (0.2,0) circle (1.8pt) node[below left,font=\small] {$x$};
\fill (10.8,0) circle (1.8pt) node[below right,font=\small] {$y$};
\fill (2.8,0) circle (1.5pt) node[below,font=\scriptsize] {$p_1^{\mathrm{in}}$};
\fill (3.8,0) circle (1.5pt) node[below,font=\scriptsize] {$p_1^{\mathrm{out}}$};
\fill (7.2,0) circle (1.5pt) node[below,font=\scriptsize] {$p_3^{\mathrm{in}}$};
\fill (8.2,0) circle (1.5pt) node[below,font=\scriptsize] {$p_3^{\mathrm{out}}$};
\fill (2.2,0.05) circle (1.5pt) node[above left,font=\scriptsize] {$b_0$};
\fill (4.4,0.05) circle (1.5pt) node[above right,font=\scriptsize] {$a_2$};
\fill (6.6,-0.05) circle (1.5pt) node[above left,font=\scriptsize] {$b_2$};
\fill (8.8,0.05) circle (1.5pt) node[above right,font=\scriptsize] {$a_4$};
\end{tikzpicture}
\caption{The new path (thick) in $\Omega_{T_i}$ for $k=4$, compared with the ambient path (dashed). It crosses each black plane along $L_{v_j}$ from $p_j^{\mathrm{in}}$ to $p_j^{\mathrm{out}}$.}
\label{fig:model_path}
\end{figure}

All these paths lie in $\Omega_{T_i}$. The vertical segments have length at most $2$ by (\ref{eq:vertical-segment}). By (\ref{eq:Klength}), Lemma~\ref{lem:height_black}(\ref{height_black_4}) and Lemma~\ref{lem:height_black}(\ref{height_black_2}), the straight segment from $p_j^{\mathrm{in}}$ to $p_j^{\mathrm{out}}$ has length at most
\[
K(1+\kappa_{\max})\bigl(\mathrm{ht}^{\mathrm{in}}_j+\mathrm{ht}^{\mathrm{out}}_j+4R\bigr).
\]
For $j$ even, Lemma~\ref{lem:height_black}(\ref{height_black_1}), Lemma~\ref{lem:height_black}(\ref{height_black_2}) and (\ref{eq:K}) give $$d(a_j,\iota_{j-1}^+(z_{j-1}))\le K(\mathrm{ht}^{\mathrm{out}}_{j-1}+2R)$$ and $$d(\iota_{j+1}^-(z_{j+1}),b_j)\le K(\mathrm{ht}^{\mathrm{in}}_{j+1}+2R)$$ with the convention that these terms vanish for $j=0$ and $j=k$ respectively. Since $\gamma'_j$ joins $\iota_{j-1}^+(z_{j-1})$ to $\iota_{j+1}^-(z_{j+1})$, we have
\[
d(a_j,b_j)\le|\gamma'_j|+K\bigl(\mathrm{ht}^{\mathrm{out}}_{j-1}+\mathrm{ht}^{\mathrm{in}}_{j+1}+4R\bigr).
\]
Summing all contributions, and using (\ref{eq:sum_gamma}) and (\ref{eq:height_bound}) together with the fact that there are at most $k\le d/2$ odd indices, we obtain
\begin{align*}
    \widehat d_{\Omega_{T_i}}(x,y) &\le A\sum_{j\ \mathrm{even}}|\gamma'_j|+C_7(k+1)\\
    &\quad+\bigl(AK+K(1+\kappa_{\max})\bigr)\sum_{j\ \mathrm{odd}}(\mathrm{ht}^{\mathrm{in}}_j+\mathrm{ht}^{\mathrm{out}}_j)\\
    &\le C_5\,\mathfrak m\,(d+1)^2
\end{align*}
for constants $C_7$ and $C_5$ depending only on $M_\varphi$, $H$ and $T_i$.
\end{proof}

\begin{thm}[Upper bounds for the stabilizers]\label{thm:upper_Hi}
Let $T_i$ be a horizontal component of $H$ and $H_i=\mathrm{Stab}_H(T_i)$.
\begin{enumerate}
\item\label{upper_Hi_1} If the dilation homomorphism $\omega_i$ is trivial, then $\Delta^G_{H_i}(n)\preceq n^2$.
\item\label{upper_Hi_2} In general, $\Delta^G_{H_i}(n)\preceq 2^n$.
\end{enumerate}
\end{thm}

\begin{proof}
Fix a white vertex $\tilde w_0$ of $T_i$ and $x_0\in\mathcal H(\tilde w_0)$. Let $h\in H_i$. Then $hx_0\in\mathcal H(h\tilde w_0)$ and $h\tilde w_0$ is a white vertex of $T_i$. By Remark~\ref{rem:remarkonHi} and Lemma~\ref{lem:svarc-milnor}, there is a finite generating set of $H_i$ with respect to which $|h|_{H_i}\le\widehat d_{\Omega_{T_i}}(x_0,hx_0)+2$. By Lemma~\ref{lem:distortionuniversal}, $d:=d(x_0,hx_0)\le\lambda_G|h|_G+\epsilon_G$. Proposition~\ref{prop:main_estimate} gives
\[
|h|_{H_i}\le C_5\,\mathfrak m\,(d+1)^2+2 ,
\]
with $\mathfrak m\le M_i<\infty$ if $\omega_i$ is trivial (see Lemma~\ref{lem:bounded_slope_products}), and $\mathfrak m\le\Theta_i^{\,d/4}$ in general.

If $|h|_G\le n$, we get $$|h|_{H_i}\le C_5M_i(\lambda_Gn+\epsilon_G+1)^2+2$$ in the first case, and $$|h|_{H_i}\le C_5(\lambda_Gn+\epsilon_G+1)^2\Theta_i^{(\lambda_Gn+\epsilon_G)/4}+2$$ in the second case. The first bound is $\preceq n^2$. The second is at most $C\,2^{Cn+C}$ for a suitable constant $C$, hence $\preceq 2^n$.
\end{proof}

\subsection{The main theorem}

We use the following elementary facts about superadditive closures.

\begin{lem}\label{lem:closure_powers}
Let $f\colon\mathbb N\to\mathbb N$ be non-decreasing.
\begin{enumerate}
\item\label{closure_powers_1} If $f(n)\le C(n+1)^2$ for all $n$, then $\overline f(n)\le4Cn^2$ for all $n\ge1$.
\item\label{closure_powers_2} If $f(n)\le C\,a^n$ for all $n$, with $a\ge2$, then $\overline f(n)\le C\,a^n$ for all $n\ge1$.
\end{enumerate}
\end{lem}

\begin{proof}
Let $n_1+\dots+n_\ell=n$ with $n_s\ge1$. For (\ref{closure_powers_1}), we have $$\sum_s(n_s+1)^2\le\bigl(\sum_s(n_s+1)\bigr)^2\le(2n)^2$$ For (\ref{closure_powers_2}), $a^x+a^y\le a^{x+y}$ for real $x,y\ge1$ and $a\ge2$, because $a^{-x}+a^{-y}\le1$. By induction, $\sum_sa^{n_s}\le a^n$.
\end{proof}

We are now ready to prove Theorem~\ref{thm:intro_main}. For convenience, we restate it.

\begin{thm}\label{thm:main}
Let $\phi$ be a linearly growing UPG automorphism of a finitely generated free group $F$, let $G=F\rtimes_\phi\mathbb Z$, and let $H\le G$ be a finitely generated subgroup which is not conjugate into a vertex group. Then
\[
\Delta^G_H(n)\ \simeq\ \begin{cases}
n & \text{if every horizontal component consists of one piece},\\
n^2 & \text{if some horizontal component does not consist of one}\\
 & \text{piece, and every $\omega_i$ is trivial},\\
2^n & \text{if some $\omega_i$ is nontrivial}.
\end{cases}
\]
In particular, $H$ is undistorted in $G$ if and only if every horizontal component of $H$ consists of one piece.
\end{thm}

\begin{proof}
The three cases are exhaustive and mutually exclusive. Indeed, if $\omega_i$ is nontrivial, then $H_i$ does not fix a vertex of $T_i$ (see Proposition~\ref{prop:dilation_basis}(\ref{prop:dilation_basis_3})), so $T_i$ does not consist of one piece by Theorem~\ref{thm:dichotomy} and Lemma~\ref{lem:one_piece}.

\textit{Lower bounds.} If some $T_i$ does not consist of one piece, then $T_i$ is unbounded by Lemma~\ref{lem:one_piece}, so by Theorem~\ref{thm:dichotomy} some element of $H_i$ acts loxodromically on $T_i$, and $\Delta^G_{H_i}(n)\succeq n^2$, and $\Delta^G_{H_i}(n)\succeq2^n$ if $\omega_i\neq1$. Proposition~\ref{prop:distortion-Hi-H} gives the same lower bounds for $\Delta^G_H$. In all cases $\Delta^G_H(n)\succeq n$, because $H$ is infinite.

\textit{Upper bounds.} Let $f$ be the function of Theorem~\ref{thm:reduction}. If every horizontal component consists of one piece, then every $H_i$ is undistorted in $G$ by Lemma~\ref{lem:one_piece} and Theorem~\ref{thm:dichotomy}, so $f(n)\le Cn+C$, hence $\overline f(n)\le 2Cn$ for $n\ge1$, and Theorem~\ref{thm:reduction} gives $\Delta^G_H(n)\preceq n$. If all $\omega_i$ are trivial, Theorem~\ref{thm:upper_Hi}(\ref{upper_Hi_1}) and the definition of $\preceq$ give $f(n)\le C(n+1)^2$ for a constant $C$, so $\overline f(n)\le4Cn^2$ by Lemma~\ref{lem:closure_powers}(\ref{closure_powers_1}), and Theorem~\ref{thm:reduction} gives $\Delta^G_H(n)\preceq n^2$. In general, Theorem~\ref{thm:upper_Hi}(\ref{upper_Hi_2}) and the definition of $\preceq$ give a constant $C'\ge1$ with $f(n)\le C'a^n$ for all $n\ge0$, where $a:=2^{C'}\ge2$. So $\overline f(n)\le C'a^n$ by Lemma~\ref{lem:closure_powers}(\ref{closure_powers_2}), and Theorem~\ref{thm:reduction} gives $\Delta^G_H(n)\preceq a^n\simeq2^n$.
\end{proof}

\begin{cor}\label{cor:spectrum}
Let $G$ be a linear free-by-cyclic group. The distortion function of every finitely generated subgroup of $G$ is equivalent to $n$, to $n^2$, or to $2^n$.
\end{cor}

\begin{proof}
A finitely generated subgroup which is conjugate into a vertex group is undistorted (see Remark~\ref{rem:whyundistorted}). Otherwise, apply Theorem~\ref{thm:main}.
\end{proof}

\section{Examples}\label{sec:examples}
The distortion functions in this section are known. Our aim is to show that they can be recovered from the horizontal components and the dilation homomorphisms, by computing these invariants and applying Theorem~\ref{thm:main}. In both cases the answer comes from the same mechanism. Trivial slope products give at most quadratic distortion, and a nontrivial dilation gives exponential distortion. Recall that the one-piece condition can be checked on the graph of groups decomposition that $H$ inherits from $G$ (see Lemma~\ref{lem:one_piece}).

\subsection{The fiber and its subgroups}\label{subsec:ex_fiber}
Let $G = F \rtimes_{\phi} \mathbb Z$ with $\phi$ polynomially growing of degree $d \ge 0$. Then $\Delta_{F}^{G}(n) \simeq n^{d+1}$ by combining \cite{Ger94} and \cite{Mac02}. The argument is as follows. Gersten \cite[Theorem~4.1, Proposition~4.2]{Ger94} shows that $\mathrm{Div}(G) \preceq \Delta_{F}^{G}(n) \preceq n^{d+1}$, where $\mathrm{Div}(G)$ is the divergence of $G$. Macura \cite{Mac02} shows that $n^{d+1} \preceq \mathrm{Div}(G)$, hence $\Delta^{G}_{F}(n) \simeq n^{d+1}$. 

\begin{exmp}\label{ex:fiber_subgroups}
 Let $\phi$ be a linearly growing UPG automorphism of a finitely generated free group $F$, and let $G=F\rtimes_\phi\mathbb Z$.

Suppose that $H\le F$ is finitely generated and not conjugate into a vertex group (the remaining case is covered by Remark~\ref{rem:whyundistorted}), and let $q\colon G\to\mathbb Z$ be the homomorphism with kernel $F$ (see Section~\ref{subsec:free-by-cyclic}). For a black vertex $b$ of $\Delta$ we have $G_b=\langle r_b\rangle\times\langle t_b\rangle$ with $r_b\in F$ and $t_b\in Ft$ (see Proposition~\ref{prop:split-over-tori}(\ref{prop:split-over-tori2})), so $q(r_b^Xt_b^Y)=Y$. Since $F$ is normal in $G$, every element of $F\cap G_{\tilde b}$ has coordinates of the form $(X,0)$. Similarly $F\cap\langle t_{\tilde w}\rangle=1$ for every white vertex $\tilde w$, since $q(t_{\tilde w})=1$. Consequently:
\begin{enumerate}
\item every white vertex of $T_H$ is of fiber type, and every black vertex of $T_H$ is of type $(H)$ or (0) (see Lemma~\ref{lem:black_types}),
\item if $\tilde b$ is a horizontal black vertex and $u$ generates $H_{\tilde b}$, then $u$ has coordinates $(a,0)$ with $a\neq0$, so $\chi_{\tilde e}(u)=a$ for every edge $\tilde e$ at $\tilde b$, and every slope at $\tilde b$ equals $1$,
\item every slope product equals $1$, and every dilation homomorphism $\omega_i$ is trivial.
\end{enumerate}
In particular, Theorem~\ref{thm:lower_bound} only gives quadratic lower bounds for subgroups of the fiber.
\end{exmp}

Since every dilation homomorphism of $H$ is trivial, the following corollary follows from Theorem~\ref{thm:main}.
\begin{cor}\label{cor:intro_fiber}
Let $G$ be as in Theorem~\ref{thm:intro_main}, and let $H\le F$ be a finitely generated subgroup of the fiber which is not conjugate into a vertex group. Then $H$ is undistorted in $G$ if every horizontal component of $H$ consists of one piece, and $\Delta^G_H(n)\simeq n^2$ otherwise.
\end{cor}

\begin{exmp}[The fiber]\label{ex:fiber}
Let $H=F$ be the fiber subgroup, and assume that the induced splitting of $F$ is nontrivial. Every white vertex is of fiber type, since $F$ meets each centre $\langle t_{\tilde w}\rangle$ trivially. For every black vertex $\tilde b$, the intersection $F\cap G_{\tilde b}$ is infinite cyclic. Hence every edge of $T_H$ is horizontal. Since no black vertex of $\Delta$ is a leaf, every black vertex has at least two horizontal edges, so $\widehat Y=T_H$. This horizontal component is unbounded, so it does not consist of one piece (see Lemma~\ref{lem:one_piece}). By Corollary~\ref{cor:intro_fiber}, we have $\Delta^G_F(n)\simeq n^2$. This is the case $d=1$ of the theorem that $\Delta^G_F(n)\simeq n^{d+1}$ when $\phi$ is polynomially growing of degree $d$ \cite{Ger94,Mac02}. 
\end{exmp}

\subsection{The hydra subgroup}\label{subsec:ex_hydra}
The following example is due to Dison and Riley \cite{DR13}. It shows how the dilation detects exponential distortion.

\begin{figure}[htb]
\centering
\begin{tikzpicture}[scale=1]
\begin{scope}
\draw[nhor] (0,0) .. controls (1.2,0.9) .. (2.4,0)
  node[pos=0.5,above=2pt,font=\scriptsize] {$e_1\colon a\mapsto a,\ t\mapsto t$};
\draw[nhor] (0,0) .. controls (1.2,-0.9) .. (2.4,0)
  node[pos=0.5,below=2pt,font=\scriptsize] {$e_2\colon a\mapsto c,\ t\mapsto c^{-1}t$};
\node[wfib,minimum size=9pt] at (0,0) {}; \node[blk] at (2.4,0) {};
\node[font=\small,left=4pt] at (0,0) {$V$};
\node[font=\small,right=3pt] at (2.4,0) {$E$};
\node[font=\small] at (1.2,1.5) {$G$};
\end{scope}
\begin{scope}[shift={(6,0)}]
\draw[hor] (0,0) .. controls (1.2,0.9) .. (2.4,0)
  node[pos=0.5,above=2pt,font=\scriptsize] {$u\mapsto u$};
\draw[hor] (0,0) .. controls (1.2,-0.9) .. (2.4,0)
  node[pos=0.5,below=2pt,font=\scriptsize] {$u\mapsto v=c^2t^{-1}$};
\node[wfib,minimum size=9pt] at (0,0) {}; \node[blk] at (2.4,0) {};
\node[font=\small,left=4pt] at (0,0) {$Q$};
\node[font=\small,right=3pt] at (2.4,0) {$\langle u\rangle$};
\node[font=\small] at (1.2,1.5) {$H$};
\end{scope}
\end{tikzpicture}
\caption{The graphs of groups of $G$ (left) and of the hydra subgroup $H$ (right) in Example~\ref{ex:hydra}. Both edges of $H$ are horizontal.}
\label{fig:hydra}
\end{figure}
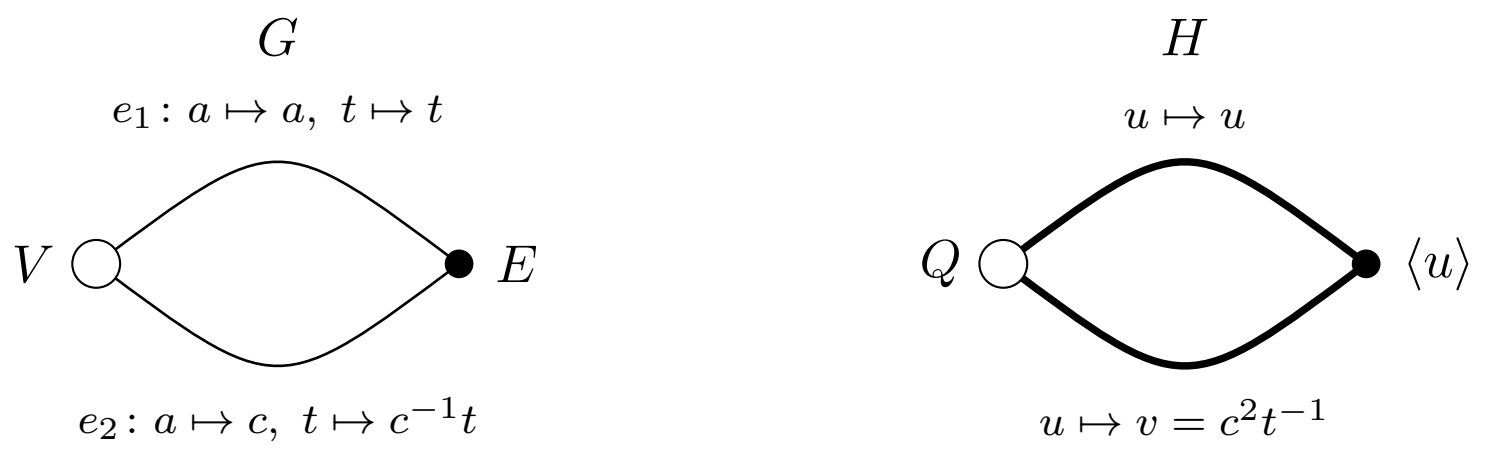

\begin{exmp}[The hydra subgroup]\label{ex:hydra}
Let $F=F(a,b)$, $\phi(a)=a$, $\phi(b)=ba$, and
\[
G=F\rtimes_\phi\langle t\rangle=\langle a,b,t\mid tat^{-1}=a,\ tbt^{-1}=ba\rangle .
\]
Since $\phi^n(b)=ba^n$, the automorphism $\phi$ is linearly growing, and it is UPG.
With $s:=t^{-1}$ this is the group $G_2$ of \cite{DR13}, and the subgroup
\[
H:=\langle u,y\rangle,\qquad u:=at^{-1},\quad y:=bt^{-1},
\]
is their hydra subgroup $H_2=\langle as,bs\rangle$, which satisfies
$\Delta^G_H(n)\simeq A_2(n)=2^n$ \cite{DR13} (here $A_2$ is the second Ackermann function, see Section~\ref{subsec:beyond}). We recover this from the dilation.

\smallskip
\noindent\emph{A splitting of $G$}:
Let $c:=bab^{-1}$. Then $tct^{-1}=(ba)a(ba)^{-1}=c$, so
\[
V:=\langle a,c\rangle\times\langle t\rangle\cong F_2\times\mathbb Z,
\qquad
E:=\langle a,t\rangle\cong\mathbb Z^2\le V .
\]
A direct computation gives $yay^{-1}=c$ and $yty^{-1}=btb^{-1}=c^{-1}t$.
Hence $G$ is the HNN extension $V*_{\theta}$ of $V$ with stable letter $y$, which conjugates
$E$ onto $\langle c,t\rangle$ via an isomorphism $\theta$ given by $a\mapsto c$, $t\mapsto c^{-1}t$.
(The two presentations are related by $b=yt$: the relation $y(at)y^{-1}=t$
is equivalent to $tbt^{-1}=ba$.)

Subdividing the loop, we obtain the splitting of
Proposition~\ref{prop:split-over-tori}, that is,  one white vertex $\bar w$ with
$G_{\bar w}=V$, one black vertex $\bar b$ with $G_{\bar b}=E$, $r_b=a$,
$t_b=t$, and two edges $e_1,e_2$, both oriented from $\bar b$ to $\bar w$, with
$G_{e_1}=G_{e_2}=E$ (the edge monomorphisms into $G_{\bar b}$ are the identity)
and edge monomorphisms $\iota_{e_i}\colon E\to V$ into the white vertex group given by
\[
\iota_{e_1}\colon\ a\mapsto a,\ t\mapsto t,\qquad
\iota_{e_2}\colon\ a\mapsto c,\ t\mapsto c^{-1}t .
\]
Thus $p_{e_1}=a$, $k(e_1)=0$ and $p_{e_2}=c$, $k(e_2)=-1$.

This splitting has the properties of Proposition~\ref{prop:split-over-tori}. Indeed, $\Delta$ is bipartite, both edges are oriented from $\bar b$ to $\bar w$, and the edge monomorphisms into $G_{\bar b}$ are the identity. The vertex groups are $G_{\bar w}=\langle a,c\rangle\times\langle t\rangle$ and $G_{\bar b}=\langle a\rangle\times\langle t\rangle$ with $t\in Ft$. The subgroup $\langle a\rangle$ is maximal cyclic in $F$, and $p_{e_1}=a$ and $p_{e_2}=c$ generate maximal cyclic subgroups of the free group $\langle a,c\rangle$, because $a,c$ is a basis of it. The centre $\langle t\rangle$ of $V$ lies in $\iota_{e_1}(E)=\langle a,t\rangle$ and in $\iota_{e_2}(E)=\langle c,t\rangle$. Moreover $k(e_1)\neq k(e_2)$, and $\langle t,\,ta\rangle=E$. Finally, two distinct black neighbours of $\bar w$ in the Bass--Serre tree have stabilizers $g\,\iota_{e_i}(E)\,g^{-1}$ and $g'\iota_{e_{i'}}(E)\,g'^{-1}$ with $g,g'\in V$. Their images in $\langle a,c\rangle$ are distinct maximal cyclic subgroups, since $\langle a\rangle$ and $\langle c\rangle$ are malnormal and $a$ is not conjugate to $c^{\pm1}$. Distinct maximal cyclic subgroups of a free group intersect trivially, so the two stabilizers intersect in $\langle t\rangle$, which is cyclic. By equivariance, the same holds at every white vertex. A topological representative as in Section~\ref{subsection:topological_representative} is obtained by taking for $X_{\bar w}$ a rose with two petals $a$ and $c$, and gluing the torus of $\bar b$ along $e_1$ and $e_2$ to these petals with the shears $k(e_1)=0$ and $k(e_2)=-1$.

In the Bass--Serre tree $T$, the black vertex $\bar b$ has exactly two edges,
$\tilde e_1=(\bar b,\bar w)$ over $e_1$ and $\tilde e_2=(\bar b,y^{-1}\bar w)$ over $e_2$.
The vertical elements of $\bar w$ and $y^{-1}\bar w$ are $t$ and $y^{-1}ty=at$,
with coordinates $(0,1)$ and $(1,1)$ (see Lemma~\ref{lem:slope_vertical}).

\smallskip
\noindent\emph{A splitting of the subgroup $H$.}
Let $v:=yuy^{-1}=c\,(c^{-1}t)^{-1}=c^2t^{-1}$ and $Q:=\langle u,v\rangle\le V$. Since $H=\langle u,y\rangle$ and $v=yuy^{-1}$, we have $H=\langle Q,y\rangle$. Let $\widehat H:=\langle Q,\hat y\mid \hat y u\hat y^{-1}=v\rangle$ be the abstract
HNN extension, and let $\widehat H\to H$ be the natural surjection.
The projection $V\to\langle a,c\rangle$, which maps $u\mapsto a$ and $v\mapsto c^2$, restricts to an isomorphism from $Q$ to $\langle a, c^2 \rangle$, which is free with basis $a, c^2$,
and therefore
$
Q\cap\langle t\rangle=1,\quad Q\cap E=\langle u\rangle,\quad
Q\cap\langle c,t\rangle=\langle v\rangle .
$
Using Britton's lemma \emph{in $G$}, which is an HNN extension
of $V$ with stable letter $y$, the surjective homomorphism $\widehat H\to H$ is injective, so 
 $H$ is the HNN extension of $Q$ with stable letter $y$
conjugating $\langle u\rangle$ onto $\langle v\rangle$, and $H\cap V=Q$ (see
Figure~\ref{fig:hydra}). Consequently, every white vertex of $T_H$ is of fiber
type, every black vertex of $T_H$ has stabilizer conjugate to
$\langle u\rangle\cong\mathbb Z$ and two horizontal edges, and
$\widehat Y=T_H$ is a single horizontal component. It is unbounded, since $y$
is loxodromic, so it does not consist of one piece.

\smallskip
\noindent\emph{Computing the dilation:}
The generator $u=at^{-1}$ of $H_{\bar b}$ has coordinates $(1,-1)$, so
\[
\chi_{\tilde e_1}(u)=1,\qquad \chi_{\tilde e_2}(u)=1+(-1)(-1)=2,\qquad
\mathfrak{sl}(\tilde e_1\to\tilde e_2)=2 .
\]

The geodesic from $\bar w$ to $y^{-1}\bar w$ passes through $\bar b$ only,
entering through $\tilde e_1$ and leaving through $\tilde e_2$, hence
\[
\omega(y^{-1})=\mathfrak S(\bar w,y^{-1}\bar w)=2\neq1 .
\]
Theorem~\ref{thm:main} gives $\Delta^G_H(n)\simeq 2^n$, in agreement with
\cite{DR13}.
\end{exmp}

\bibliographystyle{alpha}
\bibliography{distortion}

\end{document}